\documentclass[11pt,reqno]{amsart}
\usepackage[utf8]{inputenc}
\usepackage[margin=1.18in]{geometry}
\usepackage{amssymb,mathtools,graphicx,microtype}
\usepackage[backend=biber,style=alphabetic,sortcites=true,maxbibnames=99,
  isbn=false,url=false]{biblatex}
\DefineBibliographyStrings{english}{inpress={to appear}}
\usepackage[dvipsnames]{xcolor}
\usepackage{tikz}
\usetikzlibrary{arrows.meta}
\usepackage[colorlinks=true,citecolor=MidnightBlue,linkcolor=BrickRed,
  urlcolor=RoyalBlue]{hyperref}
\hypersetup{pdftitle={A new second-moment method for spectral measures of inhomogeneous non-Hermitian random matrices}}
\providecommand{\cal}[1]{\mathcal{#1}}

\theoremstyle{plain}
\newtheorem{theorem}{Theorem}[section]
\newtheorem{lemma}[theorem]{Lemma}
\newtheorem{proposition}[theorem]{Proposition}
\newtheorem{corollary}[theorem]{Corollary}
\theoremstyle{remark}
\newtheorem*{example}{Example}
\newtheorem*{remark}{Remark}
\newtheorem{numberedexample}[theorem]{Example}
\newtheorem{numberedremark}[theorem]{Remark}
\numberwithin{equation}{section}
\numberwithin{figure}{section}

\newcommand{\E}{\mathbb E}
\newcommand{\C}{\mathbb C}
\newcommand{\R}{\mathbb R}
\DeclareMathOperator{\per}{per}
\DeclareMathOperator{\Ent}{Ent}
\DeclareMathOperator{\diag}{diag}
\newcommand{\dd}{\,\mathrm d}
\newcommand{\BL}{\mathrm{BL}}

\title[A second-moment method for spectral measures]
  {A new second-moment method for spectral measures of
   inhomogeneous non-Hermitian random matrices}
\author{Yi Han}
\thanks{Yi Han: Department of Mathematics, Massachusetts Institute of Technology,
Cambridge, MA, USA\newline Email: \texttt{hanyi16@mit.edu}.}
\author{Zheye Song}
\thanks{Zheye Song: Department of Biostatistics,
University of North Carolina at Chapel Hill,
Chapel Hill, NC, USA\newline
Email: \texttt{songzy@unc.edu}.}
\hypersetup{pdfauthor={Yi Han and Zheye Song}}
\date{}

\begin{document}
\begin{abstract}
We develop a second-moment method for the empirical spectral
distribution (ESD) of non-Hermitian random matrices
\(X_N=(\sqrt{s_{ij}}\,\xi_{ij})\) with a deterministic variance profile
\(S_N=(s_{ij})\). When the entries are independent and centered, the second
moment of \(\det(zI-X_N)\) is the permanent \(\per(S_N+|z|^2I)\) of a diagonal
shift of \(S_N\), and the method approximates \(\log|\det(zI-X_N)|\)
by half its logarithm. For i.i.d.\ centered \(\xi_{ij}\) of unit variance
satisfying a Dini-type condition on their Lévy concentration function,
we prove the circular law for every doubly stochastic profile whose
largest entry tends to zero, however slowly. This condition holds for
every law with a bounded density on \(\R\) or \(\C\). We also allow the
two entries \(X_{ij}\) and \(X_{ji}\) of each reciprocal pair to be
correlated, with a nondegenerate joint law that has a finite mixed
fourth moment. If the pairs and the diagonal entries also have bounded
densities on their real linear supports, we prove
an elliptic law for every symmetric doubly stochastic profile whose
largest entry tends to zero. For i.i.d.\ \(\xi_{ij}\) as above and general profiles with vanishing largest entry and bounded average row sum, we show that the ESD is asymptotically close to the Brown measure of the matrix
\((\sqrt{s_{ij}}\,c_{ij})\) of free circular elements. The ESD
converges weakly in probability if and only if
\(N^{-1}\log\per(S_N+qI)\) converges for every \(q>0\), and these
limits determine the radial distribution of the limit. As
applications, we give sufficient conditions for convergence in terms
of tree weights, characterize the radial laws that arise from profiles
with uniformly bounded row and column sums, and show that the empirical distribution of the zeros of \(\per(zI-X_N)\) is asymptotically close to the ESD.
\end{abstract}
\maketitle

\section{Introduction}\label{sec:introduction}

Let \(\xi_{ij}\) be i.i.d.\ centered random variables with unit
variance. The classical circular law states that the empirical
spectral distribution (ESD) of \(N^{-1/2}(\xi_{ij})_{i,j=1}^N\)
converges to the uniform probability measure on the unit disk. Following Girko's work \cite{cl-girko}, Bai gave a proof under
additional moment and smoothness assumptions \cite{cl-bai}. Tao and Vu,
with an appendix by Krishnapur, proved it under the assumption of a
finite second moment \cite{cl-taovu}; see \cite{cl-bc} for a survey.
Rudelson and Tikhomirov \cite{cl-rt} and Sah, Sahasrabudhe and
Sawhney \cite{cl-sss} extended it to independent Bernoulli
sparsifications, provided the expected number of selected entries per
row diverges. Litvak et al.\ treated sparse random regular
digraphs \cite{cl-regular}.

In this paper we study inhomogeneous random matrices
\(X_N=(\sqrt{s_{ij}}\,\xi_{ij})_{i,j=1}^N\) with a deterministic,
possibly sparse, variance profile \(S_N=(s_{ij})\). A central
difficulty in Girko's Hermitization approach
\cite{cl-girko,cl-bai,cl-bc} is to control the small singular values
of \(X_N-zI\), particularly on sparse deterministic supports.

In the dense regime \(s_{ij}\le C/N\), Cook, Hachem, Najim and Renfrew
constructed deterministic equivalents of the ESD through master
equations. Under moment and admissibility assumptions, they proved that
the ESD is asymptotically close to these equivalents in probability,
and for doubly stochastic profiles they obtained the circular law
\cite{vp-chnr1,vp-chnr}. For variances uniformly comparable
to \(N^{-1}\), Alt, Erd{\H{o}}s and Kr{\"u}ger established a local
inhomogeneous circular law at the optimal bulk scale under moment and
density assumptions \cite{vp-aek}.

When the largest variance decays at a sufficiently fast polynomial
rate in \(N\), Han proved circular laws for sparse deterministic
profiles, including arbitrary doubly stochastic profiles with Gaussian entries or with subgaussian entries satisfying additional
density hypotheses \cite[Theorem~1.5]{vp-han}. For symmetric or transitively invariant profiles with real or circular
complex Gaussian entries, he later showed that bandwidths
\(N^{2/3+c}\) with fixed \(0<c<1/3\) suffice \cite[Theorem~1.6]{cl-han-coarse}, improving the bandwidth \(N^{5/6+\epsilon}\) of \cite[Theorem~1.5]{vp-han}. Under bounded-density and moment
assumptions, Han's transfer arguments for open block-tridiagonal and
periodic band models evaluate the full shifted determinant along the
band and reach bandwidths \(W_N\) that diverge arbitrarily slowly
\cite{intro-han-block,cl-han-periodic}. An additional bandwidth assumption lets the periodic full-block
results cover real subgaussian discrete laws such as centered
Bernoulli entries \cite{cl-han-periodic}.

Two questions from Tikhomirov's 2022 ICM survey motivate this paper.
Problem~9 asks whether every doubly stochastic profile whose largest
standard deviation decays polynomially yields the circular law
\cite[Section~5, Problem~9]{cl-tikhomirov-icm}. Problem~8 asks which
profiles yield the circular law for every centered unit-variance
entry law.

We propose a new approach that compares the full shifted
determinant with its second moment. For independent entries, this
second moment is the permanent \(\per(S_N+|z|^2I)\), and we
approximate \(\log|\det(zI-X_N)|\) by half its logarithm. We show that the error in this logarithmic comparison is \(o(N)\) in
\(L^1\) and evaluate the permanent asymptotically. This leads to the
following results.
\begin{enumerate}
\item For entry laws satisfying a Dini-type condition on their Lévy
concentration function, we prove
the circular law for every doubly stochastic profile whose largest
variance tends to zero, at any rate (Theorem~\ref{vp:circle}). For
these entry laws this covers Problem~9.
\item Allowing the two entries \(X_{ij}\) and \(X_{ji}\) of each
reciprocal pair to be correlated, we prove an elliptic law for
symmetric doubly stochastic profiles under assumptions on the pair law
(Theorem~\ref{ep:main}).
\item For independent entries and general profiles, we compare the
ESD with a Brown measure built at each dimension from the variance profile and a deterministic diagonal matrix (Theorem~\ref{vp:brown}). When this matrix vanishes, the shifted permanent gives a
criterion for the convergence of the ESD and for the circular law
(Theorem~\ref{vp:pressure-main} and
Corollary~\ref{vp:circular-pressure}). This addresses Problem~8 of
Tikhomirov's survey \cite[Section~5]{cl-tikhomirov-icm} for entry laws satisfying the Dini-type condition above.
\end{enumerate}
Since the permanental polynomial \(\per(zI-X_N)\) has the same second
moment as the determinant, the same comparison applies to its zeros
(Theorem~\ref{vp:permanent}).

\subsection{Main results}

For an \(N\times N\) matrix \(Y_N\), let
\(\mu_N=N^{-1}\sum_{j=1}^N\delta_{\lambda_j(Y_N)}\) denote its ESD,
where eigenvalues are counted with algebraic multiplicity. Write
\(\epsilon_N=\max_{i,j}s_{ij}\) for the largest entry variance. We
suppress the dependence of the entries of \(S_N\) on \(N\).

\subsubsection{Independent entries and the circular law}

For a complex random variable \(Z\), let
\(Q_Z(t)=\sup_{w\in\C}\Pr(|Z-w|\le t)\) be its Lévy concentration
function, and set \(\Psi_Z(t)=\int_0^t Q_Z(r)\,\mathrm dr/r\).
Our results for independent entries assume a fixed law \(\xi\) with
\begin{equation}\label{pd:law}
 \E\xi=0,\qquad \E|\xi|^2=1,\qquad \Psi_\xi(1)<\infty .
\end{equation}
We call the last hypothesis the planar Dini condition. It is implied
by any power bound \(Q_\xi(t)\le Bt^\beta\) with \(\beta>0\), and hence holds for laws with bounded density on the
real line or in the plane. It also holds for some singular continuous
laws; for instance, the middle-thirds Cantor distribution, after
centering and rescaling, satisfies the power bound with
\(\beta=\log2/\log3\). The real and imaginary parts of
\(\xi\) may be dependent.

Let \(\mu_{\rm circ}\) be the uniform probability measure on the unit
disk. For the almost sure results on independent entries, we assume
the power bound
\begin{equation}\label{pd:power}
 Q_\xi(t)\le Bt^\beta\quad(0<t\le1)
 \qquad\text{for some }B,\beta>0.
\end{equation}

\begin{theorem}[Circular law for doubly stochastic profiles]\label{vp:circle}
Let \(\xi\) satisfy \eqref{pd:law}, and let
\(X_{ij}=\sqrt{s_{ij}}\,\xi_{ij}\), where the \(\xi_{ij}\) are
independent copies of \(\xi\). Suppose that each \(S_N\) is doubly
stochastic and that \(\epsilon_N\to0\). Then:
\begin{enumerate}
\item The ESD \(\mu_N\) of \(X_N\) satisfies \(\mu_N\Rightarrow\mu_{\rm circ}\) in probability.
\item Under \eqref{pd:power}, the convergence is almost sure under
any coupling across dimensions.
\end{enumerate}
\end{theorem}

For a uniform periodic band of width \(W_N\), the condition
\(\epsilon_N\to0\) means \(W_N\to\infty\). Using the band geometry, Han's transfer argument reaches this regime
for bounded-density entries under an additional moment assumption
\cite[Theorem~2.2]{cl-han-periodic}. We do not require such geometry for
Theorem~\ref{vp:circle}, since the second moment gives the scale of the determinant for every
doubly stochastic support. The following support is far from any
band.

\begin{example}[A profile from finite geometry]
Let \(h\to\infty\) through prime powers, and consider the points and
lines of the projective plane \(\mathrm{PG}(2,h)\)
\cite{intro-projective-geometry}; there are \(N=h^2+h+1\) of each.
Enumerate both sets and define
\[
 s_{p\ell}=\frac{\mathbf1_{\{p\in\ell\}}}{h+1}.
\]
Each point lies on \(h+1\) lines, and each line contains \(h+1\)
points. Thus \(S_N\) is doubly stochastic with
\(\epsilon_N=(h+1)^{-1}\asymp N^{-1/2}\). By
Theorem~\ref{vp:circle}, the circular law holds on this deterministic
support for every entry law satisfying \eqref{pd:law}.

Distinct rows correspond to distinct points, and two distinct points
lie on exactly one common line. Hence any two rows share exactly one
nonzero column, and every \(2\times2\) submatrix has at least one zero
entry.

One can also check that any fixed positive fraction of the rows meets
all but \(o(N)\) columns.
Hence the support contains no all-zero rectangle whose two sides are
both of order \(N\). By contrast, a band of width \(o(N)\), cyclic or
not, contains such rectangles between distant sets of rows and
columns. Row and column permutations preserve both side lengths of an
all-zero rectangle, so this profile cannot be reordered into a narrow
band.
\end{example}

\begin{remark}[Independent non-identical entries]
Let \(\xi_{ij}^{(N)}\) be independent centered unit-variance
variables, not necessarily identically distributed, such that
\begin{align*}
 \lim_{K\to\infty}\sup_{N,i,j}
 \E\bigl[|\xi_{ij}^{(N)}|^2
       \mathbf1_{\{|\xi_{ij}^{(N)}|>K\}}\bigr]&=0,\\
 \lim_{\delta\downarrow0}\sup_{N,i,j}
 \Psi_{\xi_{ij}^{(N)}}(\delta)&=0.
\end{align*}
Then Theorems~\ref{vp:circle}, \ref{vp:brown}, \ref{vp:pressure-main}
and~\ref{vp:permanent}, together with
Corollary~\ref{vp:circular-pressure}, remain valid with
\(\xi_{ij}^{(N)}\) in place of \(\xi_{ij}\).
The first condition makes the mean-square approximation in
Lemma~\ref{pd:approx} uniform over the entries. In the conditioning
step of Lemma~\ref{pd:row}, we choose an entry whose coefficient has
maximal modulus and keep its index fixed as the radius varies. The
integrated bound then involves the Dini integral of a single entry
law, to which the second condition applies.
For the almost sure conclusions, we additionally assume
\eqref{pd:power} for every \(\xi_{ij}^{(N)}\), with the same
constants \(B,\beta>0\).
\end{remark}

\subsubsection{Reciprocal pairs and the elliptic law}

We next keep the profile doubly stochastic, but allow the two entries
\(X_{ij}\) and \(X_{ji}\) of each reciprocal pair to be dependent.
Let \(S_N\) be a symmetric profile, and set
\[
 (X_{ij},X_{ji})=\sqrt{s_{ij}}\,(U_{ij},V_{ij})\quad(i<j),\qquad
 X_{ii}=\sqrt{s_{ii}}\,\zeta_i .
\]
Here the \((U_{ij},V_{ij})\) are copies of a centered pair \((U,V)\)
with \(\E|U|^2=\E|V|^2=1\), and the \(\zeta_i\) are copies of a
centered unit-variance variable \(\zeta\). The pairs and the diagonal
variables are mutually independent. Write
\(\tau=\E(UV)\) and \(\kappa=\E(U\overline V)\). We assume the
following.
\begin{enumerate}
\item The pair \((U,V)\) and the variable \(\zeta\) have bounded
densities with respect to Lebesgue measure on fixed real linear
subspaces \(L\subset\C^2\) and \(L_0\subset\C\), respectively.
\item The pair law satisfies
\[
 \E|UV|^2<\infty,\qquad |\tau|<1,\qquad |\kappa|<1 .
\]
\end{enumerate}
The first assumption allows, for instance, real pairs, complex pairs
and Gaussian pairs supported on a proper real subspace.

The limiting ellipse depends only on \(\tau\). Write \(\rho=|\tau|\),
and let \(\theta=\arg(\tau)/2\) if \(\tau\ne0\) and \(\theta=0\)
otherwise. Let \(\mu_\tau\) be the uniform probability measure on
\[
 {\cal E}_\tau=e^{i\theta}
 \left\{x+iy:\frac{x^2}{(1+\rho)^2}+
                    \frac{y^2}{(1-\rho)^2}<1\right\}.
\]
When \(\tau=0\), this is the unit disk and \(\mu_\tau=\mu_{\rm circ}\).

\begin{theorem}[Elliptic law for reciprocal pairs]\label{ep:main}
Consider the reciprocal-pair model above, and suppose that each
\(S_N\) is symmetric and doubly stochastic, with \(\epsilon_N\to0\).
Then the ESD \(\mu_N\) of \(X_N\) satisfies \(\mu_N\Rightarrow\mu_\tau\) almost surely, under any coupling
across dimensions.
\end{theorem}

\subsubsection{General profiles and Brown measures}

We now return to independent entries and general variance profiles.
Without double stochasticity, the limit, if it exists, depends on the
profile, so we compare the ESD at each dimension with a deterministic
measure \(\nu_N\). This separates the random-matrix estimate from the
deterministic question of whether \((\nu_N)\) converges.

\paragraph{\textit{The comparison measure.}}
The measure \(\nu_N\) is the Brown measure of a free analogue of the
random matrix. Let \(\varphi\) be the normalized trace of a tracial von Neumann
algebra, and let \(c_{ij}\) be freely independent standard circular
elements in it, so that \(\varphi(c_{ij}^*c_{ij})=1\)
\cite{fp-vdn,fp-ns}. Replacing each \(\xi_{ij}\) by \(c_{ij}\) gives
\[
 C_{S_N}=\sum_{i,j}\sqrt{s_{ij}}E_{ij}\otimes c_{ij},
\]
where the \(E_{ij}\) are the matrix units. We equip the matrix algebra
with the normalized trace
\(\mathfrak t_N=N^{-1}\operatorname{Tr}\otimes\varphi\). We also allow a deterministic
complex diagonal matrix \(A_N=\diag(a_i)\), which we call the
diagonal field. Let \(\nu_N\) be the Brown measure of
\(A_N+C_{S_N}\) with respect to \(\mathfrak t_N\)
\cite{vp-brown-original,vp-hs}; its logarithmic potential is the Brown potential \(U_N\) defined in Section~\ref{sec:variance-profiles}. Both
the profile and the diagonal field may depend on \(N\). The next
theorem compares \(\nu_N\) with the ESD of \(A_N+X_N\).

\begin{theorem}[Brown-measure comparison]\label{vp:brown}
Let \(\xi\) satisfy \eqref{pd:law}, and let
\(X_{ij}=\sqrt{s_{ij}}\,\xi_{ij}\), where the \(\xi_{ij}\) are
independent copies of \(\xi\). Let \(\mu_N\) be the ESD of
\(A_N+X_N\). Suppose that \(\epsilon_N\to0\) and
\begin{equation}\label{vp:average-bounds}
 \sup_N\frac1N\sum_{i,j}s_{ij}<\infty,\qquad
 \sup_N\frac1N\sum_i|a_i|^2<\infty.
\end{equation}
Then:
\begin{enumerate}
\item For every bounded continuous \(f:\C\to\R\),
\(\int f\dd\mu_N-\int f\dd\nu_N\to0\) in probability.
\item Under \eqref{pd:power}, this convergence holds almost surely
for all such \(f\) simultaneously under any coupling across dimensions.
\end{enumerate}
Whenever \(\nu_N\Rightarrow\nu\), the same conclusions give
\(\mu_N\Rightarrow\nu\).
\end{theorem}

Apart from \(\epsilon_N\to0\) and \eqref{vp:average-bounds}, the
profile is arbitrary. It may be nonsymmetric or reducible and may have
zero rows. Individual row sums may grow with \(N\), and \(\epsilon_N\)
may decay much more slowly than \(1/N\). The variance bound can also be weakened as follows.

\begin{remark}[A weaker variance condition]
Theorems~\ref{vp:brown} and~\ref{vp:pressure-main}, together with
Corollary~\ref{vp:circular-pressure}, also hold with the variance
bound in \eqref{vp:average-bounds} replaced by
\[
 \lim_{r\to\infty}\limsup_{N\to\infty}
 \frac1N\sum_i\log(1+d_i/r^2)=0,\qquad d_i=\sum_j s_{ij}.
\]
This is the row-tail condition \eqref{vp:row-log-tail}.
\end{remark}

\paragraph{\textit{Deciding whether the ESD converges.}}
By Theorem~\ref{vp:brown}, it remains to decide whether \(\nu_N\)
converges. When \(A_N=0\), \(\nu_N\) is rotation invariant, and in
Section~\ref{sec:variance-profiles} we compare its logarithmic
potential at \(z\) with the normalized log permanent of
\(S_N+|z|^2I\). Varying \(|z|\) gives the following criterion.

\begin{samepage}
\begin{theorem}[Permanent pressure criterion]\label{vp:pressure-main}
Assume the hypotheses of Theorem~\ref{vp:brown} with \(A_N=0\), and
define the scalar pressure
\begin{equation}\label{vp:scalar-pressure}
 H_N(t)=\frac1N\log\per(S_N+e^t I),\qquad t\in\R .
\end{equation}
Then the ESD converges weakly in probability if and only if
\(H_N(t)\) converges for every \(t\in\R\). In this case, let
\(H(t)=\lim_{N\to\infty}H_N(t)\). The limiting law \(\nu\) is
rotation invariant, with radial distribution
\[
 \nu\{w:|w|\le r\}=H'_+(2\log r),\qquad r>0 ,
\]
where \(H'_+\) denotes the right derivative of \(H\).
Under \eqref{pd:power}, the same criterion characterizes almost sure
weak convergence under any coupling across dimensions.
\end{theorem}
\end{samepage}
\noindent As shown in Section~\ref{sec:variance-profiles}, each
\(H_N\) is convex. Hence the limiting pressure \(H\) is a real-valued
convex function, and its right derivative
\(H'_+(t)=\lim_{h\downarrow0}h^{-1}(H(t+h)-H(t))\) exists at every
\(t\).
\begin{samepage}
\begin{corollary}[A criterion for the circular law]\label{vp:circular-pressure}
Under the assumptions of Theorem~\ref{vp:brown} with \(A_N=0\),
the ESD converges weakly in probability to \(\mu_{\rm circ}\)
if and only if, for every \(q>0\),
\[
 \lim_{N\to\infty}\frac1N\log\per(S_N+qI)=
 \begin{cases}
 q-1,&0<q\le1,\\
 \log q,&q>1.
 \end{cases}
\]
Under \eqref{pd:power}, the same condition is equivalent to almost
sure convergence under any coupling across dimensions.
\end{corollary}
\end{samepage}

For entry laws satisfying \eqref{pd:law} and profiles satisfying the
assumptions above, Corollary~\ref{vp:circular-pressure} answers
Problem~8 of Tikhomirov's survey \cite[Section~5]{cl-tikhomirov-icm}.
Since these criteria involve only the variance profile, whether a
profile yields the circular law is a deterministic question about the
permanents of \(S_N+qI\).

\par\medskip
\paragraph{\textit{Sparse examples with nonconstant row sums.}}
The tree-weight criterion of Corollary~\ref{vp:tree-main} gives a
concrete sufficient condition for the convergence of the pressures
above. For instance, let \(b\) be a bounded nonnegative measurable
function of period one, and set
\(b_{N,i}=N\int_{(i-1)/N}^{i/N}b(x)\dd x\).
Writing \(d_N(i,j)\) for the distance between \(i\) and \(j\) on the
cycle of \(N\) sites, consider the weighted periodic bands
\[
 s_{ij}=\frac{b_{N,i}b_{N,j}}{2W_N+1}
 \mathbf1_{\{d_N(i,j)\le W_N\}},\qquad
 W_N\to\infty,\quad W_N=o(N).
\]
Their limiting ESD is the average over \(x\in[0,1]\) of the uniform
law on the disk of radius \(b(x)\). When \(b(x)=0\), this law is the
point mass at zero. The function \(b\) need not be continuous and may
vanish on whole regions. The largest
variance is \(O(W_N^{-1})\), so the bandwidth may diverge arbitrarily
slowly. Since \(W_N=o(N)\), the variances may be much larger than
\(N^{-1}\).
Section~\ref{sec:tree-limit} explains how tree weights identify this limit.

\par\medskip
\paragraph{\textit{Prescribing the limiting law.}}
The pressure criterion also gives an inverse result. Fix an
entry law satisfying \eqref{pd:law} and a constant \(d_*>0\). Consider
profiles with \(\epsilon_N\to0\) and row and column sums at most
\(d_*\), and take \(A_N=0\). For a rotation-invariant probability
measure \(\nu\), write \(F(q)=\nu\{|z|^2\le q\}\). By
Theorem~\ref{vp:radial-inverse}, the limits in probability of the ESD
arising from such profiles are exactly the rotation-invariant
probability measures \(\nu\) on \(\{|z|\le\sqrt{d_*}\}\) for which
\(F(q)/q\) is nonincreasing.

To realize a given \(\nu\), we start from a doubly stochastic directed
cycle of layers and change the column weights in a single block so
that they encode the quantiles of \(\nu\). In this way one can
obtain, for example, a limiting law supported on a disk around the
origin and an outer annulus, separated by an open gap (Example~\ref{vp:annular-example} and
Figure~\ref{vp:quantile-cycle}). If the modified and unmodified
random matrices are built from
the same entry variables, they differ by a matrix of rank \(o(N)\),
yet their limiting ESDs can differ.

\subsubsection{Permanental zeros}

Fyodorov conjectured a circular law for the zeros of the permanental
polynomial \(\per(zI-X_N)\) in the complex Ginibre ensemble
\cite[(1.21)]{vp-fyo}. He computed the asymptotics of
\(\E|\per(zI-X_N)|^2\) and proposed comparing the expected logarithm with the
logarithm of the second moment \cite[(1.18) and Section~6]{vp-fyo}.
Feng and Yao proved the circular law for both real and complex
Ginibre matrices \cite[Theorem~1.1]{vp-feng-yao-ginibre}. For the
independent-entry models considered here, the permanental polynomial
has the same second moment as the determinant, and the logarithmic
comparison applies to both polynomials. The next result therefore
extends the circular law for permanental zeros to real or complex
entry laws satisfying \eqref{pd:law} and to sparse doubly stochastic
variance profiles.

\begin{samepage}
\begin{theorem}[Circular law for permanental zeros]\label{vp:permanent}
Under the assumptions of Theorem~\ref{vp:circle}, let
\(\mu_N^{\mathrm{per}}\) be the normalized zero measure of
\(P_N(z)=\per(zI-X_N)\), where zeros are counted with multiplicity.
Then:
\begin{enumerate}
\item \(\mu_N^{\mathrm{per}}\Rightarrow\mu_{\rm circ}\) in probability.
\item Under \eqref{pd:power}, the convergence is almost sure under
any coupling across dimensions.
\end{enumerate}
\end{theorem}
\end{samepage}

More generally, under the hypotheses of Theorem~\ref{vp:brown} with
\(A_N=0\), the ESD and the normalized zero measure \(\mu_N^{\mathrm{per}}\) approach
the same Brown measure (Proposition~\ref{pp:esd-comparison}). In
particular, the radial laws of Theorem~\ref{vp:radial-inverse} also
arise as limits of \(\mu_N^{\mathrm{per}}\).

\subsection{The second-moment comparison}

We explain the method for the doubly stochastic model of
Theorem~\ref{vp:circle}, whose proof is given in
Section~\ref{sec:planar-dini}. Write
\[
 D_N(z)=\det(zI-X_N),\qquad
 u_N(z)=\frac1N\log|D_N(z)|.
\]
The second moment has an exact deterministic expression
(Lemma~\ref{pd:identity}),
\[
 M_N(z):=\E|D_N(z)|^2=\per(S_N+|z|^2I),
\]
and we call \(p_N(z):=(2N)^{-1}\log M_N(z)\) the second-moment
pressure.
In this notation, the probabilistic part of the proof shows that
\(u_N\) is close to \(p_N\), and the deterministic part shows that
\(p_N\) converges to the logarithmic potential \(U_{\rm circ}\) of the
circular law.

\begin{figure}[!ht]
\centering
\begin{tikzpicture}[
  st/.style={inner sep=3pt,text height=2.1ex,text depth=0.7ex},
  ref/.style={font=\normalfont\footnotesize,inner sep=1pt},
  lab/.style={font=\normalfont\footnotesize,align=left,inner sep=2pt},
  rlab/.style={lab,right=2.5pt},
  head/.style={font=\normalfont\small\itshape},
  arr/.style={-{To[length=1.5mm]},line width=0.4pt}
]
\node[head] at (3.6,0) {Probabilistic part};
\node[head] at (10.6,0) {Deterministic part};
\node[st] (A) at (1.3,-0.95) {\(r_{N,m}\to0\)};
\node[ref,anchor=north] (Ar) at (A.south) {Proposition~\ref{pd:fourth-ds}};
\node[st] (B) at (5.9,-0.95) {concentration of \(u_N^{(m)}\)};
\node[ref,anchor=north] (Br) at (B.south) {Lemma~\ref{pd:row-concentration}};
\node[st] (C) at (3.6,-3.35) {\(\E|u_N^{(m)}-p_N|\to0\)};
\node[st] (D) at (3.6,-5.4) {\(\E|u_N-p_N|\to0\)};
\draw[arr] (Ar.south) -- (C.north west);
\draw[arr] (Br.south) -- (C.north east);
\node[lab,align=center] at (3.6,-2.3) {Paley--Zygmund\\Lemma~\ref{pd:resampling}};
\draw[arr] (C) -- node[rlab] {row replacement\\planar Dini condition\\Corollary~\ref{pd:row-approx}} (D);
\node[st] (E) at (10.6,-0.95) {\(M_N=\per(S_N+|z|^2I)\)};
\node[ref,anchor=north] (Er) at (E.south) {Lemma~\ref{pd:identity}};
\node[st] (F) at (10.6,-5.4) {\(p_N\to U_{\rm circ}\)};
\draw[arr] (Er.south) -- node[rlab] {Gurvits' lower bound\\Gaussian upper bound\\Proposition~\ref{vp:ds}} (F);
\node[st] (G) at (7.1,-7.0) {\(\E\|u_N-U_{\rm circ}\|_{L^1(D_R)}\to0\)};
\node[st] (H) at (7.1,-8.6) {\(\mu_N\Rightarrow\mu_{\rm circ}\) in probability};
\draw[arr] (D.south) -- (G.north west);
\draw[arr] (F.south) -- (G.north east);
\draw[arr] (G) -- node[rlab] {Lemma~\ref{pd:root-return}} (H);
\end{tikzpicture}
\caption{The two parts of the proof of Theorem~\ref{vp:circle}. The
probabilistic part first works at fixed \(m\) with the bounded
approximants \(\eta_m\) and then returns to \(\xi\) by row replacement
as \(m\to\infty\). The deterministic part evaluates
\(p_N=(2N)^{-1}\log M_N\). References below statements point to their
proofs, and arrow labels name the tools and results used.}
\label{intro:second-moment-flow}
\end{figure}
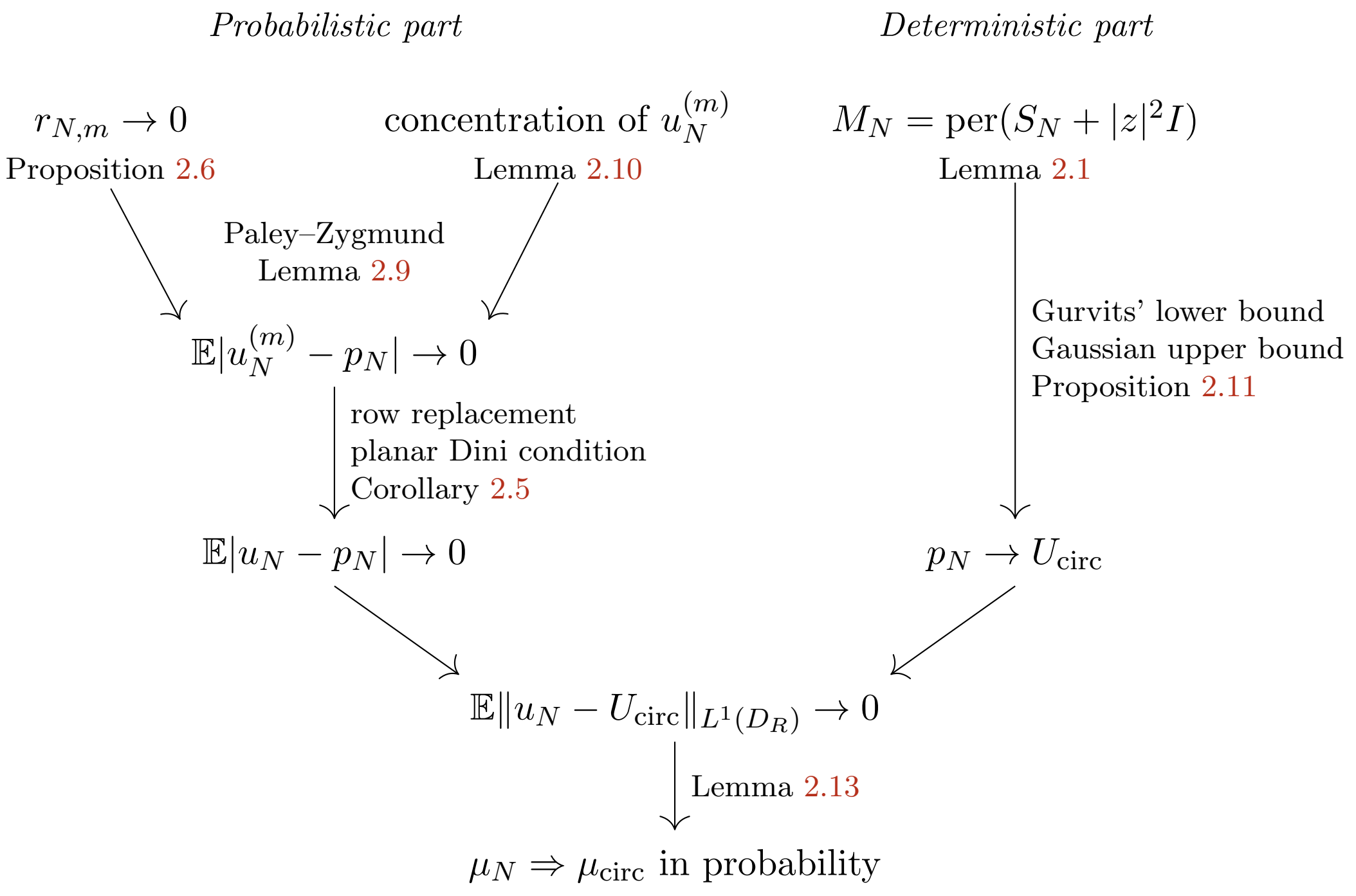

Since \(\xi\) need not have a finite fourth moment, we first work
with the bounded, smoothed approximants \(\eta_m\) of
Lemma~\ref{pd:approx}. Write \(D_N^{(m)}\) for the corresponding
determinant and \(u_N^{(m)}=N^{-1}\log|D_N^{(m)}|\). The second
moment of \(D_N^{(m)}\) is still \(M_N\), because the approximation
preserves centering and variances. Figure~\ref{intro:second-moment-flow}
shows how the estimates below fit together.

For each fixed \(m\), the probabilistic part rests on two estimates,
both uniform in \(z\):
\begin{align*}
 r_{N,m}&:=\sup_{z\in\C}\frac1N\log
  \frac{\E|D_N^{(m)}(z)|^4}{M_N(z)^2}
  \le\frac{C_m}{\log(1/\epsilon_N)}\longrightarrow0,\\
 \Pr\!\left(\left|u_N^{(m)}(z)-\E u_N^{(m)}(z)\right|>t\right)
 &\le2\exp\{-c_mN\min(t^2,t)\},\qquad t>0.
\end{align*}
The first is a bound on the relative fourth moment. We obtain it by
expanding the product of two determinants and two conjugate
determinants. Since the mixed third moments of the approximants vanish, only paired or
quadrupled entries contribute, and their contribution can be counted
by two weighted matchings. Short alternating cycles in these
matchings have small total weight because \(\epsilon_N\to0\), while
long cycles are few because they are disjoint
(Lemma~\ref{pd:matchings}). The second estimate, a concentration
bound for the logarithm, rests on a different feature of the
determinant. Once the other rows are fixed, the determinant is an
affine function of the remaining row. A small-ball bound for this affine form controls how much the logarithm can change when that row
is resampled (Lemma~\ref{pd:row-concentration}).

Together, the fourth-moment bound and the logarithmic concentration
close the gap in Jensen's inequality \(\E u_N^{(m)}\le p_N\).
Indeed, the Paley--Zygmund inequality gives
\[
 \Pr\!\left(u_N^{(m)}(z)\ge
 p_N(z)-\frac{\log2}{2N}\right)
 \ge\frac14e^{-Nr_{N,m}}.
\]
Thus the determinant reaches its second-moment scale with
probability at least \(e^{-o(N)}\). If \(p_N-\E u_N^{(m)}\) stayed
above a positive constant, concentration would bound the probability
of the same event by \(e^{-cN}\), a contradiction.
Lemma~\ref{pd:resampling} makes this comparison quantitative:
\[
 \sup_z\E|u_N^{(m)}(z)-p_N(z)|
 \le C_m\bigl(\sqrt{r_{N,m}}+r_{N,m}+N^{-1/2}\bigr)
 \longrightarrow0.
\]

To return to \(\xi\), we use row replacement
(Proposition~\ref{pd:transport} and Corollary~\ref{pd:row-approx}),
which gives, for all sufficiently large \(m\) and all
\(0<\delta\le1\),
\[
 \sup_{N,z}\E|u_N(z)-u_N^{(m)}(z)|
 \le\frac{\sqrt{\E|\xi-\eta_m|^2}}{\delta}
   +C\bigl\{\Psi_\xi(4\sqrt\delta)+\sqrt\delta\bigr\}.
\]
Here \(C\) is independent of \(N,m,z\). Mean-square approximation
controls the first term, and the planar Dini condition controls the second.
We first let \(N\to\infty\) at fixed \(m\), then let \(m\to\infty\)
at fixed \(\delta\), and finally let \(\delta\downarrow0\).
This proves \(\sup_z\E|u_N(z)-p_N(z)|\to0\). Since \(N\to\infty\) is
taken first at fixed \(m\), the fourth-moment and density bounds of
\(\eta_m\) may grow arbitrarily quickly with \(m\).

It remains to compute the deterministic scale. Proposition~\ref{vp:ds}
uses Gurvits' inequality for the lower bound and a Gaussian matrix
with the same variance profile for the upper bound. Gaussian
concentration reduces the upper bound to a regularized trace
computation. By double stochasticity, the leading tree pairings in
this computation have weight one, and the remaining pairings vanish as
\(\epsilon_N\to0\).
Removing the regularization after letting \(N\to\infty\) gives
\[
 p_N(z)\longrightarrow U_{\rm circ}(z)=
 \begin{cases}
 (|z|^2-1)/2,&|z|\le1,\\
 \log|z|,&|z|>1.
 \end{cases}
\]
This convergence also holds in local \(L^1\). Combining the two
parts yields, on every disk \(D_R=\{z:|z|<R\}\),
\[
 \E\|u_N-U_{\rm circ}\|_{L^1(D_R)}\longrightarrow0.
\]
Lemma~\ref{pd:root-return} then gives the circular law. Under the
power bound \eqref{pd:power}, Lemma~\ref{pd:row-concentration} provides the summable fluctuation estimate
\(\E\|u_N-\E u_N\|_{L^1(D_R)}^4\le C_RN^{-2}\), and the
Borel--Cantelli lemma gives the almost sure conclusion.

\begin{numberedremark}[Necessity of regularity conditions on the entry law]\label{rem:atomic-law}
For each \(0<\alpha<1/3\), the doubly stochastic variance profiles
\(S_N=(s_{ij})\) constructed in Appendix~\ref{app:atomic-profile}
satisfy \(\max_{i,j}s_{ij}<N^{-\alpha}\) for all sufficiently large
\(N\). Let \(\xi_{ij}\) be independent copies of a centered
unit-variance random variable \(\xi\) taking the values \(0,\sqrt{2},-\sqrt{2}\)
with probabilities \(1/2,1/4,1/4\), respectively; its atom at zero
makes \(\Psi_\xi(1)=\infty\). The ESD of
\(X_N=(\sqrt{s_{ij}}\,\xi_{ij})\) then converges weakly in probability
to the point mass at the origin.
Even under the profile assumptions of Tikhomirov's question
\cite[Section~5, Problem~9]{cl-tikhomirov-icm}, centering and unit
variance of the entry law alone do not guarantee the circular law.
With the additional planar Dini condition, Theorem~\ref{vp:circle}
establishes the circular law for every doubly stochastic profile
whose largest variance tends to zero, however slowly.
\end{numberedremark}

\subsection{Related work}

Earlier proofs of the circular law for band matrices control the
logarithmic determinant through the singular values. Jain, Jana,
Luh and O'Rourke combine quantitative convergence of Stieltjes
transforms with a least-singular-value bound for shifted block bands
\cite{cl-jain-band}. Tikhomirov obtains a Gaussian
least-singular-value estimate for general inhomogeneous matrices and
combines it with estimates on the singular spectrum to prove a
circular law for periodic bands \cite[Section~4]{pd-tikhomirov-pseudo}.
Building on this approach, Han uses comparison with a free model to
count small singular values for arbitrary doubly stochastic profiles,
thereby improving the sufficient bandwidths
\cite[Theorem~1.5 and Section~1.4]{vp-han}. Symmetry or transitive
invariance allows a finer Gaussian estimate
\cite[Theorems~1.6 and~2.1]{cl-han-coarse}, but these arguments still
require effective bandwidths above polynomial thresholds.

For narrower bands, Han's arguments instead use the band geometry to
compute the full determinant through transfer products. In the open
block-tridiagonal model, short transfer segments are controlled
uniformly over their boundary frames and combined by a martingale
argument \cite[Sections~1.2 and~2]{intro-han-block}.
Closing a periodic band requires further coefficients, which are
encoded by exterior powers of the transfer product. Their
size is calibrated on a short auxiliary ring, and local estimates
carry this calibration to the target ring
\cite[Sections~3--5 and~10]{cl-han-periodic}.
The band geometry also permits Bernoulli entries, which violate the planar Dini condition. Han's results cover such entries for periodic
scalar strips at polynomial bandwidth
\cite[Theorems~1.11 and~1.12]{cl-han-coarse} and for periodic
full-block bands when \(W_N/\log N\to\infty\)
\cite[Theorem~2.8 and Section~7]{cl-han-periodic}.

These transfer arguments suggest looking for a scale of the
determinant that can be calibrated without band geometry. For
independent entries, the second moment is a natural candidate, since
it is the permanent of a diagonal shift of the variance profile. The identity that expresses it as a permanent has
classical precedents. Edelman used second moments of shifted
determinants to compute the density of the non-real eigenvalues of real Gaussian
matrices \cite{intro-edelman}, and the unshifted identity underlies
the Godsil--Gutman and Barvinok estimators of the permanent
\cite{pd-godsil-gutman-graph,pd-barvinok,intro-rz}. The probabilistic step
in our argument is to show that the second moment captures the typical
logarithm of the determinant. Under the planar Dini condition, the second
moment then supplies the calibration for every doubly stochastic support with
\(\epsilon_N\to0\).

Elliptic laws for dense matrices with correlated reciprocal entries
are established in
\cite{ep-girko,ep-sommers,ep-naumov,ep-no,ep-alt-kruger}. Under our
assumptions on the pair law, Theorem~\ref{ep:main} extends these results to symmetric doubly stochastic profiles whose largest variance may vanish
arbitrarily slowly. Reciprocal dependence changes the second-moment computation, and in the logarithmic comparison we resample the two entries of each reciprocal pair together. Han also studies outliers of
non-Hermitian band matrices and of more general inhomogeneous models
\cite{ep-hanband}.

For general variance profiles, Cook, Hachem, Najim and Renfrew
construct deterministic equivalents of the ESD under a \(C/N\)
variance bound and admissibility assumptions \cite{vp-chnr1,vp-chnr}.
Under assumptions on limiting profiles, Alt and Kr\"uger give related
descriptions through deformed Brown measures
\cite{vp-ak-deformation,vp-ak-brown}. Since we compare both the logarithmic potential of the ESD and the Brown potential with the same second-moment pressure
at each dimension, we need not specify a limiting profile in advance. We also allow variances much larger than \(C/N\), disconnected supports
and zero rows. To identify the common pressure, we use the solution of
the Dyson equation \eqref{vp:dyson}. It supplies the weights in Gurvits' lower bound for the permanent.
Combined with an entropy estimate for weighted permutations, the same
solution also gives the reverse bound (Lemma~\ref{vp:flow} and
Proposition~\ref{vp:entropy}). This step is
related to variational descriptions of Brown potentials, including
Kargin's formula for matrix-valued semicircular elements
\cite[Corollary~1.5]{vp-kargin-variational}.

The radial restriction in Theorem~\ref{vp:radial-inverse} is already
implied by the density formula for the deterministic equivalents of
Cook et al.\ \cite[(5.19)]{vp-chnr}. For a family of products of rectangular complex Gaussian matrices,
Qi and Zhao obtained the uniform law on an annulus after a nonlinear
rescaling of the eigenvalue radii \cite[Theorem~2.1 and
Example~2]{vp-qi-zhao}. Kabluchko and
Zaporozhets studied a related prescription problem for zeros of random
analytic functions, through the choice of coefficient scales
\cite[Theorem~2.9]{vp-kz}. For each fixed entry law satisfying \eqref{pd:law}, every measure
meeting conditions (1) and (2) of Theorem~\ref{vp:radial-inverse}
arises as the limiting ESD of independent-entry matrices with sparse
variance profiles.

\section{Proof of the circular law}\label{sec:planar-dini}

Throughout this section \(S_N=(s_{ij})\) is doubly stochastic,
\(\epsilon_N=\max_{i,j}s_{ij}\to0\), and
\(X_{ij}=\sqrt{s_{ij}}\,\xi_{ij}\), where the \(\xi_{ij}\) are
independent copies of a law satisfying \eqref{pd:law}.
We write \([N]=\{1,\ldots,N\}\). For \(I,J\subset[N]\),
\(S_N[I,J]\) denotes the submatrix with row set \(I\) and column
set \(J\); the same convention applies to other matrices.
Complements are taken in \([N]\). We also write \(D_R=\{z\in\C:|z|<R\}\) and
\(\|f\|_{L^1(D_R)}=\int_{D_R}|f(z)|\dd A(z)\), where
\(\mathrm dA\) denotes planar Lebesgue measure. Let \(\mu_N\) be the
ESD of \(X_N\) and \(D_N(z)=\det(zI-X_N)\).
The logarithmic potential of \(\mu_N\) is
\[
 u_N(z)=\frac1N\log|D_N(z)|.
\]
Write \(M_N(z)=\E|D_N(z)|^2\) for the second moment of the
determinant. We compare \(u_N\) with the deterministic function
\[
 p_N(z)=\frac1{2N}\log M_N(z),
\]
which we call the second-moment pressure. Our goal is to show that
\(u_N\) converges to the logarithmic potential of the circular law,
\[
 U_{\rm circ}(z)=
 \begin{cases} (|z|^2-1)/2,&|z|\le1,\\ \log|z|,&|z|>1.\end{cases}
\]
The proof combines the logarithmic comparison \eqref{pd:ds-gap}
with the pressure limit of Proposition~\ref{vp:ds};
Lemma~\ref{pd:root-return} then gives convergence of the ESDs. Approximation of the entry law is needed only for the
logarithmic comparison.

\begin{lemma}[The second moment as a permanent]\label{pd:identity}
We have
\begin{equation}\label{pd:second-ds}
 M_N(z)=\per(S_N+|z|^2I).
\end{equation}
In particular, \(M_N(z)>0\) for every \(z\).
\end{lemma}
\begin{proof}
Expand the determinant according to the set \(J\) of rows in which a random entry is selected. Distinct entry monomials are orthogonal, since
an entry appearing in just one of them is centered and independent
of the others. Thus
\[
 M_N(z)
 =\sum_{J\subset[N]}|z|^{2(N-|J|)}\per(S_N[J,J])
 =\per(S_N+|z|^2I),
\]
where the empty permanent is one. Since \(S_N\) is doubly
stochastic, \(\per S_N>0\), and the term \(J=[N]\) gives \(M_N(z)\ge\per S_N\).
\end{proof}

\subsection{Small-ball probabilities}

Conditionally on all but one row, the determinant is an affine
function of the entries of that row. A small-ball bound for such affine forms will control the negative part of the logarithm
when that row is replaced. The bound must be uniform in the
coefficients, since these depend on the other rows.

For a complex random variable \(Z\), recall the Lévy concentration
function and its Dini integral
\[
 Q_Z(r)=\sup_{w\in\C}\Pr(|Z-w|\le r),\qquad
 \Psi_Z(\delta)=\int_0^\delta Q_Z(r)\,\frac{\mathrm dr}{r}.
\]
In the next lemma we track how the constant depends on the level \(r_0\) in its hypothesis, since the entry law will later be approximated.

\begin{lemma}[Small-ball bound for a linear form]\label{pd:row}
Let \(Z_1,Z_2,\ldots\) be independent copies of a centered complex
random variable \(Z\) with \(\E|Z|^2=1\). Suppose that
\(Q_Z(r_0)\le1/2\) for some \(0<r_0\le1\). Then, for all
coefficients \(b_j\) with \(\sum_j|b_j|^2=1\) and all \(0<r\le1\),
\begin{equation}\label{pd:small}
 Q_{\sum b_jZ_j}(r)\le Q_Z(\sqrt r)+C_{r_0}\sqrt r .
\end{equation}
In particular, for every \(0<\delta\le1\),
\[
 \int_0^\delta Q_{\sum_j b_jZ_j}(r)\,\frac{\mathrm dr}{r}
 \le 2\Psi_Z(\sqrt\delta)+2C_{r_0}\sqrt\delta .
\]
If \(\Psi_Z(1)<\infty\), these normalized linear forms therefore satisfy the planar Dini condition, with bounds uniform in the number of
summands and the coefficients.
\end{lemma}
\begin{proof}
If \(|b_j|\ge\sqrt r\) for some \(j\), condition on all variables
except \(Z_j\). For each \(w\), the values of \(Z_j\) with
\(|\sum_kb_kZ_k-w|\le r\) form a disk of radius \(r/|b_j|\le\sqrt r\),
which gives the first term on the right
of \eqref{pd:small}. It remains to consider the case
\(\max_j|b_j|<\sqrt r\). Let \(Z^\Delta=Z-Z'\), where \(Z'\) is an
independent copy of \(Z\).
Conditioning on \(Z'\) gives \(\Pr(|Z^\Delta|\le r_0)\le1/2\).
Since \(\E|Z^\Delta|^2=2\), we also have
\(\Pr(|Z^\Delta|>4)\le1/8\). Consequently
\[
 v:=\E\bigl[|Z^\Delta|^2{\bf1}_{|Z^\Delta|\le4}\bigr]\ge 3r_0^2/8 .
\]
One of the projections \(P=\Re,\Im\) satisfies
\(\sum_j\E[P(b_jZ^\Delta)^2{\bf1}_{|Z^\Delta|\le4}]\ge v/2\), since
the sums for the two projections add up to \(v\). Write \(Y=\sum_j P(b_j Z_j)\)
and let \(\phi_j(s)\) be the characteristic function of \(P(b_j Z_j)\).
For \(|s|\le(4\sqrt r)^{-1}\), symmetrization and
\(1-\cos x\ge x^2/4\) on \([-1,1]\) give
\[
 1-|\phi_j(s)|^2=\E[1-\cos(sP(b_jZ^\Delta))]
 \ge \tfrac14s^2\E[P(b_jZ^\Delta)^2{\bf1}_{|Z^\Delta|\le4}].
\]
By independence, \(|\phi_Y(s)|^2=\prod_j|\phi_j(s)|^2\), so
\(1-x\le e^{-x}\) and the lower bound \(v/2\) give
\(|\phi_Y(s)|\le e^{-v s^2/16}\) on this interval.

To pass back to probabilities, take
\[
 K(x)=\left(\frac{\sin(x/2)}{x/2}\right)^2
     =\int_{-1}^{1}(1-|t|)e^{itx}\dd t .
\]
Then \(K\ge0\), \(K(0)=1\) and \(k=\inf_{[-1,1]}K>0\). For
\(h=4\sqrt r\) and any \(y\in\R\),
\[
 \Pr(|Y-y|\le h)
 \le k^{-1}\E K((Y-y)/h)
 \le k^{-1}h\int_{-1/h}^{1/h}|\phi_Y(s)|\dd s
 \le C_{r_0}\sqrt r .
\]
Since a planar disk of radius \(r\) projects into an interval of
radius \(r\le h\), this proves \eqref{pd:small}; the remaining assertions
follow by integration with the substitution \(t=\sqrt r\).
\end{proof}

\subsection{Approximation by bounded laws}

To prove the fourth-moment estimate of Proposition~\ref{pd:fourth-ds}, we multiply two permutation expansions of the determinant and two of
its complex conjugate, so each entry occurs at most four times. Since
\(\xi\) need not have a finite fourth moment, we carry out this
computation for bounded approximants with vanishing mixed third
moments. Centering and the vanishing of these moments then eliminate the
terms in which an entry occurs once or three times, and the surviving
entries, which occur in pairs or quadruples, are counted by matchings.
To make the mixed third moments vanish, we multiply a bounded
approximation of \(\xi\) by an independent real random variable that has vanishing third moment
and converges to one in mean square.

\begin{lemma}[Approximation of the entry law]\label{pd:approx}
There are couplings \((\xi,\eta_m)\) with the following properties.
\begin{enumerate}
\item Each approximant is centered and has unit variance,
\[
 \E\eta_m=0,\qquad \E|\eta_m|^2=1.
\]
\item The approximants converge to \(\xi\) in mean square,
\[
 \E|\eta_m-\xi|^2\longrightarrow0.
\]
\item For each fixed \(m\), the variable \(\eta_m\) is bounded and
has a bounded density with respect to planar Lebesgue measure.
Both bounds may depend on \(m\).
\item All mixed third moments vanish,
\[
 \E\eta_m^a\overline{\eta_m}^{\,b}=0
 \qquad(a,b\in\mathbb Z_{\ge0},\ a+b=3).
\]
\item For all sufficiently large \(m\) and every \(t>0\),
\begin{equation}\label{pd:envelope}
 Q_{\eta_m}(t)\le2Q_\xi(4t).
\end{equation}
\end{enumerate}
\end{lemma}
\begin{proof}
For all sufficiently large \(m\), let \(T_m\) have the conditional
law of \(\xi\) given \(|\xi|\le m\). Couple the two by
setting \(T_m=\xi\) on that event and, on its complement,
resampling \(T_m\) from this conditional law. Then
\[
 \E|T_m-\xi|^2\le2\E[|\xi|^2{\bf1}_{|\xi|>m}]
                +2\Pr(|\xi|>m)\E|T_m|^2\longrightarrow0 .
\]
Let \(U_{\rm disk}\) be independent and uniform on the disk of radius
\(\sqrt2\), so that \(\E U_{\rm disk}=0\) and
\(\E|U_{\rm disk}|^2=1\). Define
\[
 V_m=(T_m-\E T_m+m^{-1}U_{\rm disk})/\sigma_m,\qquad
 \sigma_m^2=\E|T_m-\E T_m|^2+m^{-2}\longrightarrow1 .
\]
Each \(V_m\) is bounded, and the independent noise
\(m^{-1}U_{\rm disk}\) gives it a bounded planar density. Since \(T_m\) converges to \(\xi\) in mean square,
\(\E T_m\to0\) and \(\sigma_m\to1\), so
\(\E|V_m-\xi|^2\to0\). Conditioning on \(U_{\rm disk}\) gives
\(Q_{V_m}(t)\le Q_\xi(\sigma_m t)/\Pr(|\xi|\le m)\).

Let \(W_m\) be an independent real multiplier defined by
\[
 W_m=\sqrt{\frac{m^3+1}{m^3+m^2}}\,
 \begin{cases}
  1,&\text{with probability }m^3/(m^3+1),\\
  -m,&\text{with probability }1/(m^3+1).
 \end{cases}
\]
Then \(\E W_m^3=0\), \(\E W_m^2=1\) and \(\E W_m\to1\), so
\(\E|W_m-1|^2=2(1-\E W_m)\to0\). Set \(\eta_m=W_m V_m\).
By independence, \(\eta_m\) is centered and has unit variance and
vanishing mixed third moments. Using
\(\eta_m-\xi=(W_m-1)V_m+(V_m-\xi)\) and \(\E|V_m|^2=1\), we obtain
\[
 \E|\eta_m-\xi|^2
 \le 2\E|W_m-1|^2+2\E|V_m-\xi|^2\longrightarrow0.
\]
Conditionally on either value of \(W_m\), the variable \(\eta_m\) has
a bounded density. For large \(m\), we have \(|W_m|\ge1/2\),
\(\sigma_m\le2\) and \(\Pr(|\xi|\le m)\ge1/2\). Conditioning on
\(W_m\) in the concentration bound for \(V_m\) therefore gives
\eqref{pd:envelope}. In particular, one may use the same \(r_0\) in
Lemma~\ref{pd:row} for \(\xi\) and for all large \(m\), since
\(Q_\xi(t)\to0\) as \(t\downarrow0\) by \eqref{pd:law}.
\end{proof}

The approximants are complex even when \(\xi\) is real, because the
noise \(m^{-1}U_{\rm disk}\) smooths in the plane. Their first and
absolute second moments match those of \(\xi\), so
Lemma~\ref{pd:identity} gives the same \(M_N\).

Boundedness of \(\eta_m\) provides the finite moments needed for the
fourth-moment expansion, and its bounded planar density gives the
power bound \(Q_{\eta_m}(t)\le C_m t^2\) used for logarithmic
concentration in Lemma~\ref{pd:row-concentration}. The constants may
deteriorate with \(m\); this is harmless, since concentration is
applied at fixed \(m\) as \(N\to\infty\).

For the return to \(\xi\), the common envelope \eqref{pd:envelope}
provides a bound independent of \(C_m\). In particular,
\[
 \int_0^\delta Q_{\eta_m}(\sqrt t)\,\frac{\mathrm dt}{t}
 \le4\Psi_\xi(4\sqrt\delta).
\]
The right side is independent of \(m\) and tends to zero as
\(\delta\downarrow0\) by the planar Dini condition; this controls the
truncated logarithmic errors in row replacement uniformly in
large \(m\).

For two matrices with the same variance profile and coupled entry
laws, replacing rows one at a time reduces the difference of the log determinants to a coupling estimate and the small-ball bound of Lemma~\ref{pd:row}.

\begin{samepage}
\begin{proposition}[Row replacement for log determinants]\label{pd:transport}
Let \((\xi,\eta)\) be a coupling of two centered complex random
variables with \(\E|\xi|^2=\E|\eta|^2=1\) and
\(\Psi_\xi(1)+\Psi_\eta(1)<\infty\), and let \(r_0\in(0,1]\) satisfy
\(Q_\xi(r_0),Q_\eta(r_0)\le1/2\). Given a doubly stochastic profile
\(S_N\), form \(X_N=(\sqrt{s_{ij}}\,\xi_{ij})\) and
\(\widetilde X_N=(\sqrt{s_{ij}}\,\eta_{ij})\) from independent copies
\((\xi_{ij},\eta_{ij})\) of the coupling, and let
\(\widetilde D_N(z)=\det(zI-\widetilde X_N)\). Then for all \(z\in\C\)
and \(\delta\in(0,1]\),
\begin{equation}\label{pd:coupling}
 \frac1N\E\bigl|\log|D_N(z)|-\log|\widetilde D_N(z)|\bigr|
 \le\frac{\sqrt{\E|\xi-\eta|^2}}{\delta}
 +2\Psi_\xi(\sqrt\delta)+2\Psi_\eta(\sqrt\delta)+C_{r_0}\sqrt\delta,
\end{equation}
where \(C_{r_0}\) depends only on \(r_0\).
\end{proposition}
\end{samepage}
\begin{proof}
We first show that all logarithms below are integrable. A nonzero
polynomial \(F\) that is affine in each row is affine in each
random entry, which we call a label. If \(F\)
depends on a label \(Z\), write \(F=ZB+C\). By induction on the
number of labels, \(B\ne0\) almost surely and \(\log^-|B|\) is
integrable. Conditionally on the remaining labels, the planar Dini
condition gives
\[
 \E_Z\log^-|F|
 \le\log^-|B|+\E_Z\log^-|Z+C/B|
 \le\log^-|B|+\Psi_Z(1),
 \qquad \log^\pm t=\max(0,\pm\log t).
\]
Here \(\E_Z\) denotes averaging over the label \(Z\) alone. The
induction starts from a nonzero constant, and the positive part of
the logarithm is integrable because \(F\) has a finite second moment.
If each row's entry law satisfies the planar Dini condition, the argument allows different entry laws in different rows. By the
computation in Lemma~\ref{pd:identity}, every intermediate
determinant obtained by replacing rows has second moment
\(\per(S_N+|z|^2I)>0\), so it is not identically zero as a polynomial
in the random entries.

Replace the rows of \(X_N\) one by one by those of
\(\widetilde X_N\). Conditionally on the other rows, expansion along
the row being replaced gives the two determinants as
\[
 Y_\xi=b+\sum_jc_j\xi_j,\qquad
 Y_\eta=b+\sum_jc_j\eta_j.
\]
The coefficients are the same in both expressions, since they
depend only on the other rows. Write
\(\sigma^2=\sum_j|c_j|^2\). By the integrability argument, the conditional forms are almost
surely nonzero; if \(\sigma=0\), they are equal constants and
contribute no error.

For \(\sigma>0\), clip both logarithms below
\(\log(\delta\sigma)\). The clipped function
\(y\mapsto\log\max(|y|,\delta\sigma)\) is
\(1/(\delta\sigma)\)-Lipschitz. Independence and centering of the
coupled differences give
\[
 \E|Y_\xi-Y_\eta|^2=\sigma^2\E|\xi-\eta|^2
\]
conditionally on the other rows. Thus the expected absolute difference of the clipped logarithms is at most
\(\sqrt{\E|\xi-\eta|^2}/\delta\).
By Lemma~\ref{pd:row}, the expected error from clipping
\(\log|Y_\xi|\) is
\[
 \int_0^\delta \Pr(|Y_\xi|\le t\sigma)\,\frac{\mathrm dt}{t}
 \le 2\Psi_\xi(\sqrt\delta)+2C_{r_0}\sqrt\delta,
\]
and the same bound holds for \(Y_\eta\) with \(\Psi_\eta\) in place of
\(\Psi_\xi\).
Summing these errors over the \(N\) rows and dividing by \(N\)
proves \eqref{pd:coupling}, after enlarging \(C_{r_0}\).
\end{proof}

\begin{samepage}
\begin{corollary}[Approximation of log determinants]\label{pd:row-approx}
Let \(X_N\) and \(X_N^{(m)}=(\sqrt{s_{ij}}\,\eta_{ij}^{(m)})\) be formed
from independent copies \((\xi_{ij},\eta_{ij}^{(m)})\) of the coupling
\((\xi,\eta_m)\) of Lemma~\ref{pd:approx}, and let
\(D_N^{(m)}(z)=\det(zI-X_N^{(m)})\). Then
\[
 \lim_{m\to\infty}\sup_{N,S_N}\sup_{z\in\C}
 \frac1N\E\left|
   \log|D_N(z)|-\log|D_N^{(m)}(z)|
 \right|=0,
\]
where the first supremum ranges over all sizes \(N\) and
doubly stochastic profiles \(S_N\).
\end{corollary}
\end{samepage}
\begin{proof}
By \eqref{pd:envelope}, a common \(r_0\) serves for \(\xi\) and all
sufficiently large \(m\), and
\(\Psi_{\eta_m}(\sqrt\delta)\le2\Psi_\xi(4\sqrt\delta)\).
Proposition~\ref{pd:transport} therefore yields
\[
 \frac1N\E\left|
   \log|D_N(z)|-\log|D_N^{(m)}(z)|
 \right|
 \le\frac{\sqrt{\E|\eta_m-\xi|^2}}{\delta}
       +C\bigl\{\sqrt\delta+\Psi_\xi(4\sqrt\delta)\bigr\},
\]
uniformly in \(N,S_N,z\), where \(C\) is the same for all
sufficiently large \(m\). Let \(m\to\infty\) at fixed \(\delta\),
using mean-square convergence, and then let \(\delta\downarrow0\),
using the planar Dini condition.
\end{proof}

Under the same entry-law assumptions, the integrability argument and
the replacement bound in the proof of Proposition~\ref{pd:transport}
hold for every nonzero polynomial that is affine in each row, with the
logarithmic difference divided by the number of rows. We use this extension in Section~\ref{sec:variance-profiles} for the determinant and the permanent of \(zI-A_N-X_N\).

\subsection{The relative fourth moment}

Row replacement reduces the logarithmic comparison to the bounded
approximating law \(\eta_m\) at each fixed \(m\). For fixed \(m\) and
\(z\), the approximating determinant \(D_N^{(m)}(z)\) has second
moment \(M_N(z)\) by \eqref{pd:second-ds}. The Paley--Zygmund
inequality gives
\[
 \Pr\!\left(|D_N^{(m)}(z)|^2\ge\tfrac12 M_N(z)\right)
 \ge\frac{M_N(z)^2}{4\E|D_N^{(m)}(z)|^4}.
\]
The following proposition shows that the right side is at least
\(e^{-o(N)}\); in the next subsection we combine this with exponential
concentration of the logarithm to close the Jensen gap.

\begin{samepage}
\begin{proposition}[Relative fourth moment for doubly stochastic profiles]\label{pd:fourth-ds}
For each fixed \(m\), the approximating determinant \(D_N^{(m)}\) of
Corollary~\ref{pd:row-approx} satisfies
\begin{equation}\label{pd:ds-moment}
 \sup_{z\in\C}\frac1N\log
 \frac{\E|D_N^{(m)}(z)|^4}{M_N(z)^2}\longrightarrow0 .
\end{equation}
\end{proposition}
\end{samepage}
\begin{proof}
Fix \(m\) and \(z\), and write \(k=\E|\eta_m|^4\) and \(q=|z|^2\).
We first bound \(\E|D_N^{(m)}(z)|^4\) by an expansion with nonnegative coefficients. Let \(\chi\) be a symmetric real random
variable with unit variance and fourth moment \(k\); for example,
\(\chi\) may put mass \(1/(2k)\) at each of \(\pm\sqrt{k}\) and the
remaining mass at zero. With independent copies \(\chi_{ij}\), let
\(X_N^{\rm ref}=(\sqrt{s_{ij}}\chi_{ij})_{i,j}\). Then
\begin{equation}\label{pd:dominate}
 \E|D_N^{(m)}(z)|^4\le \E\bigl[\per(|z|I+X_N^{\rm ref})^4\bigr].
\end{equation}
Indeed, expand
\(D_N^{(m)}(z)^2\overline{D_N^{(m)}(z)}^{\,2}\)
into its four permutation sums. Terms in which some random entry
occurs once or three times have zero expectation, by centering and
the vanishing mixed third moments of \(\eta_m\). The mixed moments of
\(\eta_m\) of total degree two and four are bounded in modulus by one
and \(k\), respectively. Taking absolute values of the coefficients and applying these bounds,
we dominate the expectation of each term by that of the corresponding
term in the real permanent expansion on the right. Both polynomials
have second moment \(M_N(z)\). The reference variables serve only to
bound moments, so their atoms play no role in the logarithmic comparison.

To count the surviving terms, multiply each row of
\(|z|I+X_N^{\rm ref}\) by an independent fair sign \(\omega_i\), and denote the
resulting entries by
\[
 x_{ij}^{\rm sym}=\omega_i\bigl(|z|{\bf1}_{i=j}+(X_N^{\rm ref})_{ij}\bigr).
\]
The row sign symmetrizes the diagonal entry; the other entries are
already symmetric, so after multiplication they keep their joint law
and are independent of the sign. Thus all \(x_{ij}^{\rm sym}\) are
independent and symmetric, and the absolute value of the permanent is
unchanged. We normalize the second moments and record the
standardized fourth moments,
\[
 T_{ij}=\frac{\E|x_{ij}^{\rm sym}|^2}{1+q}
       =\frac{s_{ij}+q{\bf1}_{i=j}}{1+q},
 \qquad
 k_{ij}=\frac{\E|x_{ij}^{\rm sym}|^4}{(\E|x_{ij}^{\rm sym}|^2)^2}.
\]
Entries of zero variance have zero weight and may be assigned
\(k_{ij}=1\). The matrix \(T\) is doubly stochastic.

The square of the second-moment expansion is a sum over pairs of
permutations. We sample two permutations \(\pi\) and \(\sigma\) independently with law
\(\Pr_T(\pi)=\prod_iT_{i\pi(i)}/\per T\) and write \(\E_T\) for
expectation over the pair. Each permutation defines a perfect
matching between disjoint sets of row and column vertices, in which
\(\pi\) joins row \(i\) to column \(\pi(i)\). We superpose the matchings
of \(\pi\) and \(\sigma\), keeping both copies of any shared edge. Since
every vertex meets one edge of each matching, each connected component
is either a shared edge counted twice or
an even cycle whose edges alternate between the two matchings. A cycle
on \(\ell\) rows has \(\ell\) columns and \(2\ell\) edges.
Equivalently, shared edges are the fixed points of
\(\sigma^{-1}\pi\), and nontrivial cycles correspond to its cycles
of length at least two. Write \(\mathrm{cyc}(\pi,\sigma)\) for the number of
nontrivial cycles.

In the fourth-moment expansion, each surviving edge has multiplicity
two or four. Halving these multiplicities leaves a union of two perfect matchings, so the components are again shared edges and alternating cycles. A nontrivial cycle admits \(\binom{4}{2}=6\) assignments to the four
factors. Choosing which two factors use a specified alternating
matching forces the other two to use its complement. In the square of the second moment, by contrast, the same cycle is counted only twice, once for each
ordering of its two matchings. This gives a factor of \(6/2=3\) per
nontrivial cycle. A shared edge contributes
\(k_{ij}\), since its fourth moment replaces the square of its second
moment. This is the counting used by Karmarkar, Karp, Lipton,
Lov\'asz and Luby for random-determinant estimators of the permanent
\cite{pd-kklll}. Consequently
\begin{equation}\label{pd:matching-ratio}
 \frac{\E[\per(|z|I+X_N^{\rm ref})^4]}{M_N(z)^2}
 =\E_T\!\left[
    3^{\mathrm{cyc}(\pi,\sigma)}
    \prod_{i:\,\pi(i)=\sigma(i)}k_{i,\pi(i)}
 \right],
\end{equation}
where the common row normalization cancels from the ratio.
For example, \(\pi=\mathrm{id}\) and
\(\sigma=(1\,2)(3\,4\,5)(6)\) give a four-edge cycle, a six-edge
cycle and the shared edge \((6,6)\); their extra factor is
\(3^2k_{66}\).

To apply Lemma~\ref{pd:matchings} to the weighted component sum in
\eqref{pd:matching-ratio}, we need an exponential lower bound on the
permanent and the bound \(O(N\epsilon_N)\) on the excess weight of
shared edges, \(\sum_{i,j}(k_{ij}-1)T_{ij}^2\).

For off-diagonal edges of positive variance, \(k_{ij}=k\), and
\(T_{ij}\le\epsilon_N\). Thus
\[
 \sum_{i\ne j}(k_{ij}-1)T_{ij}^2
 \le(k-1)\epsilon_N\sum_{i\ne j}T_{ij}
 \le(k-1)N\epsilon_N.
\]
The diagonal requires a separate computation, since \(T_{ii}\) can be
of order one. The entry \(x_{ii}^{\rm sym}\) has second moment \(q+s_{ii}\)
and fourth moment \(q^2+6qs_{ii}+ks_{ii}^2\). Hence
\(k_{ii}\le\max(3,k)\), and
\[
 (k_{ii}-1)T_{ii}^2
 =\frac{4qs_{ii}+(k-1)s_{ii}^2}{(1+q)^2}
 \le s_{ii}+(k-1)s_{ii}^2\le k\epsilon_N,
\]
using \(4q/(1+q)^2\le1\) and \(s_{ii}\le\epsilon_N\le1\).
Adding the two contributions gives
\[
 \sum_{i,j}(k_{ij}-1)T_{ij}^2\le(2k-1)N\epsilon_N.
\]
Thus \(K_0=\max(3,k)\) and \(C_0=2k-1\) satisfy the moment
conditions of Lemma~\ref{pd:matchings} with \(\epsilon=\epsilon_N\),
uniformly in \(z\).
They depend only on the fixed approximating law and not on any
fourth moment of \(\xi\).

The van der Waerden bound \cite{pd-egorychev,pd-falikman,pd-vdw}
gives \(\per T\ge N!/N^N\ge e^{-N}\). Combining
\eqref{pd:dominate}, \eqref{pd:matching-ratio} and
Lemma~\ref{pd:matchings} gives, for all sufficiently large \(N\),
\[
 \sup_{z\in\C}\frac1N\log
 \frac{\E|D_N^{(m)}(z)|^4}{M_N(z)^2}
 \le\frac{C_m}{\log(1/\epsilon_N)}\longrightarrow0.
\]
Since the ratio is at least one by the Cauchy--Schwarz inequality,
this proves \eqref{pd:ds-moment}.
\end{proof}

It remains to prove the weighted cycle bound used in
Proposition~\ref{pd:fourth-ds}. We allow row sums up to one and
column sums up to two, since Section~\ref{sec:variance-profiles}
needs this generality.

\begin{samepage}
\begin{lemma}[Weighted cycle bound]\label{pd:matchings}
Let \(T\) be a nonnegative \(N\times N\) matrix with positive permanent,
row sums at most one, column sums at most two and off-diagonal
entries at most \(\epsilon\), and let \(1\le k_{ij}\le K_0\), with
\(K_0\ge3\), satisfy \(\sum_{i,j}(k_{ij}-1)T_{ij}^2\le C_0N\epsilon\).
Let \(\pi\) and \(\sigma\) be independent random permutations with law
\(\Pr_T(\pi)=\prod_i T_{i\pi(i)}/\per T\), write \(\E_T\) for
expectation over the pair, and let \(\mathrm{cyc}(\pi,\sigma)\) be the
number of nontrivial alternating cycles in their union. Then, for all sufficiently small \(\epsilon>0\),
\begin{equation}\label{pd:matching-exp}
 1\le\E_T\!\left[
     3^{\mathrm{cyc}(\pi,\sigma)}\prod_{i:\,\pi(i)=\sigma(i)}k_{i,\pi(i)}
 \right]
 \le\exp\!\left\{
    \frac{C\bigl(N+\log(1/\per T)\bigr)}{\log(1/\epsilon)}
 \right\}.
\end{equation}
Here the constant \(C\) and the smallness threshold depend only on
\(K_0,C_0\). In particular, if \(\epsilon\to0\) and
\(-\log\per T=O(N)\), the expectation is \(e^{o(N)}\).
\end{lemma}
\end{samepage}
\begin{proof}
Fix an integer \(L\ge2\), and call a cycle long if it contains at least \(L\) rows and short otherwise. Since components are disjoint, there are at most \(N/L\) long
cycles, so their combined factor is at most \(3^{N/L}\) and they need
no weighted counting.
Short cycles and shared edges will instead be controlled by their
total weights. We first prove, for every \(h\ge K_0\), the bound
\begin{equation}\label{pd:count}
 \begin{aligned}
 0&\le\frac1N\log\E_T\!\left[
     3^{\mathrm{cyc}(\pi,\sigma)}\prod_{i:\,\pi(i)=\sigma(i)} k_{i,\pi(i)}
 \right]\\
 &\le\frac{\log3}{L}
 +\frac{\log K_0}{\log h}
       \left(-\frac2N\log\per T+Ch\epsilon\,2^L\right).
 \end{aligned}
\end{equation}
To count short cycles, first sum over their column vertices. For
distinct rows \(i,j\), the quantity
\[
 B_{ij}=\sum_a T_{ia}T_{ja},\qquad B_{ii}=0,
\]
is the total weight of the two-edge paths from row \(i\) through a column to row \(j\). Since the row sums of \(T\) are at most one and its
column sums at most two, \(B\) has row sums at most two. We also have \(B_{ij}\le3\epsilon\), because each of the two terms
with \(a\in\{i,j\}\) contains an off-diagonal factor and is at most
\(\epsilon\), while the remaining terms contain the off-diagonal factor
\(T_{ia}\) and sum to at most \(\epsilon\) by the row-sum bound for
row \(j\).
Multiplying \(\ell\) such path weights and taking a trace closes
the walk. Consequently the total weight of alternating cycles on
\(\ell\ge2\) rows, with the roles of the two matchings specified, is at
most
\[
 \frac1\ell\operatorname{Tr}(B^\ell)
 \le\frac{3N\epsilon\,2^{\ell-1}}\ell.
\]
The trace also counts walks that revisit rows or columns, which only
increases it, and division by \(\ell\) removes the choice of starting
row. The final inequality follows by bounding one factor of \(B\) by
\(3\epsilon\) and summing the remaining walk using the row-sum bound.

Let \(W\) be the product of the short-cycle factors and all
shared-edge factors. A direct estimate of \(\E_T W\) would carry the
normalizing factor \((\per T)^{-2}\), which can be exponentially
large. We therefore estimate \(\E_T W^t\), where
\(t=\log h/\log K_0\ge1\), and take a \(t\)-th root at the end, which reduces this factor to
\((\per T)^{-2/t}\).
Expand \(W^t\) by writing \(x^t=1+(x^t-1)\) for each component factor
\(x\) of \(W\). Choosing the second term marks that component. The marked shared edges
\(\mathcal M_1\) and alternating cycles \(\mathcal M_2\), with the roles
of the two matchings specified, form a disjoint family \(\mathcal M\)
occupying row set \(I\) and column set \(J\). Once this family is
fixed, the two matchings may be completed independently on the
complementary minor. For an edge \(e=(i,j)\), write \(T_e=T_{ij}\)
and \(k_e=k_{ij}\). Since a shared edge is used by both matchings and a
cycle edge by only one, we obtain
\[
 (\per T)^2\,\E_T W^t
 =\sum_{\mathcal M}\bigl(\per T[I^c,J^c]\bigr)^2
   \prod_{e\in\mathcal M_1}(k_e^t-1)T_e^2
   \prod_{c\in\mathcal M_2}\left[(3^t-1)\prod_{e\in c}T_e\right].
\]
The row-sum bound gives \(\per T[I^c,J^c]\le1\) for every family.
We now bound the remaining sum by a product. The inequality
\(x^t-1\le(h-1)(x-1)/(K_0-1)\), valid for \(1\le x\le K_0\), bounds
the total weight of marked shared edges by \(O(hN\epsilon)\) and that
of marked cycles by \(O(hN\epsilon\,2^L)\). Allowing marked components
to overlap adds only nonnegative terms, and the resulting sum
factorizes over components, so \(1+x\le e^x\) gives
\[
 \E_T W^t\le(\per T)^{-2}\exp(ChN\epsilon\,2^L).
\]
Taking a \(t\)-th root and restoring the factor \(3^{N/L}\) from
long cycles gives \eqref{pd:count}. Finally choose
\(h=\epsilon^{-1/4}\) and
\(L=\lfloor\log(1/\epsilon)/(4\log2)\rfloor\). For sufficiently
small \(\epsilon\), these choices are admissible and
\(h\epsilon2^L\le\epsilon^{1/2}\). Since the row-sum bound implies
\(\per T\le1\), \eqref{pd:count} gives \eqref{pd:matching-exp}.
\end{proof}

\subsection{The logarithmic comparison}

The following comparison of \(u_N\) with \(p_N\) completes the
probabilistic part of the proof of the circular law; what remains is
the limit of the log permanent.

\begin{samepage}
\begin{proposition}[Logarithmic comparison for doubly stochastic profiles]\label{pd:compare-ds}
Suppose that \(S_N\) is doubly stochastic and \(\epsilon_N\to0\), and
that the \(\xi_{ij}\) are independent copies of a law \(\xi\)
satisfying \eqref{pd:law}. Then, for \(X_N=(\sqrt{s_{ij}}\xi_{ij})\),
\begin{equation}\label{pd:ds-gap}
 \sup_{z\in\C}\E|u_N(z)-p_N(z)|\longrightarrow0.
\end{equation}
\end{proposition}
\end{samepage}
\begin{proof}
Fix \(m\) and work first with \(\eta_m\). By \eqref{pd:second-ds},
\(p_N\) is also the second-moment pressure of the approximating law.
Since \(Q_{\eta_m}(t)\le C_mt^2\), Lemma~\ref{pd:row-concentration}
below applies, and together with Proposition~\ref{pd:fourth-ds} it
gives
\[
 \sup_{z\in\C}\E\left|
   \frac1N\log|D_N^{(m)}(z)|-p_N(z)
 \right|\longrightarrow0
 \qquad(N\to\infty)
\]
for each fixed \(m\).

For \(\xi\) itself, the triangle inequality gives
\[
 \begin{aligned}
 \sup_{z\in\C}\E|u_N(z)-p_N(z)|
 &\le\sup_{z\in\C}\frac1N\E
       \left|\log|D_N(z)|-\log|D_N^{(m)}(z)|\right|\\
 &\quad+\sup_{z\in\C}\E
       \left|\frac1N\log|D_N^{(m)}(z)|-p_N(z)\right|.
 \end{aligned}
\]
The second term vanishes as \(N\to\infty\) at fixed \(m\), and by
Corollary~\ref{pd:row-approx} the first term vanishes as
\(m\to\infty\), uniformly in \(N\) and \(z\). Letting first
\(N\to\infty\) and then \(m\to\infty\) proves \eqref{pd:ds-gap}.
\end{proof}

The concentration estimate used in
Proposition~\ref{pd:compare-ds} is Lemma~\ref{pd:row-concentration}, which we derive from a general resampling bound. That bound, Lemma~\ref{pd:resampling}, is stated for a function of general independent random elements, so that it applies
both to the rows considered here and to the stars used for
reciprocal pairs in Section~\ref{sec:reciprocal}.

\begin{lemma}[From resampling to logarithmic comparison]\label{pd:resampling}
Let \(F\) be a function of \(N\) independent random elements with
\(0<\E|F|^2<\infty\) and \(Y=\log|F|\) integrable, and let
\(Y^{(i)}\) denote the logarithm obtained by replacing element \(i\)
by an independent copy and keeping all other elements fixed. Assume
that for some \(a>0\), \(B<\infty\) and every \(i\),
\[
 \E\big[\exp(a|Y-Y^{(i)}|)\,\big|\,\text{all elements except }i\big]\le B
 \quad\text{almost surely}.
\]
Then, for \(s>0\), with constants depending only on \(a,B\),
\begin{equation}\label{pd:bern}
 \Pr(|Y-\E Y|\ge s)\le 2\exp[-c\min(s^2/N,s)],\qquad
 \E|Y-\E Y|^4\le C N^2 .
\end{equation}
If \(\E|F|^4<\infty\), set
\(r=N^{-1}\log\bigl(\E|F|^4/(\E|F|^2)^2\bigr)\). Then
\begin{equation}\label{pd:gap}
 \E\left|\frac{Y}{N}-\frac1{2N}\log\E|F|^2\right|
       \le C(\sqrt r+r+N^{-1/2}) .
\end{equation}
\end{lemma}
\begin{proof}
Let \(\mathcal F_i\) be the sigma-field generated by the first \(i\)
elements. The martingale increments are
\[
 \Delta_i=\E[Y\mid\mathcal F_i]-\E[Y\mid\mathcal F_{i-1}]
          =\E[Y-Y^{(i)}\mid\mathcal F_i],
\]
where the last equality uses independence. The conditional Jensen
inequality and the hypothesis give
\(\E[e^{a|\Delta_i|}\mid\mathcal F_{i-1}]\le B\). The linear term in
the conditional exponential moment vanishes because
\(\E[\Delta_i\mid\mathcal F_{i-1}]=0\), so Taylor expansion yields
\(\E[e^{\lambda\Delta_i}\mid\mathcal F_{i-1}]\le e^{C\lambda^2}\)
for \(|\lambda|\le a/2\). Iterating this bound and applying the
exponential Markov inequality gives the tail bound in
\eqref{pd:bern}; integrating the tail gives the fourth-moment bound.

The Paley--Zygmund inequality supplies the lower bound
\[
 \Pr\!\left(Y\ge\tfrac12\log(\E|F|^2/2)\right)\ge e^{-Nr}/4.
\]
If \(s=\tfrac12\log(\E|F|^2/2)-\E Y>0\), this event is a deviation
of at least \(s\) above the mean. Comparing its probability with
\eqref{pd:bern} gives \(\min(s^2/N,s)\le C(Nr+1)\). In the other direction, Jensen's inequality gives \(\E Y\le\tfrac12\log\E|F|^2\), and \eqref{pd:bern}
gives \(\E|Y-\E Y|\le C\sqrt N\). Combining these bounds proves
\eqref{pd:gap}, whether or not \(s>0\).
\end{proof}

For the approximating determinants, it remains to verify the
resampling hypothesis of Lemma~\ref{pd:resampling}. We do so for all polynomials affine in each row, so that the next lemma also applies to the determinant and the permanent of \(zI-A_N-X_N\).

\begin{samepage}
\begin{lemma}[Concentration for polynomials affine in each row]\label{pd:row-concentration}
Let \(F_N(z)\) be a polynomial in \(z\) and in the entries of \(N\)
rows, where all entries are independent copies of a
centered unit-variance law \(Z\). Suppose that \(F_N\) is monic in
\(z\) and affine in each row, meaning that fixing all other rows leaves a
constant plus a linear form in the entries of that row. Assume
\(Q_Z(t)\le Bt^\beta\) for \(0<t\le1\), with \(B,\beta>0\). Then the
following hold.
\begin{enumerate}
\item At each fixed \(z\) where \(F_N(z)\) is not identically zero as
a polynomial in the random entries, \(Y=\log|F_N(z)|\) satisfies the
concentration bounds \eqref{pd:bern}.
\item If moreover \(\E|F_N(z)|^4<\infty\), then \(Y\) satisfies the
logarithmic comparison \eqref{pd:gap}, with the relative
fourth-moment parameter \(r\) defined there for \(F=F_N(z)\).
\item Writing \(u_N^F(z)=N^{-1}\log|F_N(z)|\), we have
\begin{equation}\label{pd:area-fluct}
 \E\|u_N^F-\E u_N^F\|_{L^1(D_R)}^4\le C_RN^{-2}.
\end{equation}
\end{enumerate}
In (1) and (2) the constants depend only on \(B,\beta\), and all
constants are independent of the coefficients of the polynomial.
\end{lemma}
\end{samepage}
\begin{proof}
By Lemma~\ref{pd:row}, the small-ball probabilities of normalized
row forms are at most \(Cr^{\min(\beta,1)/2}\). Conditionally on the other rows, we write the affine form
as \(L=b+\sum_jc_jZ_j\) and set
\[
 \sigma^2=\sum_j|c_j|^2,\qquad
 s_L^2=\E_{\rm row}|L|^2=|b|^2+\sigma^2 .
\]
Here \(\E_{\rm row}\) averages only over the row being resampled.
For a nonzero constant form, \(|L|/s_L=1\); otherwise \(\sigma>0\).
Markov's inequality bounds the upper tail of \(\log(|L|/s_L)\) by
\(e^{-2s}\). For the lower tail, consider the event \(|L|\le s_Le^{-s}\). If
\(|b|/\sigma\le e^{s/2}\), then \(s_Le^{-s}\le\sqrt2e^{-s/2}\sigma\), so
the small-ball estimate applies.
Otherwise, for \(s\ge\log4\), the event requires \(|L-b|\ge|b|/2\),
which has probability at most \(4\sigma^2/|b|^2\le4e^{-s}\). Thus
\[
 \Pr\!\left(\left|\log\frac{|L|}{s_L}\right|>s\right)\le Ce^{-cs},
\]
uniformly in the coefficients of every nonzero affine form.

Conditionally on the other rows, \(s_L\) is the same for the row and
for its independent copy, so \(Y-Y^{(i)}\) is the difference of their
normalized logarithms \(\log(|L|/s_L)\), each with uniform exponential
tails. By the triangle inequality and the Cauchy--Schwarz inequality,
\(Y-Y^{(i)}\) satisfies the hypothesis of Lemma~\ref{pd:resampling}. The induction
in the proof of Proposition~\ref{pd:transport} supplies
integrability. Hence Lemma~\ref{pd:resampling} applies at each \(z\) where
\(F_N(z)\) is not identically zero as a polynomial in the random
entries. The exceptional points
lie among the roots of the monic polynomial obtained by setting all
labels to zero. Let \(|D_R|=\pi R^2\) denote the area of the disk.
Hölder's inequality and \eqref{pd:bern} give
\[
 \E\|u_N^F-\E u_N^F\|_{L^1(D_R)}^4
 \le |D_R|^3\int_{D_R}\E|u_N^F(z)-\E u_N^F(z)|^4\dd A(z)
 \le C_R N^{-2}.
\]
Markov's inequality and the Borel--Cantelli lemma give almost sure
concentration in local \(L^1\).
\end{proof}

\subsection{The pressure limit and the circular law}

Proposition~\ref{pd:compare-ds} reduces the problem to identifying the limit of the second-moment pressure \(p_N\). The following result, applied with \(q=|z|^2\), gives \(p_N(z)\to U_{\rm circ}(z)\).

\begin{proposition}[Doubly stochastic pressure]\label{vp:ds}
If \(S_N\) is doubly stochastic and \(\epsilon_N\to0\), then for every
\(q\ge0\)
\[
 \lim_{N\to\infty}N^{-1}\log\per(S_N+qI)=
 \begin{cases}q-1,&0\le q\le1,\\ \log q,&q>1.\end{cases}
\]
\end{proposition}
\begin{proof}
The lower bound uses a one-parameter trial matrix in Gurvits'
inequality. For the upper bound, we pass to a circular Gaussian matrix with the
same variance profile, and hence the same second moment of the
determinant; Gaussian concentration then reduces the bound to the mean
of a regularized logarithm. We identify its limit as \(N\to\infty\) at fixed regularization and let the regularization tend to zero only in the limiting free-probability expression, so the argument needs no estimate on
small singular values.

For \(q>0\) and \(0\le r\le1\), use the doubly stochastic trial matrix \(P=rS_N+(1-r)I\) in Gurvits'
permanent lower bound \cite[Theorem~2.2]{vp-bethe}. With this choice,
the bound for \(\log\per(S_N+qI)\) is a sum over all entries
\((i,j)\) of the terms
\(P_{ij}\log((S_N+qI)_{ij}/P_{ij})+h(P_{ij})\), where
\(h(t)=(1-t)\log(1-t)\). On the diagonal, the log-sum inequality separates the contributions of the weights \(s_{ii}\) and
\(q\). The logarithmic terms then contribute at least
\[
 -Nr\log r+N(1-r)\log\frac q{1-r}.
\]
Off the diagonal, \(h(P_{ij})\ge-P_{ij}\); on the diagonal,
\(h(1-r+rs_{ii})\ge h(1-r)-rs_{ii}\), since \(h'\ge-1\).
Summing these terms gives at least \(Nr\log r-Nr\), whose first term
cancels the term \(-Nr\log r\) in the preceding display. Thus
\[
 N^{-1}\log\per(S_N+qI)\ge(1-r)\log\frac q{1-r}-r .
\]
Taking \(r=\max(1-q,0)\) gives the desired lower bound. At \(q=0\) the
van der Waerden bound gives the lower bound \(-1\).

For the upper bound, let \((G_N)_{ij}=\sqrt{s_{ij}}\,g_{ij}\), where
the labels \(g_{ij}\) are independent standard circularly symmetric
complex Gaussian variables. Fix \(z\) with \(|z|^2=q\). For
\(\delta>0\) and an \(N\times N\) matrix \(Y\), we use the regularized
logarithm
\begin{equation}\label{pd:regularized-log}
 R_{N,\delta}(Y;z)
 =\frac1{2N}\log\det((Y-zI)(Y-zI)^*+\delta^2 I).
\end{equation}
Regularization gives
\(|\det(G_N-zI)|^2\le e^{2N R_{N,\delta}(G_N;z)}\) and makes
\(R_{N,\delta}(G_N;z)\) a Lipschitz function of the real Gaussian
coordinates, with squared Lipschitz constant at most
\(C\epsilon_N/(N\delta^2)\). Indeed, the squared norm of the matrix
gradient is
\[
 \frac1{N^2}\sum_j\frac{\sigma_j^2}{(\sigma_j^2+\delta^2)^2}
 \le\frac1{4N\delta^2},
\]
where the \(\sigma_j\) are the singular values of \(G_N-zI\), and
passing to the label coordinates costs a factor of at most
\(C\epsilon_N\). The Gaussian
logarithmic Sobolev inequality \cite{ep-gross} and \eqref{pd:second-ds}
give
\[
 p_N(z)\le\frac1{2N}\log\E e^{2N R_{N,\delta}(G_N;z)}
 \le \E R_{N,\delta}(G_N;z)+\frac{C\epsilon_N}{\delta^2}.
\]

At fixed \(\delta\), the limit of \(\E R_{N,\delta}(G_N;z)\) is the
regularized logarithmic potential of a standard circular element \(c\), normalized by \(\varphi(c^*c)=1\).
To see this, expand a fixed normalized trace word in \(G_N,G_N^*\)
by Gaussian pairings. Each pairing identifies indices and gives
a connected graph, with one variance factor per paired edge.
If the graph has fewer vertices than the number of pairs plus one, we
sum along a spanning tree using the row and column bounds, and
estimate at least one remaining variance by \(\epsilon_N\), so the normalized
contribution is \(O(\epsilon_N)\).

In the remaining terms, the graph is a tree and the trace walk
traverses each edge twice, in opposite directions. These are the
noncrossing pairings of opposite letters, and summing along the
leaves shows that each has weight one after normalization, which
yields the moments of the circular element
\cite[Lecture~11]{fp-ns}.

Applying this argument to the moments of the squared singular values
of \(G_N-zI\) gives weak convergence of their expected empirical
measures, since the circular element has bounded norm. For each
fixed \(z\), these expected empirical measures have first moment
\(1+|z|^2\) for every \(N\). Together with the sublinear growth of the
logarithm, this controls the positive logarithmic tail uniformly in
\(N\), and the negative part is bounded for fixed \(\delta>0\). Consequently,
\[
 \lim_{N\to\infty}\E R_{N,\delta}(G_N;z)
 =\tfrac12\varphi\log((c-z)(c-z)^*+\delta^2).
\]
Since the error \(C\epsilon_N/\delta^2\) vanishes as
\(N\to\infty\) for fixed \(\delta\), we may then let
\(\delta\downarrow0\); the free regularized
potential tends to \(U_{\rm circ}(z)\) by the circular case of
\cite[Example~5.3]{ep-bl}. This gives the upper bound for every \(q\).
\end{proof}

The case \(q=0\) is a classical consequence of the van der Waerden
and Br\'{e}gman--Minc bounds; see \cite[Sections~7.2--7.3]{pd-barvinok-notes}.
For \(q>1\), we expand \(\per(I+S_N/q)\) over families of disjoint
directed cycles; dropping the disjointness requirement gives the
elementary bound
\[
 \log\per(I+S_N/q)
 \le\sum_{k\ge1}\frac{\operatorname{Tr}S_N^k}{kq^k}
 \le-N\epsilon_N\log(1-q^{-1}),
\]
since every entry of \(S_N^k\) is at most \(\epsilon_N\).
Together with the lower bound \(\per(S_N+qI)\ge q^N\) given by the
all-diagonal term, this proves the limit for \(q>1\); the endpoint
\(q=1\) follows by monotonicity.
Thus \(0<q<1\) is the range in which the upper bound requires an
argument beyond these elementary estimates.

For the uniform profile \(J_N=(1/N)_{ij}\), expansion in principal
permanents gives the classical formula
\[
 \per(J_N+qI)=\frac{N!}{N^N}\sum_{j=0}^N\frac{(Nq)^j}{j!},
\]
whose exponential rate is the one in Proposition~\ref{vp:ds};
compare \cite[Example~3.2]{vp-bethe}. When the largest entry tends to zero, the
proposition establishes the same rate for arbitrary supports. For fixed \(q>0\), every diagonal entry of
\((S_N+qI)/(1+q)\) is at least \(q/(1+q)\). This lower bound rules
out the moderate-deviation hypothesis in McCullagh's relative
asymptotic formula \cite{pd-mccullagh}.
The universal factor \(2^{N/2}\) in the Bethe approximation
\cite{vp-anari} leaves an uncertainty of order \(N\) in the
logarithm.

Two further facts are used to pass from logarithmic estimates to
zero measures. For a random monic polynomial \(F_N\) of degree \(N\),
we write \(u_N^F=N^{-1}\log|F_N|\) and
\(p_N^F=(2N)^{-1}\log\E|F_N|^2\).

\begin{lemma}[Absolute error and the Jensen gap]\label{pd:jensen}
At points where \(0<\E|F_N|^2<\infty\)
and \(\log|F_N|\) is integrable,
\begin{equation}\label{pd:jensen-abs}
 0\le \E|u_N^F-p_N^F|-(p_N^F-\E u_N^F)
       =2\E[(u_N^F-p_N^F)_+]\le N^{-1}.
\end{equation}
Here \(a_+=\max(a,0)\).
\end{lemma}
\begin{proof}
Since \(\E e^{2N(u_N^F-p_N^F)}=1\), Markov's inequality gives
\(\Pr(u_N^F-p_N^F>s)\le e^{-2Ns}\); integrating over \(s>0\) proves
the bound.
\end{proof}

Integrating \eqref{pd:jensen-abs} over a disk \(D_R\) shows that the
integrated Jensen gap and the expected \(L^1(D_R)\) error differ by at
most \(|D_R|/N\). Taking the distributional Laplacian turns local
\(L^1\) control of potentials into convergence against smooth
compactly supported test functions; this is the mechanism behind the
replacement principle of Tao and Vu
\cite[Theorem~2.1 and Section~3.6]{cl-taovu}. The next lemma records
this for comparison measures that may vary with \(N\), using their
tightness to reach the bounded-Lipschitz distance.

Recall that the bounded-Lipschitz distance between probability
measures \(\mu,\nu\) on \(\C\) is
\[
 d_{\BL}(\mu,\nu)=
 \sup_f\left|\int f\dd\mu-\int f\dd\nu\right|,
\]
where the supremum is over all \(f:\C\to\R\) with \(|f|\le1\) and
\(|f(z)-f(w)|\le|z-w|\) for all \(z,w\). It metrizes weak
convergence of probability measures.

\begin{lemma}[From logarithmic potentials to zero measures]\label{pd:root-return}
Let \(F_N\) be a random monic polynomial of degree \(N\) with
normalized zero measure \(\mu_N^F\), and let
\(V_N\in L^1_{\mathrm{loc}}(\C)\) be deterministic functions such that
\(\nu_N=(2\pi)^{-1}\Delta V_N\) is a probability measure, where
\(\Delta=\partial_x^2+\partial_y^2\). Suppose that the measures
\(\nu_N\) are tight in the sense that
\begin{equation}\label{pd:comparison-tightness}
 \lim_{R\to\infty}\sup_N\nu_N(D_R^c)=0.
\end{equation}
Then:
\begin{enumerate}
\item If \(\E\|u_N^F-V_N\|_{L^1(D_R)}\to0\) for every \(R>0\),
then \(\E d_{\BL}(\mu_N^F,\nu_N)\to0\).
\item If \(\|u_N^F-V_N\|_{L^1(D_R)}\to0\) almost surely for every
positive integer \(R\), then
\(d_{\BL}(\mu_N^F,\nu_N)\to0\) almost surely.
\end{enumerate}
\end{lemma}
\begin{proof}
Monicity gives \(u_N^F(z)=\int\log|z-w|\dd\mu_N^F(w)\). For every
smooth compactly supported real-valued test function \(\phi\),
\[
 \int\phi\dd(\mu_N^F-\nu_N)
   =\frac1{2\pi}\int_{\C}(\Delta\phi)(u_N^F-V_N)\dd A .
\]
This tends to zero in mean absolute value in case~(1), and on one
common event for all such tests in case~(2). Given \(\varepsilon>0\),
choose \(R\) with \(\nu_N(D_R)\ge1-\varepsilon\) for every \(N\).
Testing against a smooth cutoff with values in \([0,1]\), equal to
one on \(D_R\) and to zero outside \(D_{2R}\), we find that
\(\limsup_N\E\mu_N^F(D_{2R}^c)\le\varepsilon\) in case~(1), and
\(\limsup_N\mu_N^F(D_{2R}^c)\le\varepsilon\) on the common event in case~(2).
On \(D_{2R}\) the test functions defining \(d_{\BL}\) have a finite
uniform \(\varepsilon\)-net. Approximate its members by smooth
compactly supported functions and use convergence against those tests.
The mass outside the disk contributes \(O(\varepsilon)\); letting
\(\varepsilon\downarrow0\) proves both assertions.
\end{proof}

\begin{proof}[Proof of Theorem~\ref{vp:circle}]
Proposition~\ref{vp:ds} gives \(p_N(z)\to U_{\rm circ}(z)\) at every
fixed \(z\). The bounds
\(e^{-N}\le\per S_N\le\per(S_N+|z|^2I)\le(1+|z|^2)^N\), the last
from the product of the row sums, give
\(-1/2\le p_N(z)\le\tfrac12\log(1+|z|^2)\), so the convergence also
holds in local \(L^1\). Together with \eqref{pd:ds-gap}, this gives
\(\E\|u_N-U_{\rm circ}\|_{L^1(D_R)}\to0\).
The radial derivatives of the two branches of \(U_{\rm circ}\) agree at radius one, so \(\Delta U_{\rm circ}\) has no mass on the unit circle and \((2\pi)^{-1}\Delta U_{\rm circ}=\mu_{\rm circ}\).
Lemma~\ref{pd:root-return} now gives weak convergence in probability.

For the almost sure assertion, assume \eqref{pd:power}.
By \eqref{pd:area-fluct}, \(\|u_N-\E u_N\|_{L^1(D_R)}\to0\) almost
surely for every integer \(R\). Since \(\E\|u_N-U_{\rm circ}\|_{L^1(D_R)}\to0\) also gives \(\E u_N\to U_{\rm circ}\) in local \(L^1\), we get
\(u_N\to U_{\rm circ}\) in local \(L^1\) almost surely, and
Lemma~\ref{pd:root-return} gives almost sure weak convergence of the ESDs \(\mu_N\).
\end{proof}

\section{General profiles and Brown measures}\label{sec:variance-profiles}

For general variance profiles satisfying the assumptions below, we
compare the random logarithmic potential with the Brown potential
through the second-moment pressure.
Section~\ref{sec:applications} applies the resulting convergence
criterion to specific profiles and to the construction of limiting
laws.

\subsection{Comparison and convergence}\label{sec:general-comparison}

Let \(X_N=(\sqrt{s_{ij}}\xi_{ij})\), where the \(\xi_{ij}\) are
independent copies of a random variable \(\xi\) satisfying \eqref{pd:law}, and write
\(\epsilon_N=\max_{i,j}s_{ij}\). The free matrix in
Theorem~\ref{vp:brown} is \(B_N=A_N+C_{S_N}\), formed from the same
profile \(S_N=(s_{ij})\) and the diagonal field \(A_N=\diag(a_i)\).
Here \(C_{S_N}=(\sqrt{s_{ij}}c_{ij})\), and the \(c_{ij}\) are freely
independent standard circular elements in an algebra with a trace
\(\varphi\) satisfying \(\varphi(1)=1\) and \(\varphi(ab)=\varphi(ba)\)
\cite[Lectures~1 and~5]{fp-ns}.

For an \(N\times N\) matrix \(B\) with entries in this algebra, let
\(\mathfrak t_N(B)=N^{-1}\sum_{i=1}^N\varphi(B_{ii})\) be the
normalized matrix trace. For \(\eta>0\), the regularized Brown
potential \(U_{N,\eta}\) and the Brown potential \(U_N\) of \(B_N\) are
\[
 U_{N,\eta}(z)=\tfrac12\mathfrak t_N
 \log\big((B_N-zI)(B_N-zI)^*+\eta^2I\big),\qquad
 U_N(z)=\lim_{\eta\downarrow0}U_{N,\eta}(z).
\]
The logarithm is defined by functional calculus for the positive
operator inside the trace \cite[Lecture~3]{fp-ns}; thus \(U_N(z)\) is
the logarithm of the normalized Fuglede--Kadison determinant of
\(B_N-zI\) \cite{vp-fk,vp-hs}. The Brown measure \(\nu_N\) of
\(B_N\) has logarithmic potential
\(U_N(z)=\int_{\C}\log|z-w|\dd\nu_N(w)\), so
\(\nu_N=(2\pi)^{-1}\Delta U_N\), where \(\Delta\) is the distributional
Laplacian on \(\C\) \cite{vp-brown-original}.
We write \(S=S_N\) and \(A=A_N\) when the dimension is clear.

The comparison applies to determinants and permanents alike, so we
treat both through one polynomial. Let \(\omega_N(\pi)\) be
deterministic coefficients of modulus one, indexed by permutations
\(\pi\in\mathfrak S_N\) of \([N]\), with \(\omega_N(\mathrm{id})=1\),
and define
\[
 F_N(z)=\sum_{\pi\in\mathfrak S_N}\omega_N(\pi)
             \prod_i(zI-A_N-X_N)_{i,\pi(i)}.
\]
The choices \(\omega_N(\pi)=\operatorname{sgn}(\pi)\) and
\(\omega_N(\pi)=1\) give the determinant and the permanent of
\(zI-A_N-X_N\), respectively. Since \(\omega_N(\mathrm{id})=1\), the
polynomial \(F_N\) is monic.
Let \(\mu_N^F\) be the normalized zero measure of \(F_N\), which
gives mass \(1/N\) to each zero counted with multiplicity. Its logarithmic potential is
\[
 u_N^F(z)=\frac1N\log|F_N(z)|
        =\int_{\C}\log|z-w|\dd\mu_N^F(w).
\]
For the determinant coefficients and \(A_N=0\), the function
\(u_N^F\) is the potential \(u_N\) of
Section~\ref{sec:planar-dini}.
As in Section~\ref{sec:planar-dini}, \(D_R=\{z\in\C:|z|<R\}\), \(\mathrm dA\)
is planar Lebesgue measure, and zero measures are compared in the
bounded-Lipschitz distance \(d_{\BL}\) recalled before Lemma~\ref{pd:root-return}.

Set \(q_i=|z-a_i|^2\) and \(\boldsymbol q=(q_1,\ldots,q_N)\).
We call \(q_i\) the diagonal activity, to distinguish it from the
diagonal variance \(s_{ii}\).
Expanding each diagonal factor \(z-a_i-X_{ii}\) of \(F_N(z)\) into
\(z-a_i\) and \(-X_{ii}\), we write \(F_N(z)\) as a sum of terms. In
each term, let \(I\) be the set of rows in which the term uses a random entry. Deterministic factors occur only on the diagonal, so the
permutation fixes every \(j\notin I\) and restricts to a permutation
of \(I\) with factors from \(X_N[I,I]\). Distinct monomials of this
form are orthogonal in \(L^2\) by centering and independence, as in
Lemma~\ref{pd:identity}; since the coefficients \(\omega_N(\pi)\) have modulus one, this yields
\begin{equation}\label{pd:second}
 M_N(z)=\E|F_N(z)|^2
       =\sum_{I\subset[N]}\per(S_N[I,I])\prod_{j\notin I}q_j
       =\per(S_N+\diag q_i).
\end{equation}
Write \(p_N=(2N)^{-1}\log M_N\) for the second-moment pressure. The
term with \(I=\varnothing\) equals \(\prod_i q_i\) and all terms are
nonnegative, so \(p_N\) is finite away from the diagonal centers
\(a_1,\ldots,a_N\); local \(L^1\) estimates may ignore this finite set.

\subsubsection{Comparison with Brown measures}
Write \(d_i=\sum_j s_{ij}\) for the row sums of the profile. We use
the following row-tail condition on these sums:
\begin{equation}\label{vp:row-log-tail}
 \lim_{r\to\infty}\limsup_{N\to\infty}
 \frac1N\sum_i\log(1+d_i/r^2)=0.
\end{equation}

\begin{proposition}[Comparison of zero measures with Brown measures]\label{vp:brown-log-tail}
Let the entry law satisfy \eqref{pd:law}, and let
\(A_N=\diag(a_1,\ldots,a_N)\) be a deterministic diagonal
field with
\(\sup_N N^{-1}\sum_i|a_i|^2<\infty\).
Suppose that \(\epsilon_N=\max_{i,j}s_{ij}\to0\) and that \(S_N\)
satisfies \eqref{vp:row-log-tail}. Then the Brown measures \(\nu_N\)
formed from \(S_N\) and \(A_N\) satisfy
\(\E d_{\BL}(\mu_N^F,\nu_N)\to0\). If moreover \eqref{pd:power}
holds, then \(d_{\BL}(\mu_N^F,\nu_N)\to0\) almost surely under any
coupling across dimensions.
\end{proposition}

The average variance bound in \eqref{vp:average-bounds}
implies \eqref{vp:row-log-tail}. Indeed, if
\(N^{-1}\sum_{i,j}s_{ij}\le C\), then \(\log(1+x)\le x\) gives
\(N^{-1}\sum_i\log(1+d_i/r^2)\le C/r^2\). In particular, uniformly bounded row sums suffice.

\subsubsection{The pressure criterion}
By Proposition~\ref{vp:brown-log-tail}, the existence of a weak ESD
limit reduces to weak convergence of the deterministic measures
\(\nu_N\), which can be decided from the pressure \(p_N\). Without a
diagonal field, \(p_N\) is radial and \(2p_N(e^{t/2})=H_N(t)\), the
scalar pressure of Theorem~\ref{vp:pressure-main}; for a general
field, the criterion uses local \(L^1\) convergence in the plane.

\begin{corollary}[Pressure criterion with a diagonal field]\label{vp:pressure-limit}
Under the assumptions of Proposition~\ref{vp:brown-log-tail}, the ESD of \(A_N+X_N\) converges weakly in probability if and only if the pressure \(p_N\)
converges in \(L^1_{\mathrm{loc}}(\C)\).
In that case the limiting ESD is \(\nu=(2\pi)^{-1}\Delta p\), where
\(p\) is the limit of \(p_N\).
Under \eqref{pd:power}, the same criterion characterizes almost sure
weak convergence under any coupling across dimensions.
When \(A_N=0\), the pointwise criterion and the radial formula of
Theorem~\ref{vp:pressure-main} hold under these same assumptions.
\end{corollary}

\subsection{From the random logarithm to the permanent}

The comparison of \(u_N^F\) with \(p_N\) follows the proof for doubly
stochastic profiles in Section~\ref{sec:planar-dini}, which combines
a fourth-moment estimate with row replacement. Since the fourth-moment estimate now has singularities at the diagonal centers and is no
longer uniform in \(z\), we integrate the error over disks, on which
these singularities are integrable.

\begin{proposition}[Comparison with the second-moment pressure]\label{pd:main}
Let the entry law satisfy \eqref{pd:law}. Suppose that
\(\epsilon_N=\max_{i,j}s_{ij}\to0\) and that \(S_N\) satisfies
\eqref{vp:row-log-tail}. Then for each \(R>0\),
\begin{equation}\label{pd:area}
 \sup_{A_N\ {\rm diagonal}}\E\|u_N^F-p_N\|_{L^1(D_R)}\longrightarrow0.
\end{equation}
Under \eqref{pd:power}, \(\|u_N^F-p_N\|_{L^1(D_R)}\to0\)
almost surely for each prescribed sequence of diagonal fields and
coefficients, under any coupling across dimensions.
\end{proposition}

\begin{proof}
\noindent\textbf{1. The integrated fourth-moment error.}
Fix \(m\), and let \(F_N^{(m)}\) be the polynomial built from the
approximating law \(\eta_m\) of Lemma~\ref{pd:approx}. Since
\(\eta_m\) is centered with unit variance, the computation leading
to \eqref{pd:second} shows that \(F_N^{(m)}\) also has second
moment \(M_N\).
For \(z\notin\{a_1,\ldots,a_N\}\), write
\(r_{N,m}(z)=N^{-1}\log(\E|F_N^{(m)}(z)|^4/M_N(z)^2)\).
We will prove that, for all sufficiently large \(N\) with
\(\epsilon_N>0\),
\begin{equation}\label{pd:profile-moment-bound}
 r_{N,m}(z)\le\frac{C_m}{\log(1/\epsilon_N)}
 \left[1+\frac1N\sum_i\log(1+d_i)
          +\frac1N\sum_i\log(1+|z-a_i|^{-2})\right].
\end{equation}
If \(\epsilon_N=0\), the profile vanishes, so the polynomial is
deterministic and \(r_{N,m}=0\).

We first normalize the rows and columns of the profile:
\begin{equation}\label{vp:moment-normalization}
 d_j^{\mathrm{col}}=\sum_i\frac{s_{ij}}{1+d_i},\qquad
 \widetilde s_{ij}=\frac{s_{ij}}{(1+d_i)(1+d_j^{\mathrm{col}})}.
\end{equation}
The rows and columns of \(\widetilde S_N\) sum to at most one,
its entries are at most \(\epsilon_N\), and
\(\sum_jd_j^{\mathrm{col}}\le N\).
Multiplying entry \((i,j)\) of \(zI-A_N-X_N^{(m)}\) by
\(\big[(1+d_i)(1+d_j^{\mathrm{col}})\big]^{-1/2}\) multiplies every
permutation term by the same constant
\(\prod_i\big[(1+d_i)(1+d_i^{\mathrm{col}})\big]^{-1/2}\), so the
ratio \(\E|F_N^{(m)}(z)|^4/M_N(z)^2\) is unchanged. The rescaled
matrix has entry variances \(\widetilde s_{ij}\) and diagonal
activities
\[
 \widetilde q_i=\frac{|z-a_i|^2}{(1+d_i)(1+d_i^{\mathrm{col}})}.
\]
As in the proof of Proposition~\ref{pd:fourth-ds}, we pass to the real
reference matrix, symmetrize its rows by independent signs and divide
row \(i\) by \(\sqrt{1+\widetilde q_i}\). Entry \((i,j)\) of the
resulting matrix has second moment
\[
 T_{ij}=\frac{\widetilde s_{ij}+\widetilde q_i{\bf1}_{i=j}}
                  {1+\widetilde q_i},
\]
and \(T\) has row sums at most one, column sums at most two and
off-diagonal entries at most \(\epsilon_N\).
For each symmetrized reference entry, let \(k_{ij}\) be its fourth
moment divided by \(T_{ij}^2\), and set \(k_{ij}=1\) when
\(T_{ij}=0\). The entrywise computation in the proof of
Proposition~\ref{pd:fourth-ds}, with \(q\) replaced by
\(\widetilde q_i\) in row \(i\), gives
\[
 k_{ij}\le\max(3,k),\qquad
 \sum_{i,j}(k_{ij}-1)T_{ij}^2\le(2k-1)N\epsilon_N,
 \qquad k=\E|\eta_m|^4.
\]
These are the moment conditions of Lemma~\ref{pd:matchings} with
\(\epsilon=\epsilon_N\), \(K_0=\max(3,k)\) and \(C_0=2k-1\).
The domination \eqref{pd:dominate} applies to \(F_N^{(m)}\) because
\(|\omega_N(\pi)|=1\). Applying \eqref{pd:matching-ratio} to the
real reference matrix and then Lemma~\ref{pd:matchings}, we obtain
\[
 r_{N,m}(z)\le\frac{C_m}{\log(1/\epsilon_N)}
             \left(1-\frac1N\log\per T\right).
\]
The identity permutation gives
\(\per T\ge\prod_iT_{ii}\ge\prod_i\widetilde q_i/(1+\widetilde q_i)\),
which supplies the required lower bound:
\begin{align*}
 -\frac1N\log\per T
 &\le\frac1N\sum_i\log(1+1/\widetilde q_i)\\
 &\le\frac1N\sum_i\log(1+d_i)+1
       +\frac1N\sum_i\log(1+|z-a_i|^{-2}).
\end{align*}
In the second inequality we used
\(\sum_i\log(1+d_i^{\mathrm{col}})\le\sum_id_i^{\mathrm{col}}\le N\).
This proves \eqref{pd:profile-moment-bound}, after adjusting \(C_m\).

It remains to integrate this bound over \(z\). Applying
\(1+d_i\le r^2(1+d_i/r^2)\) with a fixed large \(r\ge1\), we see that
the row-tail assumption implies
\begin{equation}\label{vp:mean-log}
 \sup_N\frac1N\sum_i\log(1+d_i)<\infty.
\end{equation}
The singularities at the diagonal centers are integrable, uniformly
in their position:
\[
 \sup_{a\in\C}\int_{D_R}\log(1+|z-a|^{-2})\dd A(z)
 \le |D_R|\log2+2\pi\log2.
\]
This follows by splitting at \(|z-a|=1\) and using polar coordinates
inside. Thus the bracket in \eqref{pd:profile-moment-bound},
integrated over \(D_R\), is bounded uniformly in \(N\) and \(A_N\),
while \(\log(1/\epsilon_N)\to\infty\); hence
\begin{equation}\label{pd:areamoment}
 \sup_{A_N}\int_{D_R}r_{N,m}(z)\dd A(z)\longrightarrow0
 \qquad\text{for each fixed }m.
\end{equation}

\par\smallskip\noindent\textbf{2. The logarithmic comparison.}
Let \(u_N^{(m)}=N^{-1}\log|F_N^{(m)}|\).
At each \(z\) off the centers, Lemma~\ref{pd:row-concentration}
gives \eqref{pd:gap} for
\(F=F_N^{(m)}(z)\), with \(r=r_{N,m}(z)\). Integrating it over
\(D_R\) and using
\(\int_{D_R}\sqrt{r_{N,m}}\le|D_R|^{1/2}
(\int_{D_R}r_{N,m})^{1/2}\), we obtain
\(\sup_{A_N}\E\|u_N^{(m)}-p_N\|_{L^1(D_R)}\to0\).
To return to the original law, we use the extension stated after Corollary~\ref{pd:row-approx}, which applies because \(F_N(z)\) and
\(F_N^{(m)}(z)\) are affine in each row and, off the centers, have
positive second moment. By Proposition~\ref{pd:transport} and the common Dini
envelope \eqref{pd:envelope}, the replacement error
\(N^{-1}\E\big|\log|F_N(z)|-\log|F_N^{(m)}(z)|\big|\)
tends to zero as \(m\to\infty\), uniformly in \(N\), \(A_N\) and
\(z\) off the centers.
Letting first \(N\to\infty\) and then \(m\to\infty\) proves
\eqref{pd:area}.

\par\smallskip\noindent\textbf{3. Almost sure convergence.}
Under \eqref{pd:power}, \eqref{pd:area-fluct} gives
\(\|u_N^F-\E u_N^F\|_{L^1(D_R)}\to0\) almost surely for every
integer \(R\). Since
\(\|\E u_N^F-p_N\|_{L^1(D_R)}\le\E\|u_N^F-p_N\|_{L^1(D_R)}\to0\) by
the mean comparison, the triangle inequality gives the almost sure
assertion.
\end{proof}

\subsection{The Brown potential and the permanent}

We now compare \(\log\per(S+\diag q_i)\) directly with \(2NU_N(z)\).
The lower bound follows from Gurvits' inequality at weights given by
the solution of the Dyson equation, and the same solution bounds the
mean log weight plus an entropy expression in the edge marginals by
\(2NU_N(z)\). We use this bound in the next subsection to prove the
reverse inequality, up to \(o(N)\) when the activities are bounded
below, and in local \(L^1\) in general.

\begin{lemma}[Dyson bounds for the permanent]\label{vp:flow}
Let \(z\) satisfy \(q_i=|z-a_i|^2>0\) for every \(i\), and let
\(x=(x_{ij})\) be a nonnegative array that vanishes where \(s_{ij}=0\) and whose
outgoing and incoming masses agree at each vertex and are at most one:
\begin{equation}\label{vp:flow-conservation}
 \sum_jx_{ij}=\sum_jx_{ji}=:r_i\le1,\qquad i\in[N].
\end{equation}
Then, with \(m_i=1-r_i\),
\begin{equation}\label{vp:direct-sandwich}
 \sum_{i,j}x_{ij}\log\frac{s_{ij}r_i}{x_{ij}}
 +\sum_i\left[m_i\log\frac{q_i}{m_i}-r_i\right]
 \le2NU_N(z)\le\log\per(S+\diag q_i).
\end{equation}
Terms with zero mass are given their continuous values.
\end{lemma}
\begin{proof}
\noindent\textbf{The Dyson equation and the regularized potential.}
Write \(Y=A+C_S-zI_N\) and fix \(\eta>0\). The positive
Green-function vectors \(u\) and \(v\) are obtained by applying
\(\varphi\) to the diagonal entries of the imaginary part of the
Hermitized resolvent of \(Y\):
\[
 u_i=\eta\,\varphi\!\left([(YY^*+\eta^2I_N)^{-1}]_{ii}\right),\qquad
 v_i=\eta\,\varphi\!\left([(Y^*Y+\eta^2I_N)^{-1}]_{ii}\right).
\]
We write \(u(\eta),v(\eta)\) when the dependence on \(\eta\)
matters, and we define the auxiliary quantities
\[
 g_i:=\eta+(Sv)_i,\qquad \ell_i:=\eta+(S^{\mathsf T}u)_i .
\]
The standard Dyson equation
\cite[Lemmas~4.1--4.2]{vp-ak-brown}
characterizes \(u,v\) as the positive solution of
\begin{equation}\label{vp:dyson}
 u_i=\frac{\ell_i}{q_i+g_i\ell_i},\qquad
 v_i=\frac{g_i}{q_i+g_i\ell_i},\qquad i\in[N].
\end{equation}
This is the equation of the operator-valued circular framework
\cite[Chapter~4]{fp-speicher-op}, with covariance maps
\(w\mapsto Sw\) and \(w\mapsto S^{\mathsf T}w\) obtained by applying
\(\varphi\) to the diagonal entries of \(C_S\diag(w)C_S^*\) and
\(C_S^*\diag(w)C_S\); it holds for every nonnegative profile.

Write \(\langle w\rangle=N^{-1}\sum_iw_i\), and read products of
vectors inside the brackets coordinatewise. Then
\begin{equation}\label{vp:stationary}
 2U_{N,\eta}(z)=\langle\log(\boldsymbol q+g\ell)\rangle
                         -\langle u(Sv)\rangle .
\end{equation}
To verify this identity, differentiate in \(\eta\) at fixed \(z\):
\[
 \frac{\mathrm d}{\mathrm d\eta}(2U_{N,\eta}(z))
 =2\eta\,\mathfrak t_N((YY^*+\eta^2I_N)^{-1})
 =2\langle u\rangle=\langle u+v\rangle .
\]
The last equality holds because \(YY^*\) and \(Y^*Y\) have the
same spectral distribution under \(\mathfrak t_N\). On the right side
of \eqref{vp:stationary}, with primes denoting \(\eta\)-derivatives,
\eqref{vp:dyson} gives
\[
 \begin{aligned}
 \frac{\mathrm d}{\mathrm d\eta}
 \big[\langle\log(\boldsymbol q+g\ell)\rangle-\langle u(Sv)\rangle\big]
 &=\langle ug'+v\ell'-u'(Sv)-u(Sv')\rangle\\
 &=\langle u+v\rangle .
 \end{aligned}
\]
For the second equality, substitute \(g'={\bf1}+Sv'\) and
\(\ell'={\bf1}+S^{\mathsf T}u'\); the terms in \(u'\) and \(v'\) then
cancel. Both sides of \eqref{vp:stationary} therefore have the same
\(\eta\)-derivative, and both are \(2\log\eta+o(1)\) as
\(\eta\to\infty\), since \(u,v=O(\eta^{-1})\) and
\(g,\ell=\eta{\bf1}+O(\eta^{-1})\); so they coincide.

\noindent\textbf{The bound for edge marginals.}
Fix \(x\) satisfying \eqref{vp:flow-conservation}, and let
\(v=v(\eta)\). Since \(x\) conserves mass at every vertex, reweighting the edges by
\(v_j/v_i\) does not change the \(x\)-weighted sum of logarithms:
\begin{equation}\label{vp:tilt-cancellation}
 \sum_{i,j}x_{ij}(\log v_j-\log v_i)
 =\sum_i\log v_i\left(\sum_jx_{ji}-\sum_jx_{ij}\right)=0.
\end{equation}
For nonnegative arrays \((a_j)\) and \((b_j)\) with positive totals
such that \(a_j=0\) wherever \(b_j=0\), we use the log-sum inequality
\begin{equation}\label{vp:log-sum}
 \sum_j a_j\log\frac{a_j}{b_j}
 \ge\left(\sum_j a_j\right)\log\frac{\sum_j a_j}{\sum_j b_j}.
\end{equation}
In each row with \(r_i>0\), we take \(a_j=x_{ij}\) and
\(b_j=s_{ij}v_j/v_i\). Then
\begin{equation}\label{vp:tilted-edge-bound}
 \sum_jx_{ij}\log\frac{s_{ij}v_jr_i}{v_ix_{ij}}
 \le r_i\log\frac{(Sv)_i}{v_i}
 \le r_i\log(q_i+g_i\ell_i),
\end{equation}
since \((Sv)_i\le g_i\) and, by \eqref{vp:dyson},
\(g_i/v_i=q_i+g_i\ell_i\).
Rows with \(r_i=0\) contribute zero. By
\eqref{vp:tilt-cancellation}, the left side of
\eqref{vp:direct-sandwich} is therefore at most
\begin{equation}\label{vp:action-upper}
 \begin{aligned}
 &\sum_i\left[r_i\log(q_i+g_i\ell_i)
                  +m_i\log\frac{q_i}{m_i}-r_i\right]\\
 &\quad\le\sum_i\left[\log(q_i+g_i\ell_i)
                    -\frac{g_i\ell_i}{q_i+g_i\ell_i}\right]
 =2NU_{N,\eta}(z)-\eta\sum_i u_i
 \le2NU_{N,\eta}(z).
 \end{aligned}
\end{equation}
For the first inequality in \eqref{vp:action-upper}, write
\(r_i=1-m_i\) and use
\(m_i\log(q_i/[m_i(q_i+g_i\ell_i)])
\le q_i/(q_i+g_i\ell_i)-m_i\), which follows from \(\log t\le t-1\)
and remains valid for \(m_i=0\) by continuity.
The equality uses \eqref{vp:stationary} and the identity
\(g_i\ell_i/(q_i+g_i\ell_i)=u_i(Sv)_i+\eta u_i\), which follows from
\eqref{vp:dyson} and \(g_i=\eta+(Sv)_i\).
Letting \(\eta\downarrow0\) gives the first inequality in
\eqref{vp:direct-sandwich}.

\noindent\textbf{The permanent lower bound.}
For any \(x\) satisfying \eqref{vp:flow-conservation}, the matrix
\(P=x+\diag m_i\) is doubly stochastic. Gurvits' lower bound
\cite[Theorem~2.2]{vp-bethe}, applied to \(B=S+\diag q_i\), gives
\[
 \log\per B\ge\sum_{i,j}\left[
 P_{ij}\log(B_{ij}/P_{ij})+(1-P_{ij})\log(1-P_{ij})\right].
\]
On the diagonal, \eqref{vp:log-sum} bounds the first term from below by
\(x_{ii}\log(s_{ii}/x_{ii})+m_i\log(q_i/m_i)\).
For the second term, write \(h(t)=(1-t)\log(1-t)\); we use
\(h(t)\ge-t\) off the diagonal and
\(h(m_i+x_{ii})\ge h(m_i)-x_{ii}=r_i\log r_i-x_{ii}\) on the diagonal.
Consequently, \(\log\per B\) is at least the left side of
\eqref{vp:direct-sandwich}, and it remains to choose \(x\) so that
this left side equals \(2NU_N(z)\).

First let \(A=0\), so that every activity equals \(q=|z|^2>0\), and
fix a nonzero irreducible block of \(S\) with index set \(J\). We write
\(n=|J|\) and \(S_J=S[J,J]\), and let \(\rho_J\) be the Perron root of
\(S_J\). The block potential \(U_J\) is defined as \(U_N\), with
\(S_J\) in place of \(S\) and the normalized trace
\(\mathfrak t_J(B)=n^{-1}\sum_{i\in J}\varphi(B_{ii})\) in place of
\(\mathfrak t_N\). By \cite[Theorem~2.2 and Lemmas~4.1--4.3]{vp-chnr1},
as \(\eta\downarrow0\) the
Green-function vectors \(u(\eta),v(\eta)\) of the block and the
auxiliary quantities \(g(\eta),\ell(\eta)\) have limits
\begin{equation}\label{vp:dyson-limit}
 \begin{gathered}
 u^0:=\lim_{\eta\downarrow0}u(\eta),\qquad
 v^0:=\lim_{\eta\downarrow0}v(\eta),\\
 g^0:=\lim_{\eta\downarrow0}g(\eta)=S_Jv^0,\qquad
 \ell^0:=\lim_{\eta\downarrow0}\ell(\eta)=S_J^{\mathsf T}u^0.
 \end{gathered}
\end{equation}
These limits are positive for \(0<q<\rho_J\) and zero for
\(q\ge\rho_J\); they satisfy \(\sum_{i\in J}u_i^0=\sum_{i\in J}v_i^0\)
and
\begin{equation}\label{vp:dyson-zero}
 u_i^0=\frac{\ell_i^0}{q+g_i^0\ell_i^0},\qquad
 v_i^0=\frac{g_i^0}{q+g_i^0\ell_i^0}.
\end{equation}
Letting \(\eta\downarrow0\) in \eqref{vp:stationary} for the block
\(J\) expresses its Brown potential through the limits
\eqref{vp:dyson-limit}:
\begin{equation}\label{vp:zero-brown-potential}
 2nU_J(z)=\sum_{i\in J}
       [\log(q+g_i^0\ell_i^0)-u_i^0g_i^0].
\end{equation}
For \(q<\rho_J\), take
\begin{equation}\label{vp:dyson-candidate}
 x_{ij}=s_{ij}u_i^0v_j^0,\qquad
 m_i=\frac{q}{q+g_i^0\ell_i^0}.
\end{equation}
With this choice, the row and column masses agree:
\begin{equation}\label{vp:dyson-flow}
 \sum_{j\in J}x_{ij}=u_i^0g_i^0
 =\frac{g_i^0\ell_i^0}{q+g_i^0\ell_i^0}
 =v_i^0\ell_i^0=\sum_{j\in J}x_{ji}=:r_i=1-m_i.
\end{equation}
Substituting \(x_{ij}=s_{ij}u_i^0v_j^0\) and using these row
and column sums, we find that the left side of
\eqref{vp:direct-sandwich}, restricted to \(J\), is
\begin{equation}\label{vp:attainment-goal}
 \begin{aligned}
 &\sum_{i\in J}\left[
 r_i\log\frac{r_i}{u_i^0v_i^0}+m_i\log\frac q{m_i}-r_i\right]\\
 &\qquad=\sum_{i\in J}[\log(q+g_i^0\ell_i^0)-r_i]=2nU_J(z),
 \end{aligned}
\end{equation}
since \(r_i/(u_i^0v_i^0)=q/m_i=q+g_i^0\ell_i^0\).
For \(q\ge\rho_J\), take \(x=0\) on \(J\). Then \(m_i=1\) for
\(i\in J\), and the restricted left side and \(2nU_J(z)\) both equal
\(n\log q\), the latter by \eqref{vp:zero-brown-potential} because
the limits vanish.

Now suppose that \(S\) is reducible. Its support graph has an edge
\(i\to j\) whenever \(s_{ij}>0\), and we order the strongly connected
components of this graph so that every edge between different
components points forward. Then
\(C_S-zI\) is block upper triangular, so by
\cite[Proposition~2.24]{vp-hs} the logarithmic Fuglede--Kadison
determinant splits as
\[
 2NU_N(z)=\sum_J2|J|U_J(z).
\]
We define \(x\) on each component as
above and set \(x_{ij}=0\) when \(i\) and \(j\) lie in different
components. A singleton component with \(s_{ii}=0\) receives \(x=0\)
and contributes \(\log q\) to both sides. The combined array
satisfies \eqref{vp:flow-conservation}, and its left side in
\eqref{vp:direct-sandwich} equals \(2NU_N(z)\). Gurvits' bound
therefore proves the second inequality in \eqref{vp:direct-sandwich}
for \(A=0\).

Finally, we remove the assumption \(A=0\). Factor
\(zI-A-C_S=D(I-D^{-1}C_S)\), where \(D=\diag(z-a_i)\) is invertible. By phase invariance of the circular
elements, \(D^{-1}C_S\) may be identified with the free matrix of
the profile \(s'_{ij}=s_{ij}/q_i\). Multiplicativity of the
Fuglede--Kadison determinant \cite{vp-fk,vp-hs} and row scaling
of the permanent give
\begin{equation}\label{vp:field-factor}
 \begin{aligned}
 2NU_N(z)&=\sum_i\log q_i+2NU_{S'}(1),\\
 \log\per(S+\diag q_i)&=\sum_i\log q_i+\log\per(I+S').
 \end{aligned}
\end{equation}
Here \(U_{S'}\) denotes the zero-field Brown potential for \(S'\).
The second inequality for the field \(A\) thus follows from the
zero-field case, applied to \(S'\) at the point \(1\).
\end{proof}

By phase invariance, the potential depends on the diagonal field
only through the activities \(q_i\). Formula \eqref{vp:stationary} also controls the dependence of the
potential on the profile and the activities. We use this control to
remove diagonal variances and to raise small activities.

\begin{corollary}[Stability of the Brown potential]\label{vp:potential-stability}
At fixed \(z\), the Brown potential is nondecreasing in each
variance \(s_{ij}\) and each activity \(q_i>0\). Moreover, if \(U_N'\) is
the potential after the activities are increased to \(q_i'\ge q_i\),
then
\begin{equation}\label{vp:activity-stability}
 0\le2N(U_N'-U_N)\le\sum_i\log(q_i'/q_i).
\end{equation}
\end{corollary}
\begin{proof}
Fix \(\eta>0\). When \eqref{vp:stationary} is differentiated in
\(\log q_i\) or \(s_{ij}\), the derivatives of \(u,v,g,\ell\) cancel
as in its verification, and the remaining terms give
\[
 \frac{\partial(2NU_{N,\eta})}{\partial\log q_i}
 =\frac{q_i}{q_i+g_i\ell_i}\in[0,1],
\]
and
\begin{equation}\label{vp:brown-differential}
 \frac{\partial(2NU_{N,\eta})}{\partial s_{ij}}=u_iv_j\ge0.
\end{equation}
Integrating the first identity in \(\log q_i\) gives monotonicity in
the activities and the bound \eqref{vp:activity-stability} for
\(U_{N,\eta}\); the second identity gives monotonicity in the
variances. These conclusions extend to zero variances by
continuity, and they pass to \(U_N\) as \(\eta\downarrow0\).
\end{proof}

\subsection{Entropy of a weighted permutation}

For the upper bound on the log permanent, we draw a permutation with
probability proportional to its weight; the log permanent is then
its mean log weight plus its entropy.
Lemma~\ref{vp:perm-entropy} bounds the entropy by an expression in
the edge marginals, with an error caused by the requirement that
distinct rows select distinct columns. The first inequality in
\eqref{vp:direct-sandwich} then bounds the mean log weight plus this
expression by \(2NU_N(z)\). We first assume that the activities are
bounded below, so that the identity permutation controls the
relative entropy of the row marginals.

Write \(\Ent\) for Shannon entropy and
\(D(p\Vert q)=\sum_jp_j\log(p_j/q_j)\) for relative entropy.

\begin{proposition}[Permanent comparison for activities bounded below]\label{vp:entropy}
Suppose that \(S_N\) satisfies \eqref{vp:row-log-tail} and that
\(\epsilon_N=\max_{i,j}s_{ij}\to0\), and fix \(q_0>0\). Then,
uniformly over diagonal fields and points \(z\) with
\(q_i=|z-a_i|^2\ge q_0\) for every \(i\),
\[
 0\le N^{-1}\bigl[\log\per(S_N+\diag q_i)-2NU_N(z)\bigr]
 \longrightarrow0 .
\]
\end{proposition}
\begin{lemma}[Entropy bound through edge marginals]\label{vp:perm-entropy}
Let \(\pi\) be a random permutation of \([n]\) with an arbitrary law,
and let \(I_i={\bf1}_{\pi(i)\ne i}\), \(r_i=\E I_i\), \(m_i=1-r_i\) and
\(K_{ij}=\operatorname{Cov}(I_i,I_j)\). We define the edge array \(x\)
by \(x_{ij}=\Pr(\pi(i)=j)\) for \(i\ne j\) and \(x_{ii}=0\). Since
\(\pi\) moves \(i\) if and only if some other vertex is mapped to
\(i\), the row sum and the column sum of \(x\) at vertex \(i\) both
equal \(r_i\). The entropy expression of \(x\) is
\[
 \Ent_{\rm edge}(x)=-\sum_{i\ne j}x_{ij}\log x_{ij}-\sum_i m_i\log m_i
             +\sum_i r_i\log r_i-\sum_i r_i .
\]
Then, for each \(\delta>0\), with \(\|\cdot\|_{\rm F}\) the Frobenius norm,
\begin{equation}\label{vp:entropyineq}
 \Ent(\pi)\le\Ent_{\rm edge}(x)+2n\delta+
       \delta^{-2}\|x\|_{\rm F}\|K\|_{\rm F}
       +\delta^{-1}\sqrt n\,\|x\|_{\rm F}.
\end{equation}
\end{lemma}
\begin{proof}
\textbf{The random reveal order.}
Let \(\mathcal A=\{i:I_i=1\}\) be the set of vertices moved by
\(\pi\). Since \(\mathcal A\) is determined by \(\pi\),
\(\Ent(\pi)=\Ent(\mathcal A)+\Ent(\pi\mid \mathcal A)\). By
subadditivity, the first term is at most
\(-\sum_i(r_i\log r_i+m_i\log m_i)\). To estimate the second term,
we give each row an independent uniform priority in \([0,1]\) and
reveal the images of the rows in \(\mathcal A\) in order of
increasing priority \cite{vp-radh,vp-anari}. The priorities are
independent of \(\pi\), so conditioning on them leaves
\(\Ent(\pi\mid \mathcal A)\) unchanged. Since \(\pi\) fixes every
vertex outside \(\mathcal A\), it maps \(\mathcal A\) onto itself.

\noindent\textbf{The weights of available columns.}
When the image of row \(i\) is revealed, the available columns are
those in \(\mathcal A\) that have not yet been selected. We give such
a column \(j\) the weight \(w_{ij}=x_{ij}/(r_j+\delta)\); the divisor
compensates for the probability \(r_j\) that column \(j\) lies in
\(\mathcal A\). Writing \(W_i\) for the sum of the
available weights, we use \(w_{ij}/W_i\) as a reference distribution
for \(\pi(i)\). Since \(w_{i,\pi(i)}>0\) almost surely, nonnegativity of
relative entropy shows that the conditional entropy of \(\pi(i)\), given
\(\mathcal A\), the priorities and the images revealed earlier, is at
most the conditional mean of \(-\log(w_{i,\pi(i)}/W_i)\).
Now condition on \(\pi\) and on the priority \(t\) of a row
\(i\in\mathcal A\). The column \(\pi(i)\) is always available, and
every other column \(j\in\mathcal A\) is still available with
probability \(1-t\), the probability that its preimage
\(\pi^{-1}(j)\ne i\) has a larger priority. Thus the expected available weight is
\((1-t)\sum_{j\in \mathcal A}w_{ij}+tw_{i,\pi(i)}\).
By concavity of the logarithm,
\(\E[\log W_i\mid I_i=1,t]\le\log\E[W_i\mid I_i=1,t]\).
For \(r_i>0\), set \(\alpha_i=\E[w_{i,\pi(i)}\mid I_i=1]\) and
\(\beta_i=\E[\sum_{j\in \mathcal A}w_{ij}\mid I_i=1]\). Then
\(\E[W_i\mid I_i=1,t]=(1-t)\beta_i+t\alpha_i\), and therefore
\[
 \Ent(\pi\mid \mathcal A)\le-\sum_{i\ne j}x_{ij}\log w_{ij}
       +\sum_{i:r_i>0}r_i\int_0^1\log((1-t)\beta_i+t \alpha_i)\dd t .
\]
The two means have the explicit forms
\[
 \alpha_i=\frac1{r_i}\sum_j\frac{x_{ij}^2}{r_j+\delta},\qquad
 \beta_i=\frac1{r_i}\sum_j\frac{x_{ij}\E(I_iI_j)}{r_j+\delta}\ge \alpha_i .
\]
Here \(\beta_i\ge\alpha_i\) because \(\pi(i)\in\mathcal A\); set
\(\alpha_i=\beta_i=0\) when \(r_i=0\).

\noindent\textbf{The loss from previously selected columns.}
For \(0<a\le b\) and \(h=a/b\),
\[
 \int_0^1\log((1-t)b+ta)\dd t=\log b-1-\frac{h\log h}{1-h},
\]
with the continuous values at the endpoints. We apply this with
\(a=\alpha_i\) and \(b=\beta_i\). The term \(\log b-1\) is the source
of \(-\sum_i r_i\) in \(\Ent_{\rm edge}\), and the correction
\(-h\log h/(1-h)\) is at most \(\sqrt h\), since
\(-2u\log u\le1-u^2\) with \(u=\sqrt h\).
The integral is nondecreasing in \(b\), so increasing \(\beta_i\) to
\(\beta_i+\delta\) gives an upper bound on it and keeps the
denominator of the square-root error at least \(\delta\). As
\(r_i \alpha_i\le\delta^{-1}\sum_jx_{ij}^2\), the Cauchy--Schwarz
inequality gives
\[
 \sum_i r_i\sqrt{\alpha_i/(\beta_i+\delta)}
 \le\delta^{-1}\sum_i\sqrt{r_i\sum_jx_{ij}^2}
 \le\delta^{-1}\sqrt n\,\|x\|_{\rm F}.
\]

\noindent\textbf{The covariance error.}
To bound the main term \(r_i\log(\beta_i+\delta)\), write
\(\E(I_iI_j)=r_ir_j+K_{ij}\). For \(r_i>0\) this gives
\[
 \beta_i=\sum_j\frac{x_{ij}r_j}{r_j+\delta}
       +\frac1{r_i}\sum_j\frac{x_{ij}K_{ij}}{r_j+\delta}.
\]
The first sum is at most \(r_i\), since \(\sum_jx_{ij}=r_i\). The
tangent-line bound for the logarithm at \(r_i+\delta\) therefore shows
\[
 \sum_i r_i\log(\beta_i+\delta)\le\sum_i r_i\log(r_i+\delta)
  +\sum_{i,j}\frac{x_{ij}K_{ij}}{(r_i+\delta)(r_j+\delta)}.
\]
By the Cauchy--Schwarz inequality, the last sum is at most
\(\delta^{-2}\|x\|_{\rm F}\|K\|_{\rm F}\).
Since the column sums of \(x\) are the \(r_j\), we have
\(-\sum x_{ij}\log w_{ij}=-\sum x_{ij}\log x_{ij}
+\sum_j r_j\log(r_j+\delta)\). Combining this identity with the bound
on \(\Ent(\mathcal A)\) and the estimates above proves
\eqref{vp:entropyineq}, the remaining error being
\(2\sum_{i:r_i>0}r_i\log(1+\delta/r_i)\le2n\delta\).
\end{proof}

\begin{proof}[Proof of Proposition~\ref{vp:entropy}]
The lower bound is the second inequality in
\eqref{vp:direct-sandwich}. For the upper bound, the error terms of
Lemma~\ref{vp:perm-entropy} are small when the indicators \(I_i\) of
that lemma are weakly correlated and the edge probabilities are diffuse;
Steps~2 and~3 arrange both. Fix \(0<c<1\). We prove that
\begin{equation}\label{vp:entropy-target}
 \limsup_{N\to\infty}
 \frac{\log\per(S_N+\diag q_i)-2NU_N(z)}N
 \le c\log2,
\end{equation}
uniformly over \(q_i\ge q_0\); letting \(c\downarrow0\) after
\(N\to\infty\) then gives the upper bound in the proposition.

\noindent\textbf{1. The log permanent as mean log weight plus entropy.}
We first separate the diagonal variances from the diagonal
activities, so that a fixed point \(i\) of a permutation contributes
exactly the factor \(q_i\) to its weight.
Let \(S^0\) be \(S\) with its diagonal removed, and let \(U_N^0\)
be the Brown potential of \(A_N+C_{S^0}\). Variance monotonicity in
Corollary~\ref{vp:potential-stability} and the expansion of the
permanent over permutations give
\[
 U_N^0(z)\le U_N(z),\qquad
 \per(S+\diag q_i)\le\per(S^0+\diag q_i)
                           \prod_i(1+s_{ii}/q_i).
\]
The second holds because the diagonal variances affect the weight of
a permutation only at its fixed points, where they multiply \(q_i\)
by \(1+s_{ii}/q_i\); the logarithmic cost is at most
\(\sum_i\log(1+s_{ii}/q_i)\le\sum_is_{ii}/q_i\le N\epsilon_N/q_0=o(N)\).
Since \(S^0\) also satisfies the hypotheses, these two bounds show
that it suffices to prove \eqref{vp:entropy-target} with \(S^0\) in
place of \(S\). From
now on we write \(S\) for \(S^0\) and \(d_i=\sum_j s_{ij}\) for its
row sums. When \(S=0\), both \(\log\per(\diag q_i)\) and
\(2NU_N(z)\) equal \(\sum_i\log q_i\), so we may assume below that
\(\epsilon_N>0\), and, since \(\epsilon_N\to0\), that
\(\epsilon_N<q_0\).

Let \(B=S+\diag q_i\) and \(Z=\per B\), and sample a random
permutation \(\pi\) with probability proportional to its weight:
\begin{equation}\label{vp:weighted-entropy}
 W(\pi)=\prod_i B_{i,\pi(i)},\qquad
 \Pr(\pi)=\frac{W(\pi)}Z,\qquad
 \log Z=\E\log W(\pi)+\Ent(\pi).
\end{equation}
The mean log weight is determined by the edge marginals
\(\Pr(\pi(i)=j)\). The entropy is harder, because distinct vertices
must have distinct images; Lemma~\ref{vp:perm-entropy} bounds it in
terms of the edge marginals, sharply enough to combine with
\eqref{vp:direct-sandwich}.

\noindent\textbf{2. Revealing a small fraction of the indicators.}
The error terms in Lemma~\ref{vp:perm-entropy} involve the
covariances of the indicators \(I_i={\bf1}_{\pi(i)\ne i}\). To make
these covariances small at a low cost in entropy, we reveal the
indicators of randomly sampled vertices, but not their images.
With \(m=\lceil cN\rceil\), we sample \(m\) i.i.d.\ uniform indices,
independent of \(\pi\), and let \(\mathcal R_j\) record the first \(j\)
sampled indices and their indicators. We keep the record
\(\mathcal R=(T,\mathcal R_T)\), where \(T\) is uniform on
\(\{0,\ldots,m-1\}\) and independent of everything else.
Since the indices and \(T\) are independent of \(\pi\), and each
revealed indicator carries at most \(\log2\) of information,
\begin{equation}\label{vp:reveal-cost}
 \Ent(\pi)-\Ent(\pi\mid\mathcal R)\le m\log2
                                      =cN\log2+O(1).
\end{equation}

To see why this record weakens the correlations, write \({\cal I}\)
for conditional mutual information averaged over the record, and
consider the entropy budget \(\sum_i \Ent(I_i\mid\mathcal R_j)\),
which is at most \(N\log2\) for \(j=0\). Revealing the indicator of
the next sampled index, which is uniform on \([N]\), reduces this sum
by \(N^{-1}\sum_{i,\ell}{\cal I}(I_i;I_\ell\mid\mathcal R_j)\). Since
the budget is nonnegative, the total reduction over
the first \(m\) steps is at most \(N\log2\). Averaging over the
uniform choice of \(T\) gives
\[
 \frac1{N^2}\sum_{i,\ell}{\cal I}(I_i;I_\ell\mid\mathcal R)
 \le\frac{\log2}{m}.
\]
Fix a value of the record and indices \(i,\ell\), let \(P\) be the
conditional joint law of \((I_i,I_\ell)\), and let \(Q\) be the product
of its marginals. With
\(K^{\mathcal R}_{i\ell}=\operatorname{Cov}(I_i,I_\ell\mid\mathcal R)\),
the four entries of \(P-Q\) are
\(K^{\mathcal R}_{i\ell},-K^{\mathcal R}_{i\ell},
-K^{\mathcal R}_{i\ell},K^{\mathcal R}_{i\ell}\),
so \(\|P-Q\|_{\rm TV}=\tfrac12\sum|P-Q|=2|K^{\mathcal R}_{i\ell}|\),
and Pinsker's inequality gives
\[
 D(P\Vert Q)\ge2\|P-Q\|_{\rm TV}^2
                  =8(K^{\mathcal R}_{i\ell})^2.
\]
Since \(D(P\Vert Q)\) is the mutual information of \(I_i\) and
\(I_\ell\) under this conditional law, averaging over the records
gives \(\E[(K^{\mathcal R}_{i\ell})^2]\le
\tfrac18{\cal I}(I_i;I_\ell\mid\mathcal R)\).
Summing over \(i,\ell\) and using the preceding information bound yields
\begin{equation}\label{vp:pin}
 \E\|K^{\mathcal R}\|_{\rm F}^2
 \le\frac{N^2\log2}{8m}=O_c(N).
\end{equation}

\noindent\textbf{3. Diffuse conditional edge probabilities.}
Let \(p_i\) be the \(i\)-th row of \(B\) normalized to a probability
vector, and let \(y_i^{\mathcal R}\) be the conditional law of
\(\pi(i)\) given the record:
\[
 p_i=(B_{ij}/(d_i+q_i))_{j=1}^N,\qquad
 y_i^{\mathcal R}=(\Pr(\pi(i)=j\mid\mathcal R))_{j=1}^N.
\]
We first bound the relative entropy of \(y_i^{\mathcal R}\) with
respect to \(p_i\). By subadditivity of entropy,
\(\E\sum_i\Ent(y_i^{\mathcal R})\ge \Ent(\pi\mid\mathcal R)\).
Using \eqref{vp:weighted-entropy} and \eqref{vp:reveal-cost}, we obtain
\begin{align}
 \E\sum_iD(y_i^{\mathcal R}\Vert p_i)
 &=\sum_i\log(d_i+q_i)-\E\log W(\pi)
                            -\E\sum_i\Ent(y_i^{\mathcal R})\notag\\
 &\le\sum_i\log(d_i+q_i)-\log Z+m\log2\notag\\
 &\le\sum_i\log(1+d_i/q_0)+m\log2=O(N).
 \label{vp:KL}
\end{align}
In the last line we used \(Z\ge\prod_iq_i\), which comes from the
identity permutation, together with \(q_i\ge q_0\), \(m\le N\) and
\eqref{vp:mean-log}.

Let \(x^{\mathcal R}\) be the off-diagonal part of
\(y^{\mathcal R}\), that is, \(x_{ij}^{\mathcal R}=y_{ij}^{\mathcal R}\)
for \(i\ne j\) and \(x_{ii}^{\mathcal R}=0\). Let \(\mathcal L\) be
the set of pairs \((i,j)\) with
\(x_{ij}^{\mathcal R}\ge\sqrt{\epsilon_N/q_0}\). Every such pair has
\(i\ne j\), so \(p_{ij}\le\epsilon_N/q_0\), and hence
\(x_{ij}^{\mathcal R}/p_{ij}\ge\sqrt{q_0/\epsilon_N}\) and
\(x_{ij}^{\mathcal R}\log(x_{ij}^{\mathcal R}/p_{ij})
\ge\tfrac12x_{ij}^{\mathcal R}\log(q_0/\epsilon_N)\).
Since \(y\log(y/p)\ge-p/e\) and \(\sum_{i,j}p_{ij}=N\), the terms of
\(\sum_iD(y_i^{\mathcal R}\Vert p_i)\) with \((i,j)\notin\mathcal L\)
add up to at least \(-N/e\). Consequently
\(\tfrac12\log(q_0/\epsilon_N)
\sum_{(i,j)\in\mathcal L}x_{ij}^{\mathcal R}
\le\sum_iD(y_i^{\mathcal R}\Vert p_i)+N/e\).
To bound the Frobenius norm, use \(x^2\le x\) on \(\mathcal L\) and
\(x^2\le\sqrt{\epsilon_N/q_0}\,x\) outside it.
As \(\sum_{i,j}x_{ij}^{\mathcal R}\le N\), taking expectations gives
\begin{equation}\label{vp:diffuse}
 \E\|x^{\mathcal R}\|_{\rm F}^2
 \le N\sqrt{\epsilon_N/q_0}
 +\frac{2\E\sum_iD(y_i^{\mathcal R}\Vert p_i)+2N/e}
        {\log(q_0/\epsilon_N)}=o(N).
\end{equation}
Like the relative-entropy bound \eqref{vp:KL}, this bound is uniform
over \(q_i\ge q_0\).

\noindent\textbf{4. The entropy estimate and the Brown potential.}
We apply Lemma~\ref{vp:perm-entropy} to the conditional law of
\(\pi\) given each value of the record, with covariance matrix
\(K^{\mathcal R}\) and edge array \(x^{\mathcal R}\), and average
over the records. With \(c,\delta>0\) fixed, the Cauchy--Schwarz
inequality, \eqref{vp:pin} and \eqref{vp:diffuse} show that both
norm-error terms in \eqref{vp:entropyineq} are
\(o(N)\) after averaging. Hence
\[
 \Ent(\pi\mid\mathcal R)
 \le\E\Ent_{\rm edge}(x^{\mathcal R})+2N\delta+o(N).
\]
For each value of the record, the conditional mean log weight plus
\(\Ent_{\rm edge}(x^{\mathcal R})\) equals the left side of
\eqref{vp:direct-sandwich} with \(x=x^{\mathcal R}\), because in
\(W(\pi)\) a fixed point \(i\) contributes the factor \(q_i\) and an
image \(i\mapsto j\ne i\) contributes \(s_{ij}\). As noted in
Lemma~\ref{vp:perm-entropy}, \(x^{\mathcal R}\) satisfies
\eqref{vp:flow-conservation}. Since it also vanishes where
\(s_{ij}=0\), Lemma~\ref{vp:flow} bounds this sum by \(2NU_N(z)\).
Combining this bound with \eqref{vp:weighted-entropy} and
\eqref{vp:reveal-cost} gives
\[
 \log Z-2NU_N(z)\le cN\log2+2N\delta+o(N).
\]
By Step~1, restoring the diagonal variances adds \(o(N)\) to this
bound. Dividing by \(N\) and taking \(\limsup_{N\to\infty}\) with
\(c\) and \(\delta\) fixed, we obtain the upper bound
\(c\log2+2\delta\), uniformly over \(q_i\ge q_0\). Letting
\(\delta\downarrow0\) gives \eqref{vp:entropy-target}, and letting
\(c\downarrow0\) completes the proof.
\end{proof}

To remove the lower bound \(q_i\ge q_0\), we replace each activity
\(q_i=|z-a_i|^2\) by \(\max(q_i,t)\) for a fixed threshold \(t>0\).
This changes nothing unless
\(|z-a_i|<\sqrt t\), so the cost is concentrated near the centers,
and integrating it over \(z\) gives a local \(L^1\) comparison
uniform in the diagonal field.

\begin{proposition}[The pressure and the Brown potential]\label{vp:potential-comparison}
Suppose that \(S_N\) satisfies \eqref{vp:row-log-tail} and that
\(\epsilon_N\to0\). Then, for every \(R>0\),
\begin{equation}\label{vp:deterministic}
 \sup_{A_N\ {\rm diagonal}}\|p_N-U_N\|_{L^1(D_R)}\longrightarrow0.
\end{equation}
\end{proposition}
\begin{proof}
\noindent\textbf{Raising small activities.}
Fix \(t>0\) and a point \(z\) off the centers \(a_i\), so that all
\(q_i>0\), and set \(q_i'=\max(q_i,t)\). Since the permanent is
nondecreasing in each activity,
\(\per(S_N+\diag q_i)\le\per(S_N+\diag q_i')\).
Let \(U_N'\) denote the Brown potential at \(z\) with the raised
activities.
Write \(\log^+x=\max\{\log x,0\}\), so that
\(\log(q_i'/q_i)=\log\max\{1,t/q_i\}=\log^+(t/|z-a_i|^2)\). By
Corollary~\ref{vp:potential-stability}, \(U_N'-U_N\) is at most
\((2N)^{-1}\sum_i\log(q_i'/q_i)\).
Since \(q_i'\ge t\), Proposition~\ref{vp:entropy} with \(q_0=t\)
supplies numbers \(e_N(t)\to0\) for fixed \(t\), uniformly over the
raised activities, such that
\[
 \log\per(S_N+\diag q_i')-2NU_N'\le Ne_N(t).
\]
Combining this estimate with the monotonicity of the permanent, the
bound on \(U_N'-U_N\), and the second inequality in
\eqref{vp:direct-sandwich}, which gives \(p_N\ge U_N\), we obtain
\[
 0\le p_N(z)-U_N(z)\le\tfrac12 e_N(t)+
 \frac1{2N}\sum_i\log^+\frac t{|z-a_i|^2}.
\]

\noindent\textbf{Integrating the error.}
Since \(\log^+(t/|z-a|^2)\) vanishes for \(|z-a|\ge\sqrt t\), polar
coordinates about \(a\in\C\) give
\[
 \int_{\C}\log^+(t/|z-a|^2)\dd A(z)
 =2\pi\int_0^{\sqrt t}r\log(t/r^2)\dd r=\pi t .
\]
Integrating the pointwise bound over \(D_R\) and using this identity
for each center, we obtain
\begin{equation}\label{vp:floor}
 \sup_A\|p_N-U_N\|_{L^1(D_R)}
 \le \tfrac{\pi R^2}{2}e_N(t)+\tfrac{\pi}{2}t .
\end{equation}
Letting first \(N\to\infty\) with \(t\) fixed, and then
\(t\downarrow0\), proves \eqref{vp:deterministic}.
\end{proof}

\subsection{Proofs of the comparison and the pressure criterion}

\begin{proof}[Proof of Proposition~\ref{vp:brown-log-tail}]
\textbf{Local comparison.}
Proposition~\ref{pd:main} compares \(u_N^F\) with \(p_N\) in local
\(L^1\), and Proposition~\ref{vp:potential-comparison} compares
\(p_N\) with \(U_N\). By the triangle inequality,
\begin{equation}\label{vp:unconditional}
 \E\|u_N^F-U_N\|_{L^1(D_R)}\longrightarrow0
 \qquad(R>0),
\end{equation}
and, under \eqref{pd:power}, the same local comparison holds almost
surely.

\par\medskip\noindent\textbf{Tightness.}
Lemma~\ref{pd:root-return} also requires the comparison measures
\(\nu_N\) to be tight. We derive this from \eqref{vp:row-log-tail}
and the bound \(\sup_NN^{-1}\sum_i|a_i|^2<\infty\).

Bounding the permanent by the product of its row sums gives
\(p_N(z)\le(2N)^{-1}\sum_i\log(|z-a_i|^2+d_i)\).
We average both sides over the circle \(|z|=r\), on which \(U_N\)
has mean \(\int\log\max\{r,|w|\}\dd\nu_N(w)\) by Jensen's circular-mean formula \eqref{vp:circular-average} and
\(|z-a_i|^2\) has mean \(r^2+|a_i|^2\). Using \(U_N\le p_N\), concavity of the logarithm, and the bound
\(\log(1+x+y)\le\log(1+x)+y\) for
\(x,y\ge0\), we obtain
\begin{equation}\label{vp:comparison-tail}
 \int\log^+(|w|/r)\dd\nu_N(w)
 \le\frac1{2N}\sum_i\log(1+d_i/r^2)
       +\frac1{2Nr^2}\sum_i|a_i|^2.
\end{equation}
By the two assumptions just recalled, both terms on the right tend
to zero as \(r\to\infty\) after \(\limsup_N\). This
logarithmic tail bound, used again in the pressure criterion, implies
tightness, since \(\log^+(|w|/r)\ge\log2\) on \(\{|w|\ge2r\}\).
With \eqref{vp:unconditional} and tightness,
Lemma~\ref{pd:root-return} gives
\(\E d_{\BL}(\mu_N^F,\nu_N)\to0\), and its almost sure version gives
\(d_{\BL}(\mu_N^F,\nu_N)\to0\) almost surely under \eqref{pd:power}.
\end{proof}

\begin{proof}[Proof of Theorem~\ref{vp:brown}]
The variance bound in \eqref{vp:average-bounds} implies
\eqref{vp:row-log-tail}, so Proposition~\ref{vp:brown-log-tail}
applies. With \(\omega_N(\pi)=\operatorname{sgn}(\pi)\), the polynomial
\(F_N\) is \(\det(zI-A_N-X_N)\), so \(\mu_N^F\) is the ESD of
\(A_N+X_N\).
Since \((\nu_N)\) is tight, approximation on compact sets extends
the bounded-Lipschitz comparison to bounded continuous test
functions, which gives the stated conclusions.
\end{proof}

\begin{proof}[Proof of Corollary~\ref{vp:pressure-limit}]
By Proposition~\ref{vp:brown-log-tail}, weak ESD convergence is
equivalent to weak convergence of the tight deterministic measures
\(\nu_N\). Proposition~\ref{vp:potential-comparison} gives \(p_N-U_N\to0\) in local \(L^1\), and the proof of Proposition~\ref{vp:brown-log-tail} gives the logarithmic tail bound
\[
 \lim_{r\to\infty}\limsup_{N\to\infty}
 \int\log^+(|w|/r)\dd\nu_N(w)=0.
\]
Thus \(\log^+|w|\) is uniformly integrable with respect to the
\(\nu_N\), and by standard potential theory
\cite[Lemma~2.14 and Proposition~2.17]{vp-log-potentials},
\(U_N\to U_\nu\) in local \(L^1\) whenever
\(\nu_N\Rightarrow\nu\), where
\(U_\nu(z)=\int\log|z-w|\dd\nu(w)\).
Together with \(p_N-U_N\to0\), this gives \(p_N\to U_\nu\) in local
\(L^1\), proving necessity.

Conversely, if \(p_N\to p\) in local \(L^1\), then \(U_N\to p\), and
taking distributional Laplacians gives
\(\nu_N\to(2\pi)^{-1}\Delta p\) in distributions. By tightness, this
limit is a probability measure and is also the weak limit of the
\(\nu_N\), which proves sufficiency and identifies the limit. Under
\eqref{pd:power}, the almost sure comparison in
Proposition~\ref{vp:brown-log-tail} gives the same equivalence for
almost sure convergence.
The final assertion of the corollary concerns the case \(A_N=0\); we
prove it below together with Theorem~\ref{vp:pressure-main}.
\end{proof}

\begin{proof}[Proof of Theorem~\ref{vp:pressure-main}]
Throughout this proof \(A_N=0\). By
Proposition~\ref{vp:brown-log-tail}, we may replace ESD convergence
by weak convergence of the tight deterministic measures \(\nu_N\).
Since \(A_N=0\) and the circular elements are rotation invariant, these measures are rotation invariant, and hence so are their weak limits. Recall that \(H_N(t)=2p_N(e^{t/2})\). The expansion
\(\per(S_N+e^tI)=\sum_{k=0}^N b_ke^{kt}\), with \(b_k\ge0\),
shows that \(H_N\) is convex and \(1\)-Lipschitz.

\noindent\textbf{Only if: ESD convergence implies pressure convergence.}
Suppose that the ESD converges weakly in probability to \(\nu\).
Then \(\nu_N\Rightarrow\nu\), and the proof of
Corollary~\ref{vp:pressure-limit} gives
\(p_N\to U_\nu\) in local \(L^1\).

To identify the limit of \(H_N\), we compute \(U_\nu\) on circles.
By rotation invariance of \(\nu\), \(U_\nu(r)\) equals the average
of \(U_\nu\) over the circle \(|z|=r\). Jensen's circular-mean
formula \cite[(3.13)]{vp-log-potentials} states that, for \(r>0\),
\begin{equation}\label{vp:circular-average}
 \frac1{2\pi}\int_0^{2\pi}\log|re^{i\theta}-w|\dd\theta
 =\log\max\{r,|w|\}.
\end{equation}
Integrating \eqref{vp:circular-average} against \(\mathrm d\nu(w)\) gives
\(U_\nu(r)=\int\log\max\{r,|w|\}\dd\nu(w)\).
With \(r=e^{t/2}\) we have
\(2\log\max\{r,|w|\}=\max\{t,\log|w|^2\}\), so we define
\begin{equation}\label{vp:radial-pressure-limit}
 H(t):=2U_\nu(e^{t/2})
 =\int\max\{t,\log|w|^2\}\dd\nu(w),
\end{equation}
with \(\log0=-\infty\).

We now prove that \(H_N(t)\to H(t)\) for every \(t\).
Both \(H_N\) and \(H\) are \(1\)-Lipschitz, the latter because each
\(t\mapsto\max\{t,\log|w|^2\}\) is. Suppose that \(|H_N(t_0)-H(t_0)|\ge\delta\) for some
\(\delta>0\) and infinitely many \(N\). For each such \(N\), the
Lipschitz bounds give \(|H_N(t)-H(t)|\ge\delta/2\) whenever
\(|t-t_0|\le\delta/4\). Since \(2p_N(z)=H_N(2\log|z|)\) and, by
rotation invariance, \(2U_\nu(z)=H(2\log|z|)\), it follows that
\(|p_N(z)-U_\nu(z)|\ge\delta/4\) on the ring where
\(|2\log|z|-t_0|\le\delta/4\). This ring does not depend on \(N\)
and has positive area, which contradicts the local \(L^1\)
convergence \(p_N\to U_\nu\). Hence
\(H_N(t)\to H(t)\) for every \(t\), which proves the only-if
direction.

\begin{samepage}
\noindent\textbf{The limiting pressure determines the measure.}
The function \(H\) is finite and convex, so its one-sided
derivatives exist at every \(t\). For each \(w\), the right
difference quotients of \(t\mapsto\max\{t,\log|w|^2\}\) lie in
\([0,1]\) and converge to \({\bf1}_{\{|w|^2\le e^t\}}\). By bounded
convergence,
\begin{align}
 H'_+(t)
 &=\lim_{h\downarrow0}\int
 \frac{\max\{t+h,\log|w|^2\}-\max\{t,\log|w|^2\}}h\dd\nu(w)
 \notag\\
 &=\int {\bf1}_{\{|w|^2\le e^t\}}\dd\nu(w)
 =\nu\{|w|\le e^{t/2}\}.\label{vp:radial-pressure-slope}
\end{align}
These closed-disk masses determine the rotation-invariant measure
\(\nu\); its mass at the origin is obtained by letting
\(t\to-\infty\).

\end{samepage}

\noindent\textbf{If: pressure convergence implies ESD convergence.}
Suppose that \(H_N(t)\to H(t)\) for every \(t\), and let \(\nu'\) be
a subsequential weak limit of \((\nu_N)\). The argument of the only-if direction, applied along
this subsequence, shows that the pressure of \(\nu'\), defined as in
\eqref{vp:radial-pressure-limit}, equals \(H\). By
\eqref{vp:radial-pressure-slope} and rotation invariance, \(\nu'\) is
determined by \(H\), so all subsequential limits coincide. Tightness and Proposition~\ref{vp:brown-log-tail} then
give ESD convergence in probability.

Under \eqref{pd:power}, sufficiency for almost sure convergence
follows from the almost sure comparison in
Proposition~\ref{vp:brown-log-tail}; necessity holds because almost
sure convergence implies convergence in probability.
\end{proof}

\begin{proof}[Proof of Corollary~\ref{vp:circular-pressure}]
The proof of Theorem~\ref{vp:pressure-main} identifies the limiting
pressure as \(H(\log q)=\int\log\max\{q,|w|^2\}\dd\nu(w)\).
For the uniform law on the unit disk, \(|w|^2\) is uniform on
\([0,1]\), so this integral equals \(q-1\) for \(0<q\le1\) and
\(\log q\) for \(q>1\). Conversely, if the limiting pressure has this form, then for every \(r>0\) we have \(H'_+(2\log r)=\min(r^2,1)\), the mass of the closed disk of radius \(r\) under \(\mu_{\rm circ}\).
Theorem~\ref{vp:pressure-main} now gives the stated
equivalence in both modes of convergence.
\end{proof}

\section{Applications}\label{sec:applications}

We give three applications of the Brown-measure comparison. For
profiles with bounded row and column sums, tree weights give a criterion for convergence of the ESD, and under a fixed bound on these
sums we characterize the radial laws that arise as limits. Finally,
we treat the zeros of the permanental polynomial.

\subsection{Tree-weight criterion}\label{sec:tree-limit}

By the pressure criterion, convergence of the ESD is a
deterministic question about the variance profile. When the row and
column sums are uniformly bounded, tree weights yield a sufficient condition that can be checked without evaluating permanents. For a finite tree \(T\) with
vertex set \(V(T)\) and edge set \(E(T)\),
and a symmetric \(N\times N\) profile \(S'=(s'_{ij})\), define the tree weight
\[
 t_N(T;S')=\frac1N\sum_{\phi:V(T)\to[N]}
 \prod_{\{v,w\}\in E(T)}s'_{\phi(v)\phi(w)} .
\]
The labeling \(\phi\) need not be injective. The tree is fixed as
\(N\) grows, and the criterion below requires no uniformity over
trees of growing size.

For example, write \(d_i'=\sum_j s'_{ij}\) for the row sums. Summing
first over the labels of the leaves, one finds that a single edge has
weight \(N^{-1}\sum_i d_i'\), a path with two edges has weight
\(N^{-1}\sum_i(d_i')^2\), and, more generally, a star with
\(k\) edges has weight \(N^{-1}\sum_i(d_i')^k\). If \(S'\) is doubly
stochastic, removing the leaves one at a time gives
\(t_N(T;S')=1\) for every tree.

\begin{samepage}
\begin{corollary}[Tree-weight criterion]\label{vp:tree-main}
Assume the hypotheses of Theorem~\ref{vp:brown} with \(A_N=0\), and
let the row and column sums of \(S_N\) be uniformly bounded. If
\(t_N(T;(S_N+S_N^{\mathsf T})/2)\) converges for every finite tree
\(T\) and the profile is asymptotically symmetric, in the sense that
\begin{equation}\label{vp:asymptotic-symmetry}
 \frac1N\sum_{i,j=1}^N |s_{ij}-s_{ji}|\longrightarrow0,
\end{equation}
then the ESD converges weakly in probability to a rotation-invariant
law determined by the limiting tree weights. Under
\eqref{pd:power}, the convergence holds almost surely under any
coupling across dimensions.
\end{corollary}
\end{samepage}

How do the tree weights determine the limit? Through the pressure.
Under the hypotheses of the corollary, the scalar pressure \(H_N(t)=N^{-1}\log\per(S_N+e^tI)\) converges to a limit \(H(t)\), and Theorem~\ref{vp:pressure-main} recovers the squared-radius distribution function \(F(q)=H'_+(\log q)\) of the limiting law.

Zhu \cite[Theorem~3.2]{vp-zhu} uses the same tree weights to
identify spectral moments of Hermitian Wigner-type matrices whose
normalized variances are bounded by \(C/N\). Here the tree weights
determine a limiting
Brown measure for sparse non-Hermitian matrices whose maximal
variance tends to zero.

\begin{proof}[Proof of Corollary~\ref{vp:tree-main}]
By Corollary~\ref{vp:pressure-limit}, it suffices to prove that
\(p_N\) converges in \(L^1(D_R)\) for every \(R>0\). Steps~1 and~2 treat a symmetric profile
\(S_N\), and Step~3 reduces the asymptotically symmetric case to
the symmetric one. Let \(d_*\ge1\) be a uniform bound for the row
and column sums. Write \(S=S_N\) and \(C_S=(\sqrt{s_{ij}}c_{ij})\),
and let \(\mathfrak t_N\) be the normalized trace from
Section~\ref{sec:variance-profiles}.

\par\medskip\noindent\textbf{1. Tree limits control the logarithm after regularization.}
We begin with traces of polynomials in \(C_S\) and \(C_S^*\).
For example, with \(d_i=\sum_j s_{ij}\),
\[
 \mathfrak t_N(C_SC_S^*)=\frac1N\sum_i d_i,\qquad
 \mathfrak t_N((C_SC_S^*)^2)=\frac2N\sum_i d_i^2.
\]
These are the weight of the tree with one edge and twice the weight
of the path with two edges, respectively. For a general product of
factors \(C_S\) and \(C_S^*\), the free Wick rule
\cite[Lecture~11]{fp-ns} expresses the trace as a sum over noncrossing pairings. Each contributing pairing yields a tree. An
innermost pair contributes a factor \(s_{ij}\), and summing over its
intermediate index removes a leaf. Since \(S\) is symmetric, edge
orientations do not matter, and the resulting index sums are
exactly the tree weights defined above. Hence the trace of every polynomial in \(C_S\) and
\(C_S^*\) converges to a limit determined by the limits of the tree
weights.

To pass to the regularized logarithm, we bound the norm of
\(C_S\). Summing over the leaves bounds each tree weight with \(k\)
edges by \(d_*^k\), and there are at most \(4^k\) noncrossing
pairings, so
\(\mathfrak t_N((C_SC_S^*)^k)\le(4d_*)^k\); taking
\(2k\)-th roots and letting \(k\to\infty\) gives
\(\|C_S\|\le2\sqrt{d_*}\).
Hence, for each fixed \(z\), the spectrum of
\((C_S-zI)(C_S-zI)^*\) lies in a compact interval that does not
depend on \(N\). For fixed \(\eta>0\), the function
\(x\mapsto\log(x+\eta^2)\) is a uniform limit of polynomials on this
interval, so the convergence of the polynomial traces implies
convergence of
\[
 U_{N,\eta}(z)=\tfrac12\mathfrak t_N
 \log\big((C_S-zI)(C_S-zI)^*+\eta^2I\big).
\]
The norm bound and \(\eta>0\) make \(U_{N,\eta}\) bounded on \(D_R\)
uniformly in \(N\), so dominated convergence gives convergence in
\(L^1(D_R)\).

\par\medskip\noindent\textbf{2. Symmetry lets us remove the regularization uniformly.}
To pass to the Brown potentials
\(U_N=\lim_{\eta\downarrow0}U_{N,\eta}\), we bound the
regularization error uniformly in \(N\).
Since \(S\) is symmetric, exchanging \(u\) and \(v\) leaves the
Dyson system \eqref{vp:dyson} unchanged, and uniqueness of its
positive solution gives \(u=v\) and \(g=\ell\). For \(z\ne0\), the
equations therefore reduce to
\[
 u_i=\frac{g_i}{|z|^2+g_i^2}\le\frac1{2|z|},\qquad
 g_i=\eta+\sum_j s_{ij}u_j.
\]
As \(\partial_\eta U_{N,\eta}(z)=N^{-1}\sum_i u_i\), integrating
first in \(\eta\) and then over \(z\) gives
\[
 0\le U_{N,\eta}(z)-U_N(z)\le\frac{\eta}{2|z|},\qquad
 \|U_{N,\eta}-U_N\|_{L^1(D_R)}\le\pi R\eta,
\]
where we used \(\int_{D_R}|z|^{-1}\dd A(z)=2\pi R\).
Consequently, for any two dimensions \(N\) and \(N'\),
\[
 \|U_N-U_{N'}\|_{L^1(D_R)}
 \le2\pi R\eta+
 \|U_{N,\eta}-U_{N',\eta}\|_{L^1(D_R)}.
\]
The first term is small once \(\eta\) is small, and for fixed
\(\eta\) the last term tends to zero by Step~1. Hence \((U_N)\) is
Cauchy in local \(L^1\).
By Proposition~\ref{vp:potential-comparison}, \(p_N-U_N\to0\) in the
same space, so for symmetric profiles \(p_N\) converges as well, to a
limit determined by the tree weights.

\par\medskip\noindent\textbf{3. Two symmetric profiles enclose the given one.}
For the given asymptotically symmetric profile, set
\[
 S_N^-=(\min(s_{ij},s_{ji})),\qquad
 S_N^+=(\max(s_{ij},s_{ji})).
\]
These profiles are symmetric, their row sums are at most \(2d_*\),
and, with \(\bar S_N=(S_N+S_N^{\mathsf T})/2\), they satisfy
\(S_N^-\le S_N\le S_N^+\) and \(S_N^-\le\bar S_N\le S_N^+\)
entrywise.
To see that \(S_N^-\) and \(S_N^+\) have the same tree limits as
\(\bar S_N\), fix a tree with \(k\ge1\) edges and change the edge
factors one at a time. Removing the edge being changed splits the
tree into two subtrees rooted at its endpoints; with the two endpoint
labels fixed, summing over their leaves bounds their combined
contribution by \((2d_*)^{k-1}\). Summing over the \(k\) edges gives
\[
 |t_N(T;S_N^\pm)-t_N(T;\bar S_N)|
 \le\frac{k(2d_*)^{k-1}}N\sum_{i,j}|s_{ij}-s_{ji}|
 \longrightarrow0.
\]
Since the limits in Steps~1--2 depend only on the tree limits, the
pressures \(p_N^-\) and \(p_N^+\) converge in local \(L^1\) to the
same limit \(p\).
The permanent is monotone in the entries, so
\(p_N^-\le p_N\le p_N^+\), and hence
\[
 |p_N-p|\le|p_N^--p|+|p_N^+-p|.
\]
Integrating over \(D_R\) gives the local \(L^1\) convergence of
\(p_N\).
Corollary~\ref{vp:pressure-limit} now gives convergence of the ESD
in probability, and almost surely under \eqref{pd:power}.
Its proof also gives weak convergence of the Brown measures, and
these measures and their limit are rotation invariant because
\(A_N=0\).
\end{proof}

\begin{example}[A weighted periodic band]
Let the integer bandwidths satisfy \(W_N\to\infty\) and
\(W_N=o(N)\), and let \(b\) be a
bounded nonnegative measurable function of period one.
We discretize \(b\) by its cell averages and set
\begin{equation}\label{vp:weighted-band-profile}
 b_{N,i}=N\int_{(i-1)/N}^{i/N}b(x)\dd x,\qquad
 s_{ij}=\frac{b_{N,i}b_{N,j}}{2W_N+1}
 \mathbf1_{\{d_N(i,j)\le W_N\}},
\end{equation}
where \(d_N(i,j)\) is the distance between \(i\) and \(j\) on the
cycle of \(N\) sites.
This profile is symmetric, its row and column sums are at most
\(\|b\|_\infty^2\), and its maximal variance is \(O(W_N^{-1})\).
For every fixed tree \(T\) with \(k\) edges,
\begin{equation}\label{vp:weighted-band-trees}
 t_N(T;S_N)\longrightarrow\int_0^1 b(x)^{2k}\dd x.
\end{equation}

The limiting ESD is the average over \(x\in[0,1]\) of the uniform law
on the disk of radius \(b(x)\). When \(b(x)=0\), this law is the point
mass at zero.
Here \(b\) need not be continuous; if it is Riemann integrable, the
cell averages may be replaced by the sampled values \(b(i/N)\).
\end{example}
\begin{proof}
The condition \(W_N\to\infty\) gives
\(\max_{i,j}s_{ij}=O(W_N^{-1})\to0\), as required by
Corollary~\ref{vp:tree-main}. Since \(W_N=o(N)\), every labeling of
a fixed tree with \(k\) edges that contributes to the tree weight has
all its labels within normalized distance \(kW_N/N\to0\) of the label of a fixed vertex.

\par\medskip\noindent\textbf{1. Continuous weights.}
For large \(N\), take the fixed vertex as the root and label the tree at random. The root receives a
uniform label, and each child independently receives a uniform label
among the \(2W_N+1\) sites within distance \(W_N\) of its parent's
label. Write \(I_v\) for the label of \(v\). Each admissible labeling
has probability \(1/[N(2W_N+1)^k]\), so
\[
 t_N(T;S_N)=
 \E\prod_{\{v,w\}\in E(T)}b_{N,I_v}b_{N,I_w}.
\]
If \(b\) is continuous, uniform continuity and
\((kW_N+1)/N\to0\) allow us to replace every factor by
\(b(I_{\rm root}/N)\) with an error uniform in the labels.
Since each edge supplies two factors,
\[
 t_N(T;S_N)=\frac1N\sum_{i=1}^N b(i/N)^{2k}+o(1)
 \longrightarrow\int_0^1 b(x)^{2k}\dd x.
\]

\par\medskip\noindent\textbf{2. Measurable weights.}
Let \(b_*\ge1\) be an upper bound for \(b\), and approximate \(b\) in
\(L^1\) by a continuous function \(\widetilde b\) of period one with
\(0\le \widetilde b\le b_*\). The cell averages
\(\widetilde b_{N,i}\) of \(\widetilde b\) satisfy
\(N^{-1}\sum_i|b_{N,i}-\widetilde b_{N,i}|\le\|b-\widetilde b\|_{L^1}\).
Each label \(I_v\) is uniformly distributed, because a step to a
uniform site in the periodic band preserves the uniform law.
Replacing the \(2k\) factors one at a time therefore changes the
tree weight by at most
\(2kb_*^{2k-1}\|b-\widetilde b\|_{L^1}\), uniformly in \(N\).
The target integral changes by at most the same amount, so applying
Step~1 to \(\widetilde b\) and letting \(\|b-\widetilde b\|_{L^1}\to0\)
proves \eqref{vp:weighted-band-trees}.
If \(b\) is Riemann integrable, then
\(N^{-1}\sum_i|b(i/N)-b_{N,i}|\to0\), and the same replacement
argument shows that the sampled values may be used instead of the
cell averages.

\par\medskip\noindent\textbf{3. Limiting measure.}
Let \(c\) be a standard circular element in an algebra with trace
\(\varphi\). The model \(b(x)c\), averaged over \(x\in[0,1]\), assigns
to each noncommutative polynomial \(P\) the trace
\[
 \int_0^1\varphi\bigl(P(b(x)c,(b(x)c)^*)\bigr)\dd x.
\]
By the free Wick rule and \eqref{vp:weighted-band-trees}, the traces
of polynomials in \(C_{S_N}\) and \(C_{S_N}^*\) converge to these
averaged traces. The regularization argument in the
proof of Corollary~\ref{vp:tree-main} therefore identifies the
limiting ESD with the Brown measure of this averaged model.
For each \(x\), the Brown measure of \(b(x)c\) is the uniform law on
the disk of radius \(b(x)\) \cite[Example~5.3]{ep-bl}, or
\(\delta_0\) when \(b(x)=0\), and the claimed mixture follows.\qedhere
\end{proof}

\subsection{Realizable radial limits}\label{vp:radial-design}

Tree weights give the limit for a prescribed sequence of profiles.
Conversely, under a fixed bound on the row and column sums, we
characterize the ESD limits that can be realized by choosing the
profiles. Condition (2) below follows from concavity of the radial Brown potential in the squared radius. For
sufficiency, we arrange the profile as a long directed cycle of
layers and encode the squared-radius quantiles of the target law in the column weights of a single block.

\begin{theorem}[Realizable radial limits]\label{vp:radial-inverse}
Fix \(d_*>0\) and an entry law satisfying \eqref{pd:law}, and let
\(\nu\) be a probability measure on \(\C\). Consider variance
profiles in the model of Theorem~\ref{vp:brown} with \(A_N=0\), \(\epsilon_N\to0\) and row and column sums at most \(d_*\). A sequence of such
profiles can be chosen so that the ESD converges weakly in
probability to \(\nu\) if and only if the following conditions hold.
\begin{enumerate}
\item \(\nu\) is rotation invariant and supported on the closed disk
of radius \(\sqrt{d_*}\).
\item The squared-radius distribution function
\(F(q)=\nu\{z:|z|^2\le q\}\) of \(\nu\) satisfies
\begin{equation}\label{vp:radial-ratio}
 q\longmapsto F(q)/q\quad\text{is nonincreasing on }(0,\infty).
\end{equation}
\end{enumerate}
Moreover, for any sequence of such profiles,
every subsequential limit in probability of the ESD satisfies
\textup{(1)} and \textup{(2)}.
Under \eqref{pd:power}, the profiles can be chosen so that the ESD also converges almost surely under any coupling across dimensions.
\end{theorem}

We prove sufficiency first, by the construction below, and
necessity after Lemma~\ref{vp:radial-potential}.

\begin{figure}[htbp]
\centering
\definecolor{figblue}{HTML}{1E61A5}
\definecolor{figred}{HTML}{CF5246}
\begin{tikzpicture}[
  lab/.style={font=\normalfont\footnotesize},
  layer/.style={draw,line width=0.8pt,minimum width=6.4mm,minimum height=5.6mm,
    inner sep=1pt,font=\normalfont\footnotesize},
  arr/.style={-{Latex[length=1.9mm,width=1.5mm]},line width=0.8pt},
  eig/.style={figblue,opacity=0.55}
]
\node[layer] (n1) at (0,1.45) {\(1\)};
\node[layer] (n2) at (1.55,0.55) {\(2\)};
\node[layer] (n3) at (1.3,-1.05) {\(3\)};
\node[layer] (nL) at (-1.3,-1.05) {\(L\)};
\node[layer] (nM) at (-1.55,0.55) {\(L{+}1\)};
\draw[arr] (n1) to[bend left=22] node[lab,above right=-1pt] {\(1/m\)} (n2);
\draw[arr] (n2) to[bend left=18] node[lab,right] {\(1/m\)} (n3);
\draw[line width=0.8pt,dotted] (n3) to[bend left=28] (nL);
\draw[arr,figred,line width=1.6pt] (nL) to[bend left=18]
  node[lab,left,text=figred] {\(w_j/m\)} (nM);
\draw[arr] (nM) to[bend left=22] node[lab,above left=-1pt] {\(1/m\)} (n1);
\begin{scope}[xshift=5.6cm,yshift=0.15cm]
\input{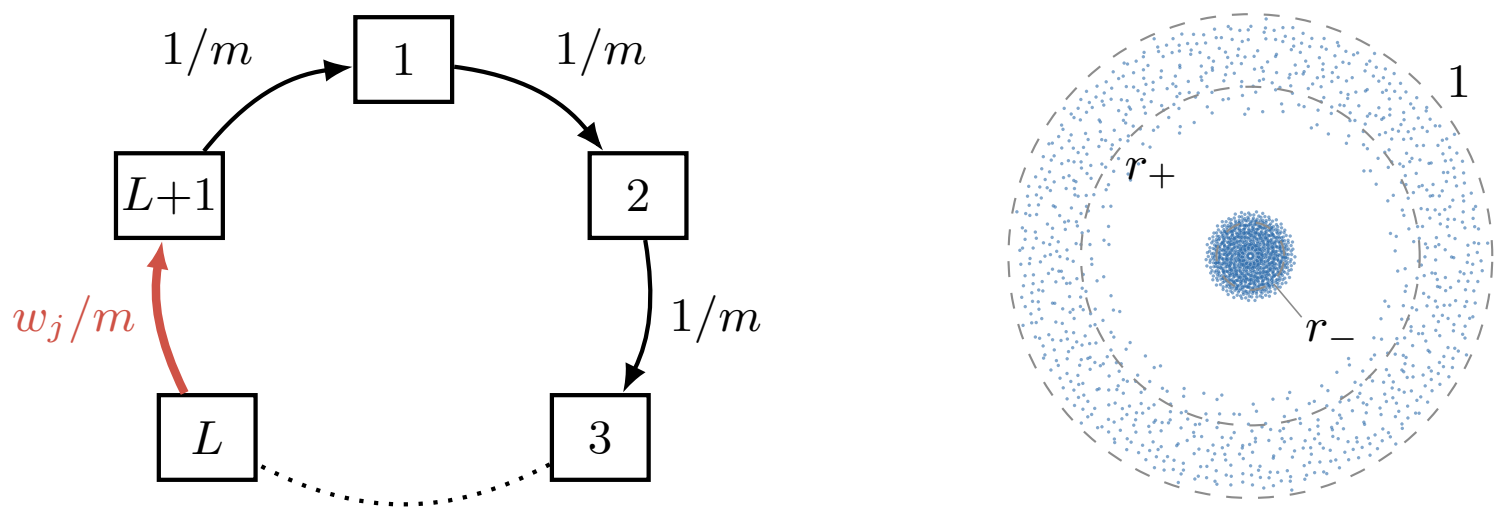}
\foreach \r in {0.14,0.70,1}
  \draw[black!45,line width=0.4pt,dashed] (0,0) circle (\r*1.6);
\draw[black!45,line width=0.3pt] (-50:0.15*1.6) -- (-50:0.33*1.6)
  node[lab,anchor=north west,inner sep=0.5pt,text=black] {\(r_-\)};
\node[lab] at (140:0.53*1.6) {\(r_+\)};
\node[lab] at (40:1.6+0.2) {\(1\)};
\end{scope}
\end{tikzpicture}
\caption{Left, the cycle of layers in the proof of
Theorem~\ref{vp:radial-inverse}, with the variances on each arrow.
Right, the eigenvalues of one sample in
Example~\ref{vp:annular-example} with \(r_-=0.14\), \(r_+=0.70\),
\(L=12\), \(m=150\) and complex Gaussian entries. The dashed circles
have radii \(r_-\), \(r_+\) and \(1\).}
\label{vp:quantile-cycle}
\end{figure}

For \(d_*=1\), the construction starts from a doubly stochastic profile
whose \(L+1\) layers of size \(m\) form a directed cycle, so the ESD tends to the circular law. A given target law is then
realized by changing the variances in a single \(m\times m\) block.
If the matrices before and after this change are
built from the same variables \(\xi_{ij}\), their difference is
supported in that block, so its rank is at most
\(m=N/(L+1)\), which is \(o(N)\) as \(L\to\infty\). Yet this change can turn the limiting support from a disk into a disk
and a separate annulus, as in Figure~\ref{vp:quantile-cycle}. For Hermitian matrices, by contrast, a perturbation of rank \(o(N)\)
does not change the limiting ESD, by the rank inequality
\cite[Lemma~3.8]{vp-zhu}.

\begin{proof}[Sufficiency in Theorem~\ref{vp:radial-inverse}]
Dividing \(S_N\) by \(d_*\) divides the matrix, and hence its
eigenvalues, by \(\sqrt{d_*}\), and it changes the squared-radius
distribution function to \(q\mapsto F(d_*q)\). This change preserves
\eqref{vp:radial-ratio}, so we may take \(d_*=1\).

\par\medskip\noindent\textbf{1. The target pressure and the profile.}
Let \(F(q)=\nu\{z:|z|^2\le q\}\), and define the corresponding
squared-radius quantile by
\[
 Q(p)=\inf\{q\ge0:F(q)\ge p\},\qquad 0<p\le1.
\]
By \eqref{vp:radial-pressure-limit}, the limiting pressure of
\(\nu\) is
\begin{equation}\label{vp:target-quantile-pressure}
 H(t)=\int_0^1\max\{t,\log Q(p)\}\dd p.
\end{equation}
It therefore suffices to construct profiles whose pressures converge
to this function.
Since \(F(1)=1\) by the support condition, \eqref{vp:radial-ratio}
gives \(F(q)\ge q\) on \((0,1)\), and hence \(Q(p)\le p\).
Condition \eqref{vp:radial-ratio} also makes \(F\) continuous on
\((0,\infty)\), so \(F(Q(p))=p\) whenever \(Q(p)>0\). It follows that
the ratio \(Q(p)/p\) is nondecreasing, including on the initial
interval where \(Q\) vanishes.

Let \(L,m\to\infty\) be integers. The index set consists of
\(L+1\) layers of \(m\) vertices each, indexed by pairs \((h,i)\) with
\(1\le h\le L+1\) and \(1\le i\le m\), so that \(N=(L+1)m\).
We take the equally spaced points \(p_j\in(0,1]\) and define the
weight \(w_j\) of column \(j\) by
\begin{equation}\label{vp:quantile-weights}
 p_j=\frac{m-j+1}{m},\qquad
 w_j=\left(\frac{Q(p_j)}{p_j}\right)^L,
 \qquad 1\le j\le m.
\end{equation}
The profile has the following entries, with \(h+1\) read
cyclically:
\begin{equation}\label{vp:cycle-profile}
 s_{(h,i),(h+1,j)}=
 \begin{cases}
 w_j/m,&h=L,\\
 1/m,&h\ne L,
 \end{cases}
 \qquad\text{all other entries are zero.}
\end{equation}

Since \(Q(p)/p\) is nondecreasing and at most one, and \(p_j\)
decreases with \(j\), we have \(1\ge w_1\ge\cdots\ge w_m\ge0\).
Hence all row and column sums are at most one, and
\(\max s_{ij}=1/m\).
For the circular law, \(Q(p)=p\), so every \(w_j\) equals one and
the profile is doubly stochastic.
The denominators \(p_j^L\) in \eqref{vp:quantile-weights} are chosen
to cancel the contributions of layers \(1,\ldots,L\) to the
permanent, as the count below shows.

\par\medskip\noindent\textbf{2. Exact cancellation in the permanent.}
In a permutation with a nonzero contribution to \(\per(S_N+qI)\),
each vertex is either fixed or mapped into the next layer. Since
fixed vertices occupy their own columns, the moving vertices of each
layer are mapped bijectively onto those of the next, so every layer
contains the same number \(k\) of moving vertices.
For fixed \(k\), we count the contribution in three steps.
\begin{enumerate}
\item \emph{Choose the moving vertices.}
Each of the layers \(1,\ldots,L\) offers \(\binom{m}{k}\) choices,
which gives \(\binom{m}{k}^L\) in total.
In layer \(L+1\), a set \(I\) of moving vertices carries the factor
\(\prod_{j\in I}w_j\). Summing over the choices of \(I\) gives
\[
 e_k(w)=\sum_{\substack{I\subset[m]\\|I|=k}}\prod_{j\in I}w_j,
 \qquad e_0(w)=1.
\]
\item \emph{Match consecutive layers.}
Once the sets of moving vertices are fixed, the permutation acts on
each of the \(L+1\) links, including the closing link from layer
\(L+1\) to layer \(1\), by one of \(k!\) independently chosen
bijections, which gives \((k!)^{L+1}\).
\item \emph{Multiply the entry factors.}
Each of the \((L+1)k\) moving vertices contributes a variance
\(1/m\), apart from the factors \(w_j\) already included in
\(e_k(w)\); together these give \(m^{-(L+1)k}\). Each of the
\((L+1)(m-k)\) fixed vertices contributes a factor \(q\), which gives
\(q^{(L+1)(m-k)}\).
\end{enumerate}
The total contribution of the permutations with \(k\) moving
vertices in each layer is therefore
\[
 q^{(L+1)(m-k)}
 \binom{m}{k}^{L}e_k(w)
 (k!)^{L+1}m^{-(L+1)k}.
\]
Since \(\binom{m}{k}k!/m^k=\prod_{j=1}^k p_j\), we can combine
\(L\) of the factors \(k!/m^k\) with \(\binom{m}{k}^L\), which leaves
a single factor \(k!/m^k\). Summing over \(k=0,\ldots,m\) gives
\begin{equation}\label{vp:quantile-count}
 \per(S_N+qI)=\sum_{k=0}^m
 q^{(L+1)(m-k)}
 \left(\prod_{j=1}^k p_j\right)^L
 \frac{k!}{m^k}\,e_k(w).
\end{equation}
The weights are nonincreasing, so the largest term of \(e_k(w)\)
is \(\prod_{j=1}^k w_j\), and \(e_k(w)\) has at most \(2^m\) terms.
Hence
\[
 \prod_{j=1}^k w_j\le e_k(w)\le2^m\prod_{j=1}^k w_j.
\]
By the definition of \(w_j\), the cancellation is exact:
\begin{equation}\label{vp:quantile-cancellation}
 \left(\prod_{j=1}^k p_j\right)^L
 \prod_{j=1}^k w_j
 =\prod_{j=1}^k Q(p_j)^L.
\end{equation}

\par\medskip\noindent\textbf{3. The limiting pressure.}
The remaining factor satisfies
\(e^{-m}\le k!/m^k\le1\) for \(0\le k\le m\).
Consequently, the \(k\)-th term in \eqref{vp:quantile-count} lies
between \(e^{-m}\) and \(2^m\) times
\[
 q^{(L+1)(m-k)}\prod_{j=1}^k Q(p_j)^L.
\]
The logarithm of a sum of \(m+1\) nonnegative terms differs from the
logarithm of its largest term by at most \(\log(m+1)\).
Since \(Q(p_j)\) is nonincreasing in \(j\), the indices \(j\) with
\(Q(p_j)^L>q^{L+1}\) form an initial segment, and the largest of these
expressions is the one whose product runs exactly over this
segment.
Setting \(q=e^t\) and dividing by \(N=(L+1)m\), we therefore obtain
\[
 H_N(t)=\frac1m\sum_{j=1}^m
 \max\left\{t,\frac{L}{L+1}\log Q(p_j)\right\}
 +O\left(\frac1{L+1}\right).
\]
Replacing \(L/(L+1)\) by one changes the right side by at most
\(|t|/(L+1)\). With \(\log0=-\infty\) inside the maximum, the resulting sums are
Riemann sums of a bounded monotone function of \(p\), and they
converge to \eqref{vp:target-quantile-pressure}.
Since \(H'_+(t)=F(e^t)\) by \eqref{vp:radial-pressure-slope}, Theorem~\ref{vp:pressure-main} identifies
the limit of the ESD as \(\nu\).

For an arbitrary dimension \(N\), we take
\(m=\lfloor\sqrt N\rfloor\) and \(L+1=\lfloor N/m\rfloor\), and we append to the construction \(N-(L+1)m\) indices whose rows and columns vanish. There are fewer than \(m=o(N)\) such indices, so they do not
change the limiting measure. The maximal variance is
\(1/m=O(N^{-1/2})\), and under \eqref{pd:power} the same
pressure criterion gives almost sure convergence.
\end{proof}

The concavity used in the proof of necessity comes from the
following variational formula, which holds for every finite profile.
For \(a\ge0\) and \(q>0\), let
\begin{equation}\label{vp:scalar-maximum}
 G(a,q)=\begin{cases}
   \log q,&a\le q,\\
   \log a-1+q/a,&a>q,
 \end{cases}
\end{equation}
and for a profile \(S\) and \(\theta\in\R^N\) write
\(d_i(\theta)=\sum_j s_{ij}e^{\theta_j-\theta_i}\).

\begin{lemma}[A variational formula for the radial potential]\label{vp:radial-potential}
Let \(S\) be a nonnegative \(N\times N\) profile, let \(A_N=0\), and let
\(z\in\C\) with \(q=|z|^2>0\). Then, for every \(\eta>0\), the
regularized Brown potential satisfies
\begin{equation}\label{vp:gauge-dual}
 2NU_{N,\eta}(z)=\inf_{\theta\in\R^N}
 \sum_i\left[G(d_i(\theta)+\eta e^{-\theta_i},q)+\eta e^{\theta_i}\right].
\end{equation}
The infimum is attained at \(\theta_i=\log v_i\), where \(v\)
is the Green-function vector from \eqref{vp:dyson}.
\end{lemma}
\begin{proof}
We verify \eqref{vp:gauge-dual} by exhibiting a minimizer. The
diagonal similarity
\(S\mapsto\diag(e^{-\theta_i})S\diag(e^{\theta_i})\) turns the row
sums of \(S\) into \(d_i(\theta)\).

\noindent\textbf{Convexity of the objective.}
The function \(t\mapsto G(e^t,q)\) has the nonnegative,
nondecreasing derivative \((1-qe^{-t})_+\), so it is convex and
nondecreasing. The function
\[
 \log(d_i(\theta)+\eta e^{-\theta_i})
 =\log\left(\sum_j s_{ij}e^{\theta_j}+\eta\right)-\theta_i
\]
is convex in \(\theta\), being a log-sum-exp minus a linear function.
Hence each term \(G(d_i(\theta)+\eta e^{-\theta_i},q)\), a convex
nondecreasing function of a convex function, is convex. So is
\(\eta e^{\theta_i}\), and hence the sum in \eqref{vp:gauge-dual}.

\noindent\textbf{A minimizer from the Dyson solution.}
Fix \(\eta>0\), and let \(u,v\) be the Green-function vectors and \(g,\ell\) the auxiliary quantities from \eqref{vp:dyson}. At the point
\(\theta_i=\log v_i\),
\[
 d_i(\theta)+\eta e^{-\theta_i}
 =\frac{g_i}{v_i}=q+g_i\ell_i,\qquad
 \partial_aG(q+g_i\ell_i,q)
 =\frac{g_i\ell_i}{(q+g_i\ell_i)^2}=u_iv_i.
\]
By these identities, the partial derivative of the sum in
\eqref{vp:gauge-dual} with respect to \(\theta_k\) at this point is
\[
 v_k(S^{\mathsf T}u)_k-u_kg_k+\eta v_k
 =v_k\ell_k-u_kg_k=0.
\]
By convexity of the objective, this critical point is a global
minimizer. Its objective value is
\(\sum_i[\log(q+g_i\ell_i)-u_ig_i+\eta v_i]\),
which equals \(2NU_{N,\eta}(z)\) by \eqref{vp:stationary} and the
identity \(\sum_i u_i=\sum_i v_i\).
\end{proof}

\begin{proof}[Necessity in Theorem~\ref{vp:radial-inverse}]
By rescaling, we may assume that \(d_*=1\). Let \(\nu_N\) be the
Brown measure associated with \(S_N\); since \(A_N=0\), it is
rotation invariant.

\par\smallskip\noindent\textbf{The support condition.}
Take \(\theta=0\) and \(q=1\) in \eqref{vp:gauge-dual}.
Since \(d_i(0)\le1\) and \(G(a,1)=0\) for \(a\le1\), letting \(\eta\downarrow0\) gives \(U_N(1)\le0\). By rotation invariance and Jensen's circular-mean formula \eqref{vp:circular-average}, \(2U_N(1)\) is a nonnegative integral, so
\[
 0=2U_N(1)=\int\log\max\{1,|w|^2\}\dd\nu_N(w).
\]
Since the integrand is nonnegative and strictly positive outside the
closed unit disk, \(\nu_N\) is supported on that disk, and so is every
weak limit of \((\nu_N)\).

\par\smallskip\noindent\textbf{The radial condition.}
Since \(G(a,\cdot)\) is concave for each \(a\ge0\),
\eqref{vp:gauge-dual} exhibits \(2NU_{N,\eta}\) as an infimum of
concave functions of \(q=|z|^2\); letting \(\eta\downarrow0\) shows
that \(q\mapsto2U_N(\sqrt q)\) is concave.
Suppose that \(\nu_N\Rightarrow\nu\) along a subsequence. Jensen's circular-mean formula \eqref{vp:circular-average} identifies the pointwise limit of \(2U_N(\sqrt q)\) along
this subsequence as
\(\int\log\max\{q,|w|^2\}\dd\nu(w)\), for \(q>0\).
This limit is concave in \(q\), and its right derivative is
\(F(q)/q\). Hence \(F(q)/q\) is
nonincreasing, which proves \eqref{vp:radial-ratio}.
Finally, Theorem~\ref{vp:brown} transfers these conclusions from the
Brown measures to the limits of the ESD.
\end{proof}

\subsubsection{A radial law with an annular gap}\label{sec:annular-gap}

\begin{numberedexample}[A disk separated from an annulus]\label{vp:annular-example}
Fix \(0<r_-<r_+<1\) and define \(\lambda=(r_-/r_+)^2\).
Let \(Q(p)=\lambda p\) for \(0<p\le r_+^2\) and \(Q(p)=p\) for
\(r_+^2<p\le1\). Then \(Q(p)/p\) is nondecreasing, and the
corresponding distribution function satisfies \eqref{vp:radial-ratio}.
In the construction from the proof of
Theorem~\ref{vp:radial-inverse}, the weights
\eqref{vp:quantile-weights} are \(w_j=\lambda^L\) for the columns
with \(p_j\le r_+^2\) and \(w_j=1\) for the remaining columns. The limiting law is
\[
 \nu(\mathrm dz)=
 \frac{1}{\pi\lambda}{\bf1}_{\{|z|<r_-\}}\dd A(z)
 +\frac1\pi{\bf1}_{\{r_+<|z|<1\}}\dd A(z).
\]
Its support therefore has an open annular gap, and the origin
carries no atom.
\end{numberedexample}

\begin{numberedremark}[Comparison with earlier profile assumptions]\label{vp:profile-assumptions}
The profiles in Example~\ref{vp:annular-example} violate
Assumption~A1 of \cite[Theorem~2.3]{vp-chnr1}, since
\(N\max_{i,j}s_{ij}\to\infty\).
\end{numberedremark}

\subsection{Permanental zeros}\label{sec:permanental-zeros}

The characteristic polynomial and the permanental polynomial have the
same second moment, and hence the same second-moment pressure, so
the comparison in Section~\ref{sec:variance-profiles} approximates
both zero measures by the same deterministic measure.

\begin{proposition}[Zeros of the permanent and the determinant]\label{vp:phase-brown}\label{pp:esd-comparison}
Assume the hypotheses of Proposition~\ref{vp:brown-log-tail}.
Let \(\mu_N\) be the ESD of \(A_N+X_N\),
and let \(\mu_N^{\mathrm{per}}\) be the normalized zero measure of
\[
 P_N(z)=\per(zI-A_N-X_N).
\]
Let \(\nu_N\) be the Brown measure formed from the full profile
\(S_N\) and the diagonal field \(A_N\). Then
\begin{equation}\label{pp:common-brown}
 \E d_{\BL}(\mu_N^{\mathrm{per}},\nu_N)
 +\E d_{\BL}(\mu_N^{\mathrm{per}},\mu_N)\longrightarrow0.
\end{equation}
Under \eqref{pd:power}, both distances also tend to zero almost
surely, under any coupling across dimensions. If, moreover,
\(\nu_N\Rightarrow\nu\), then both random measures converge weakly to
\(\nu\) in probability, and almost surely under \eqref{pd:power}.
The same conclusions hold for the normalized zero measure \(\mu_N^F\) of
every polynomial \(F_N\) defined in Section~\ref{sec:general-comparison}.
\end{proposition}
\begin{proof}
In the definition of \(F_N\) in Section~\ref{sec:general-comparison},
the choice \(\omega_N(\pi)=\operatorname{sgn}(\pi)\) gives
\(\det(zI-A_N-X_N)\), and the choice \(\omega_N(\pi)=1\) gives
\(P_N(z)\). Both choices satisfy \(|\omega_N(\pi)|=1\) and
\(\omega_N(\mathrm{id})=1\), so Proposition~\ref{vp:brown-log-tail}
compares both zero measures \(\mu_N\) and
\(\mu_N^{\mathrm{per}}\) with the same Brown measure \(\nu_N\).
The triangle inequality then gives
\eqref{pp:common-brown} and its almost sure version. The same
comparison with \(\nu_N\) applies to every fixed choice of
unit-modulus coefficients.
\end{proof}

\begin{proof}[Proof of Theorem~\ref{vp:permanent}]
For a doubly stochastic profile with \(A_N=0\),
Theorem~\ref{vp:circle} gives \(\mu_N\Rightarrow\mu_{\rm circ}\)
in probability, and almost surely under \eqref{pd:power}.
Proposition~\ref{vp:phase-brown} gives
\(d_{\BL}(\mu_N^{\mathrm{per}},\mu_N)\to0\) in the same two modes,
and the asserted convergence of \(\mu_N^{\mathrm{per}}\) follows from
the triangle inequality.
\end{proof}

\section{Proof of the elliptic law}\label{sec:reciprocal}

Let \(S=S_N\) be symmetric and doubly stochastic, with
\(\epsilon_N=\max_{i,j}s_{ij}\to0\), and let \(X_N\) have the
reciprocal-pair law of Theorem~\ref{ep:main}. As in Section~\ref{sec:planar-dini}, write
\[
 \begin{gathered}
 D_N(z)=\det(zI-X_N),\qquad M_N(z)=\E|D_N(z)|^2,\\
 u_N(z)=N^{-1}\log|D_N(z)|,\qquad p_N(z)=(2N)^{-1}\log M_N(z).
 \end{gathered}
\]
Recall the parameters \(\tau=\E(UV)\), \(\kappa=\E(U\overline V)\) and \(\rho=|\tau|\) of the pair model, and let \(\theta\) be the angle in
the definition of \(\mu_\tau\), so that \(\tau=\rho e^{2i\theta}\).
Multiplying \(X_N\) and \(z\) by \(e^{-i\theta}\) replaces \(\tau\) by
\(\rho\), leaves \(\kappa\) and the moduli of the determinants
unchanged, and maps \({\cal E}_\tau\) onto the horizontal ellipse
\(\mathcal E_\rho\). It therefore suffices to identify the limit in
these rotated coordinates as the uniform probability measure \(\mu_\rho\) on \(\mathcal E_\rho\), whose semiaxes are \(1+\rho\) and \(1-\rho\). We write \(U_\rho\) for the logarithmic potential of \(\mu_\rho\) and \(K_\rho=[-2\sqrt\rho,2\sqrt\rho]\) for the focal segment of
\(\mathcal E_\rho\). The second-moment estimates of
Sections~\ref{sec:ep-gaussian} and~\ref{sec:ep-scale} are stated in
the rotated coordinates, whereas the concentration and replacement
estimates of Section~\ref{sec:ep-concentration} are stated in the
original variables.

The proof follows the three steps of the circular-law proof. We
identify the second-moment pressure, close the Jensen gap for bounded
approximants and return to the original entry law.
Two features of the model require new estimates. In the permutation
expansion of \(D_N(z)\), each transposition \((i\,j)\) contributes a
factor with mean \(-\tau s_{ij}\), and these mean transpositions
produce a matching polynomial in the second moment. The dependent
reciprocal entries must also be resampled together, which makes the
conditional determinant quadratic in the resampled variables rather
than affine as in Section~\ref{sec:planar-dini}.
The cycle-weight computations and the small-ball bound are in
Appendix~\ref{app:stars}.

\subsection{Reduction to a Gaussian second moment}\label{sec:ep-gaussian}

As a reference, we use the Gaussian pair
\begin{equation}\label{ep:canonical}
 U^{\mathrm G}=g,\qquad
 V^{\mathrm G}=\tau\bar g+\sqrt{1-|\tau|^2}\,h,
\end{equation}
where \(g\) and \(h\) are independent standard circular Gaussians,
with \(\E|g|^2=1\) and \(\E g^2=0\). Let \(G_N^\tau\) be the matrix
built from this pair on the profile \(S\), with independent
centered Gaussian diagonal variables of unit variance, and write
\(Z_S^\tau(z)=\E|\det(zI-G_N^\tau)|^2\).
For \(J\subset[N]\), define the matching polynomial
\[
 m_J^\tau(z)=\sum_{\mathfrak m\ {\rm matching\ on}\ J}
 (-\tau)^{|\mathfrak m|}z^{|J|-2|\mathfrak m|}\prod_{\{i,j\}\in\mathfrak m}s_{ij},
 \qquad m_\varnothing^\tau(z)=1.
\]
Here \(\mathfrak m\) runs over sets of disjoint pairs in \(J\), and
\(|\mathfrak m|\) is the number of pairs. In the second moment,
\(m_J^\tau(z)\) collects the terms in which each vertex of \(J\) is
either a fixed point with factor \(z\) or lies in a mean
transposition.

From here to the end of Section~\ref{sec:ep-scale} we work in the
rotated coordinates, keeping the same symbols for the rotated matrix
and shift, so that \(\tau=\rho\).

\begin{proposition}[Gaussian reduction]\label{ep:cycle-comparison}
In the rotated coordinates,
\begin{equation}\label{ep:Z}
 Z_S^\rho(z)=\sum_{I\subset[N]}\per(S[I,I])\,|m_{I^c}^\rho(z)|^2,
\end{equation}
and, for each \(z\notin K_\rho\),
\begin{equation}\label{ep:cyclecompare}
 N^{-1}|\log M_N(z)-\log Z_S^\rho(z)|\longrightarrow0.
\end{equation}
\end{proposition}
\begin{proof}
Once the mean transpositions are collected, the remaining Gaussian
monomials are orthogonal, and their squared norms add up to the
principal permanents in \eqref{ep:Z}. For a general pair law, only
the variance of a transposition and the weights of oppositely
oriented cycles change. Summed over
the cycles \(c\) of a fixed length \(k\ge2\), the variance products
\(\prod_{e\in c}s_e\) are at most \(N\epsilon_N/k\), and for \(k\ge3\)
each cycle weight differs from its Gaussian value \(1\) by
\(O(|\kappa|^k)\). Away
from \(K_\rho\), the deletion bound for matching polynomials makes
this comparison uniform along an interpolation between
\(Z_S^\rho(z)\) and \(M_N(z)\).
Appendix~\ref{app:cycle-comparison} gives the exact expansion and
the resulting estimate.
\end{proof}

\subsection{The second-moment scale}\label{sec:ep-scale}

For each fixed auxiliary law, that is, a bounded approximant from
Lemma~\ref{ep:bounded-pairs} or a Gaussian law, we need two
properties. The second-moment pressure \(p_N\) should approach
\(U_\rho\), and the fourth moment should have the same exponential
scale as \(M_N^2\).
As in Section~\ref{sec:planar-dini}, we measure the second property by
\[
 r_N(z)=\frac1N\log\frac{\E|D_N(z)|^4}{M_N(z)^2}
\]
whenever \(M_N(z)>0\) and the fourth moment is finite.
By Cauchy--Schwarz, \(r_N(z)\ge0\), and the second property means
that \(r_N(z)\to0\).

\begin{proposition}[The second-moment scale for auxiliary laws]\label{ep:aux-pressure}
Let the entry law be either a fixed bounded approximating law from
Lemma~\ref{ep:bounded-pairs} or a fixed Gaussian law satisfying the
assumptions of Theorem~\ref{ep:main}. Then, for each
\(z\notin K_\rho\),
\begin{equation}\label{ep:relative}
 p_N(z)\to U_\rho(z),\qquad r_N(z)\to0 .
\end{equation}
\end{proposition}

The proof combines the fourth-moment upper bound of
Lemma~\ref{ep:upper} with the lower bound for the Gaussian second
moment in
Lemma~\ref{ep:lower}, which Proposition~\ref{ep:cycle-comparison}
transfers to \(M_N\). Since \(M_N(z)^2\le\E|D_N(z)|^4\), the two bounds
force both limits in \eqref{ep:relative}.

\subsubsection{The fourth-moment upper bound}

\begin{lemma}[Fourth-moment upper bound]\label{ep:regularized}\label{ep:upper}
Suppose that the pair and diagonal laws are fixed, bounded and
satisfy a logarithmic Sobolev inequality, or that they are Gaussian
as above. Then
\begin{equation}\label{ep:fourthupper}
 \limsup_N N^{-1}\log\E|D_N(z)|^4\le4U_\rho(z)
 \qquad(z\in\C).
\end{equation}
\end{lemma}
\begin{proof}
\noindent\textbf{1. The regularized mean.}
We use the regularized logarithm \(R_{N,\eta}(X_N;z)\) of \eqref{pd:regularized-log}, with \(\eta\) in place of \(\delta\), and compare \(X_N\) with a dense,
centered, real Gaussian elliptic matrix \(\widehat G_n\). The
entries of \(\widehat G_n\) have variance \(1/n\), and the two entries
of each reciprocal pair have covariance \(\rho/n\). Distinct unordered
pairs and the diagonal entries are independent. Let
\(U_{\rho,\eta}(z)\) be the limit of
\(\E R_{n,\eta}(\widehat G_n;z)\) as \(n\to\infty\).
Grouping the walk count of Proposition~\ref{vp:ds} by unordered
pairs gives
\begin{equation}\label{ep:fixed}
 \E R_{N,\eta}(X_N;z)\longrightarrow U_{\rho,\eta}(z)
 \qquad(\eta>0).
\end{equation}
Indeed, only doubled trees survive in this count; in such a tree, a pair of
opposite letters contributes \(s_{ij}\) and a pair of equal letters
\(\rho s_{ij}\), while \(\kappa\) and the marginal pseudovariances do not enter. The resulting moments are exactly the tree
moments of the
dense Gaussian model, and the logarithmic tails at fixed \(\eta\) are controlled as in that proposition.

To estimate the regularization error, we use the classical
\(2\times2\) Hermitized equations for this Gaussian limit; see
\cite[(3.4)--(3.6) and Section~4.1]{ep-alt-kruger}. Write \(iv\),
with \(v>0\), for the diagonal entries of their solution and \(b\in\C\) for its upper-right entry. The equations then read
\[
 \frac{v}{v^2+|b|^2}=\eta+v,\qquad
 \frac{b}{v^2+|b|^2}+\rho\bar b=-z.
\]
The first equation gives \(v^2\le v^2+|b|^2=v/(\eta+v)<1\), so \(v<1\).
Since \(U_{\rho,\eta}(z)\to U_\rho(z)\) as \(\eta\downarrow0\), integrating \(\partial_\eta U_{\rho,\eta}=v\) in \(\eta\) gives
\begin{equation}\label{ep:free}
 0\le U_{\rho,\eta}(z)-U_\rho(z)\le\eta.
\end{equation}

\par\smallskip\noindent\textbf{2. From the mean to the fourth moment.}
The logarithmic Sobolev calculation in
Proposition~\ref{vp:ds} gives
\[
 \log\E e^{4N R_{N,\eta}(X_N;z)}
 \le4N\E R_{N,\eta}(X_N;z)+C_{\rm law}N\epsilon_N/\eta^2 .
\]
Since \(u_N(z)\le R_{N,\eta}(X_N;z)\), the left side is an upper
bound for \(\log\E|D_N(z)|^4\). Dividing by \(N\), letting \(N\to\infty\) at
fixed \(\eta\) by \eqref{ep:fixed} and \(\epsilon_N\to0\), and then
letting \(\eta\downarrow0\) by \eqref{ep:free}, we obtain
\eqref{ep:fourthupper}.
\end{proof}

\subsubsection{The lower bound for the Gaussian second moment}
Since the terms of \eqref{ep:Z} are nonnegative, the sum of any subfamily of them bounds \(Z_S^\rho(z)\) from below.

\begin{lemma}[Lower bound for the Gaussian second moment]\label{ep:lower}
For each \(z\notin K_\rho\),
\[
 \liminf_N N^{-1}\log Z_S^\rho(z)\ge2U_\rho(z).
\]
\end{lemma}
\begin{proof}
For \(\rho=0\), \eqref{ep:Z} reduces to \(Z_S^0(z)=\per(S+|z|^2I)\) and Proposition~\ref{vp:ds} applies, so we assume
\(0<\rho<1\). Write \(U_{\rm sc,a}\) for the logarithmic
potential of the centered semicircle law of variance \(a\).
The subset calculation in Appendix~\ref{app:matching-pressure}
proves, for every fixed \(0<\lambda<1\),
\begin{equation}\label{ep:lambda-lower}
 \liminf_N\frac1N\log Z_S^\rho(z)
 \ge -1+\lambda-\lambda\log\lambda
             +2\lambda U_{\rm sc,\rho\lambda}(z).
\end{equation}
This calculation counts a family of sets \(J=I^c\) of size close to
\(\lambda N\), uses the semicircle limit for the roots of their
matching polynomials, and bounds the complementary permanents from
below after repairing the diagonal. The repair changes
\(\log Z_S^\rho(z)\) by \(o(N)\).

Now let \(z=x+iy\) lie inside the ellipse, and choose \(\lambda\)
and \(\vartheta\) by
\[
 \lambda=\frac{x^2}{(1+\rho)^2}+\frac{y^2}{(1-\rho)^2}\in(0,1),
 \qquad z=\sqrt\lambda(e^{i\vartheta}+\rho e^{-i\vartheta}).
\]
The semicircle calculation \eqref{ep:semipot} in the appendix gives
\(2U_{\rm sc,\rho\lambda}(z)=\rho\cos(2\vartheta)+\log\lambda\),
so the right side of \eqref{ep:lambda-lower} becomes
\[
 -1+\lambda+\rho\lambda\cos(2\vartheta)
 =-1+\frac{x^2}{1+\rho}+\frac{y^2}{1-\rho}
 =2U_\rho(z).
\]
If \(z\notin\mathcal E_\rho\), the
full-set analogue of \eqref{ep:matching-uniform} gives
\(N^{-1}\log|m_{[N]}^\rho(z)|\to U_{\rm sc,\rho}(z)=U_\rho(z)\).
In this case the term of \eqref{ep:Z} with \(I=\varnothing\) alone
gives \(Z_S^\rho(z)\ge|m_{[N]}^\rho(z)|^2\), which suffices.
\end{proof}

\subsubsection{Combining the bounds}

\begin{proof}[Proof of Proposition~\ref{ep:aux-pressure}]
Since \(2p_N(z)=N^{-1}\log M_N(z)\), Proposition~\ref{ep:cycle-comparison}
and Lemma~\ref{ep:lower} give
\(\liminf_N 2p_N(z)\ge2U_\rho(z)\).
On the other hand, \(M_N(z)^2\le\E|D_N(z)|^4\), so
Lemma~\ref{ep:upper} implies
\[
 \limsup_N 2p_N(z)
 \le\tfrac12\limsup_N\frac1N\log\E|D_N(z)|^4
 \le2U_\rho(z).
\]
Thus \(p_N(z)\to U_\rho(z)\). Finally,
\[
 0\le r_N(z)=\frac1N\log\E|D_N(z)|^4-4p_N(z),
\]
and by Lemma~\ref{ep:upper} and the convergence of \(p_N(z)\), the
upper limit of the right side is at most
\(4U_\rho(z)-4U_\rho(z)=0\), so \(r_N(z)\to0\).
\end{proof}

It remains to close the Jensen gap for the entry law of
Theorem~\ref{ep:main} itself, by logarithmic concentration and
replacement.

\subsection{Concentration and approximation}\label{sec:ep-concentration}

We now return to the original, unrotated variables, so
\(\tau=\E(UV)\) may be complex. To preserve the dependence within
reciprocal pairs during replacement, we resample whole pairs,
grouped into stars. The \(i\)-th \emph{star} is the group
\[
 \bigl(\zeta_i,(U_{ij},V_{ij})_{j>i}\bigr),\qquad 1\le i\le N.
\]
Every pair and every diagonal variable belongs to exactly one star.
Since the pairs and the diagonal variables are mutually independent,
so are the \(N\) stars. The entries coming from star \(i\) lie in row \(i\) or column \(i\),
and a permutation term uses at most one entry from each, with a
diagonal entry counting for both. Hence, conditionally on the other stars, the
determinant is a polynomial of degree at most two in the variables of
star \(i\).

We therefore replace the linear-form bound of Lemma~\ref{pd:row} by
a small-ball bound for quadratic polynomials, uniform in the
coefficients and in the number of pairs. This bound is proved in
Appendix~\ref{app:star-proof} and yields the next two propositions.

\begin{proposition}[Logarithmic concentration]\label{ep:log-concentration}
Let the pair and diagonal laws be as in Theorem~\ref{ep:main}, and
suppose that the determinant at \(z\) is a nonzero polynomial on the product of the real supports. Then \(Y=\log|D_N(z)|\) is integrable and,
for \(s>0\),
\begin{equation}\label{ep:fluct}
 \Pr(|Y-\E Y|>s)\le2e^{-c\min(s^2/N,s)},\qquad
 \E|Y-\E Y|^4\le C N^2 .
\end{equation}
If, in addition, \(\E|D_N(z)|^4<\infty\), then \(r_N(z)\) controls the distance between \(u_N(z)\) and
\(p_N(z)\):
\begin{equation}\label{ep:mean}
 \E|u_N(z)-p_N(z)|\le C(\sqrt{r_N(z)}+r_N(z)+N^{-1/2}).
\end{equation}
\end{proposition}

\begin{proposition}[Transport of logarithms under star replacement]\label{ep:log-transport}
Let \((U,V)\) and \((\widetilde U,\widetilde V)\) be two pair laws,
and let \(\zeta\) and \(\widetilde\zeta\) be two diagonal laws.
Assume that both entry models satisfy the centering, unit-variance,
density and mixed-moment hypotheses of Theorem~\ref{ep:main}, with
pair densities on the same real linear subspace \(L\subset\C^2\) and
diagonal densities on the same subspace \(L_0\subset\C\), and that
the two pair laws have the same parameter
\(\tau=\E(UV)=\E(\widetilde U\widetilde V)\).

Given a coupling \((U,V,\widetilde U,\widetilde V)\) of the pair
laws and a coupling \((\zeta,\widetilde\zeta)\) of the diagonal laws,
we take a copy of the coupled four-tuple for each \(i<j\) and a copy
of the coupled diagonal variables for each \(i\), all independent, and
define, on the same
profile \(S_N\) for both models,
\[
 \begin{aligned}
 (X_{ij},X_{ji})&=\sqrt{s_{ij}}(U_{ij},V_{ij}),&
 (\widetilde X_{ij},\widetilde X_{ji})
   &=\sqrt{s_{ij}}(\widetilde U_{ij},\widetilde V_{ij})\quad(i<j),\\
 X_{ii}&=\sqrt{s_{ii}}\,\zeta_i,&
 \widetilde X_{ii}&=\sqrt{s_{ii}}\,\widetilde\zeta_i.
 \end{aligned}
\]
The normalized log determinants are
\[
 u_N(z)=\frac1N\log|\det(zI-X_N)|,\qquad
 \widetilde u_N(z)=\frac1N\log|\det(zI-\widetilde X_N)| .
\]
The coupling errors are
\[
 d^2=\E\!\left[|U-\widetilde U|^2+|V-\widetilde V|^2+
                 |UV-\widetilde U\widetilde V|^2\right],
 \qquad d_0^2=\E|\zeta-\widetilde\zeta|^2 .
\]
Then, at each \(z\) for which \(\det(zI-X_N)\) is a nonzero polynomial on the product of the real supports,
\begin{equation}\label{ep:transport}
 \E|u_N(z)-\widetilde u_N(z)|\le C(d+d_0)^{1/49}
 \qquad (d+d_0\le1).
\end{equation}
The constant \(C\) depends only on the supports and on common upper
bounds for the densities and the mixed moments.
\end{proposition}

\begin{lemma}[Approximation of a reciprocal pair]\label{ep:bounded-pairs}
Let the pair and diagonal laws be as in Theorem~\ref{ep:main}. Then
there are bounded approximating laws \((U_m,V_m)\) and \(\zeta_m\),
\(m\ge1\), on the same supports \(L\) and \(L_0\), with the same means
and real covariances, coupled with the given pair and diagonal
variables. In this coupling,
\(U_m\), \(V_m\), \(U_mV_m\) and \(\zeta_m\) converge in \(L^2\) to
\(U\), \(V\), \(UV\) and \(\zeta\), respectively. The densities of the
approximating laws and the mixed moments \(\E|U_mV_m|^2\) are bounded
uniformly in \(m\). Each approximating law satisfies a logarithmic
Sobolev inequality, with a constant that may depend on \(m\).
\end{lemma}
The proofs of the concentration and transport propositions are given
in Appendix~\ref{app:star-resampling}, and
Appendix~\ref{app:coupling} constructs the approximants. Since the
approximants have the same covariances as the original variables,
\(\tau\) and \(\kappa\) do not change. The coupling errors \(d\) and
\(d_0\) of Proposition~\ref{ep:log-transport} tend to zero by the
\(L^2\) convergence in the lemma. The logarithmic Sobolev constant may depend on \(m\) because
the inequality is used only at fixed \(m\).

\begin{proof}[Proof of Theorem~\ref{ep:main}]
We work in the rotated coordinates, where \(\tau=\rho\). The rotated
pair and diagonal laws satisfy the same hypotheses, so the results of
Section~\ref{sec:ep-concentration} apply to them.
Let \(z\notin K_\rho\); in particular \(z\ne0\). The determinant
equals \(z^N\ne0\) when all the pair and diagonal variables vanish, so it is a nonzero polynomial on the product of the real supports. Let \(u_N^{(m)}\) denote
the normalized log determinants built from the bounded approximants
of Lemma~\ref{ep:bounded-pairs}. For fixed \(m\),
Proposition~\ref{ep:aux-pressure} and \eqref{ep:mean} give
\(\E|u_N^{(m)}(z)-U_\rho(z)|\to0\) as \(N\to\infty\). Since the
transport estimate \eqref{ep:transport} is uniform in \(N\), letting
\(m\to\infty\) gives \(\E|u_N(z)-U_\rho(z)|\to0\) for every \(z\) off
the focal segment.

The focal segment has zero area, so this convergence holds for
almost every \(z\). To pass to integrals over \(z\), we use the elementary bound
\[
 \int_{D_R}|\log|z-w||^2\dd A(z)\le C_R(1+|w|^2).
\]
Writing \(\mu_N\) for the ESD of \(X_N\), we obtain from Jensen's
inequality and the Schur decomposition that
\begin{align}
 \E\int_{D_R}|u_N(z)|^2\dd A(z)
 &\le C_R\left(1+\E\int|w|^2\dd\mu_N(w)\right)\notag\\
 &\le C_R\left(1+N^{-1}\E\|X_N\|_{\rm F}^2\right)
 =2C_R.\label{pd:potential-l2}
\end{align}
This uniform \(L^2\) bound gives uniform integrability on the
product of the probability space and \(D_R\), so the pointwise
convergence above upgrades to
\(\E\|u_N-U_\rho\|_{L^1(D_R)}\to0\).
The centered fourth-moment bound in \eqref{ep:fluct} also gives almost sure convergence in local \(L^1\), by Hölder's inequality and the Borel--Cantelli lemma, as in the proof of Lemma~\ref{pd:row-concentration}.
Lemma~\ref{pd:root-return} then yields convergence of the empirical
measures in mean bounded-Lipschitz distance and almost surely, and
rotating back by \(e^{i\theta}\) gives the claimed limit \(\mu_\tau\).
\end{proof}

\appendix
\section{Technical estimates for reciprocal pairs}\label{app:stars}

This appendix proves the technical estimates used in
Section~\ref{sec:reciprocal}. Throughout, we use the reciprocal-pair
model of that section with a symmetric doubly stochastic profile
\(S=S_N\) satisfying \(\epsilon_N=\max_{i,j}s_{ij}\to0\).
Sections~\ref{app:cycle-comparison} and~\ref{app:matching-pressure}
identify the second-moment scale. Sections~\ref{app:star-proof}
and~\ref{app:star-resampling} control the logarithm of the determinant
when a star is resampled, and Section~\ref{app:coupling} constructs
the bounded approximants.

\subsection{Cycle weights and the Gaussian reduction}\label{app:cycle-comparison}
Recall \(\tau=\E(UV)\), \(\kappa=\E(U\overline V)\) and \(\rho=|\tau|\).
We compute how the pair law enters \(M_N(z)=\E|D_N(z)|^2\), and then
prove the comparison \eqref{ep:cyclecompare} of Proposition~\ref{ep:cycle-comparison}. The matching polynomials \(m_J^\tau\) are defined in
Section~\ref{sec:reciprocal}.

In the permutation expansion of \(D_N(z)\), a transposition
\((i\,j)\) uses both entries of the reciprocal pair \(\{i,j\}\), whose product \(X_{ij}X_{ji}=s_{ij}U_{ij}V_{ij}\) has mean \(\tau s_{ij}\). The centered
product \(UV-\tau\) has variance
\[
 \gamma=\E|UV-\tau|^2=\E|UV|^2-|\tau|^2>0.
\]
Indeed, since \(U\) and \(V\) have unit variance, each coordinate linear form is nonzero on the real linear support
\(L\subset\C^2\), so \(UV\) is a nonconstant quadratic on \(L\). Its
level set \(\{UV=\tau\}\) is Lebesgue-null in \(L\), and since
\((U,V)\) has a density on \(L\), \(UV\) is not almost surely equal to
\(\tau\).

To state the second moment, we weight each permutation by its cycles.
We write \(s_e=s_{ij}\) for an edge \(e=(i,j)\). For a directed
simple cycle \(c=(i_1,\ldots,i_k)\), \(k\ge3\), we let \(n_+(c)\) count
its steps from a smaller to a larger index, including \(i_k\to i_1\),
and set \(n_-(c)=k-n_+(c)\). The weights are
\[
 w(c)=1+\Re(\kappa^{n_+(c)}\bar\kappa^{\,n_-(c)}),\qquad
 W(\pi)=\gamma^{C_2(\pi)}
        \prod_{\substack{c\subset\pi\\ |c|\ge3}}w(c).
\]
Here \(C_2(\pi)\) is the number of transpositions in the cycle
decomposition of \(\pi\), \(|c|\) is the length of a cycle \(c\) of
\(\pi\), and fixed points carry weight one. Since \(|\kappa|<1\), every
\(w(c)\) is positive. Write \(\mathfrak S(I)\) for the set of
permutations of \(I\).

\begin{proposition}[Second moment for reciprocal pairs]\label{ep:second-moment}
For every \(N,z\),
\begin{equation}\label{ep:second}
 M_N(z)
 =\sum_{I\subset[N]} |m_{I^c}^\tau(z)|^2
       \sum_{\pi\in\mathfrak S(I)} W(\pi)\prod_{i\in I}s_{i,\pi(i)} .
\end{equation}
\end{proposition}
\begin{proof}
In the permutation expansion of \(D_N(z)\), each cycle of length at
least two contributes minus the product of its entries, and each fixed
point \(i\) contributes \(z\) or \(-X_{ii}\). We split the
product \(s_{ij}U_{ij}V_{ij}\) of each transposition into its mean
\(\tau s_{ij}\) and its centered part \(s_{ij}(U_{ij}V_{ij}-\tau)\).
For a set \(I^c\) of vertices, the terms that cover \(I^c\) by
factors \(z\) and mean transpositions sum to \(m_{I^c}^\tau(z)\).
The vertices of the complementary set \(I\) are covered by random
loops \(-X_{ii}\), centered transpositions and cycles of length at
least three. Two such random monomials with different unordered edge
supports are orthogonal, since an edge that appears in only one of the
two supports supplies an independent centered factor.

A centered transposition contributes \(\gamma s_{ij}^2\) to
\(M_N(z)\). A cycle of length at least three has two orientations, each
with squared norm \(\prod_{e\in c}s_e\). For a triangle, after the
variance factors are removed, the two orientations are
\(U_{12}U_{23}V_{13}\) and \(V_{12}V_{23}U_{13}\), with inner product
\(\kappa^2\bar\kappa\). In general, each increasing step of \(c\)
contributes a factor \(\kappa\) and each decreasing step a factor
\(\bar\kappa\). Hence the two orientations together contribute
\(2w(c)\prod_{e\in c}s_e\). The sum in \eqref{ep:second} lists
both orientations, so each receives the weight \(w(c)\).
The only fourth-order moment that enters is \(\E|UV|^2\), through
\(\gamma\), so \(\E|UV|^2<\infty\) suffices.

For the Gaussian reference pair \eqref{ep:canonical}, we have
\(\kappa=0\) and \(\gamma=1\), since
\[
 U^{\mathrm G}V^{\mathrm G}-\tau
 =\tau(|g|^2-1)+\sqrt{1-|\tau|^2}\,gh .
\]
The two terms on the right are uncorrelated and have variances
\(|\tau|^2\) and \(1-|\tau|^2\). Thus every cycle weight equals one,
and \eqref{ep:second} becomes the expansion of \(Z_S^\tau(z)\).
\end{proof}

\begin{proof}[Cycle-weight comparison in Proposition~\ref{ep:cycle-comparison}]
\textbf{Deletion controls the interpolating partition function.}
We rotate so that \(\tau=\rho\), and fix \(z\) with
\(\operatorname{dist}(z,K_\rho)\ge h>0\).
Fix \(J\subset[N]\), and give each pair \(\{u,v\}\subset J\) with \(s_{uv}>0\) the weight
\(\sqrt{\rho s_{uv}}\). By the path-tree representation
\cite{ep-hl,ep-godsil}, every root of \(m_J^\rho\) is an eigenvalue of
the weighted adjacency matrix of a path tree of some connected
component of the resulting graph on \(J\). Since
each row of \(S\) sums to one, the parent-to-child map of a path tree
has norm at most \(\sqrt\rho\), so the adjacency matrix has norm at
most \(2\sqrt\rho\), and all roots of \(m_J^\rho\) lie in \(K_\rho\).
Moreover, \(m_{J\setminus\{i\}}^\rho(z)/m_J^\rho(z)\) is the root
entry of the resolvent at \(z\) of the path tree rooted at \(i\); as
this tree has spectrum in \(K_\rho\), the ratio has modulus at most
\(1/h\). Hence
\begin{equation}\label{ep:matching-deletion}
 |m_J^\rho(z)|\ge h\,|m_{J\setminus\{i\}}^\rho(z)|.
\end{equation}
For \(0\le t\le1\) and \(\Lambda\subset[N]\), define the positive sum
\[
 \mathcal Z_t(\Lambda)=\sum_{I\subset\Lambda}|m_{\Lambda\setminus I}^\rho(z)|^2
       \sum_{\pi\in\mathfrak S(I)}
       W(\pi)^t\prod_{i\in I}s_{i,\pi(i)} .
\]
It interpolates between \(\mathcal Z_0([N])=Z_S^\rho(z)\) and
\(\mathcal Z_1([N])=M_N(z)\).
For \(i\in\Lambda\), we keep only the sets \(I\) that do not contain
\(i\) and apply \eqref{ep:matching-deletion} with
\(J=\Lambda\setminus I\); this gives
\(\mathcal Z_t(\Lambda)\ge h^2\mathcal Z_t(\Lambda\setminus\{i\})\)
for every \(t\in[0,1]\).

\noindent\textbf{Short and long cycles give separate errors.}
Normalizing the terms of \(\mathcal Z_t([N])\) gives a probability
law on the pairs \((I,\pi)\), and the \(t\)-derivative of
\(\log\mathcal Z_t([N])\) is the average of
\(C_2(\pi)\log\gamma+\sum_{|c|\ge3}\log w(c)\) under this law.
Fix \(k\ge2\) and a \(k\)-cycle \(c\). Deleting the vertices of \(c\)
and applying the deletion bound for the interpolating partition
functions \(k\) times shows that \(\pi\) contains \(c\) with
probability at most \(C_\gamma h^{-2k}\prod_{e\in c}s_e\). Summed over all \(k\)-cycles,
these variance products are at most
\(\operatorname{Tr}S^k/k\le N\epsilon_N/k\), since the trace also
counts closed walks that repeat vertices, and every entry of \(S^k\) is
at most \(\epsilon_N\). Hence, for fixed \(k_0\ge3\), the cycles of length
less than \(k_0\) contribute \(O_{h,k_0,\gamma,\kappa}(N\epsilon_N)\)
to the average. The cycles of \(\pi\) of length at least \(k_0\) are
disjoint, so \(\pi\) has at most \(N/k_0\) of them, and each satisfies
\(|\log w(c)|\le C_\kappa|\kappa|^{|c|}\le C_\kappa|\kappa|^{k_0}\).
Integrating the logarithmic derivative over \(t\in[0,1]\) gives
\[
 \frac1N|\log M_N(z)-\log Z_S^\rho(z)|
 \le C_{h,k_0,\gamma,\kappa}\epsilon_N
       +\frac{C_\kappa|\kappa|^{k_0}}{k_0}.
\]
Letting first \(N\to\infty\) and then \(k_0\to\infty\) proves
\eqref{ep:cyclecompare}.
\end{proof}

\subsection{A lower bound from matching polynomials}\label{app:matching-pressure}
We prove the lower bound \eqref{ep:lambda-lower} used in the proof
of Lemma~\ref{ep:lower}, and then evaluate its semicircle term.

\begin{proof}[Details for Lemma~\ref{ep:lower}]
Fix \(0<\rho<1\), \(z\notin K_\rho\) and \(0<\lambda<1\).
Write \(h=\operatorname{dist}(z,K_\rho)>0\), \(t=1-\lambda\)
and \(h_{\rm bin}(\lambda)=-\lambda\log\lambda-t\log t\).
In \eqref{ep:Z}, the term indexed by \(I=J^c\) is
\(|m_J^\rho(z)|^2\per(S[I,I])\).

\noindent\textbf{1. Counting typical subsets.}
Choose a random subset \(J\subset[N]\) by including each vertex
independently with probability \(\lambda\), and let
\(\delta=\epsilon_N^{1/4}\) and
\(\beta_N=2e^{-1/\sqrt{\epsilon_N}}\).
Call a row \(i\) bad if
\(\bigl|\sum_{j\in J}s_{ij}-\lambda\bigr|>\delta\), and good
otherwise. Since \(\sum_j s_{ij}^2\le\epsilon_N\), Hoeffding's
inequality bounds the probability that a given row is bad by
\(2e^{-2/\sqrt{\epsilon_N}}\). By Markov's inequality, with
probability tending to one there are at most \(N\beta_N\) bad rows. Moreover,
\(\Pr(\bigl||J|/N-\lambda\bigr|>\delta)\le2e^{-2N\delta^2}\to0\),
because \(\epsilon_N\ge1/N\).
Let \(\mathcal J_N\) be the family of subsets \(J\) that have at most
\(N\beta_N\) bad rows and satisfy \(\bigl||J|/N-\lambda\bigr|\le\delta\).
Each \(J\in\mathcal J_N\) has probability
\(\lambda^{|J|}t^{N-|J|}=\exp[-Nh_{\rm bin}(\lambda)+O_\lambda(N\delta)]\).
Since the total probability of \(\mathcal J_N\) tends to one,
\(\#\mathcal J_N\ge\exp[Nh_{\rm bin}(\lambda)-o(N)]\).

\par\smallskip\noindent\textbf{2. Uniform control of the matching roots.}
Let \(B=(b_{ij})\) be a symmetric \(n\times n\) matrix with zero
diagonal, nonnegative entries and row sums at most one. Let \(m_B\)
be its matching polynomial and \(\nu_B\) its normalized root measure.
We write \(\varepsilon=\max_{i,j}b_{ij}\) and, for \(a\in[0,1]\),
\[
 \Delta=\frac1n\sum_i\left|\sum_jb_{ij}-a\right|.
\]
With \(\mathrm{Cat}_k=(k+1)^{-1}\binom{2k}{k}\), we claim that
\begin{equation}\label{ep:matching-moments}
 \left|\int x^{2k}\dd\nu_B(x)-\mathrm{Cat}_k a^k\right|
 \le\mathrm{Cat}_k\left(k\Delta+\binom{k+1}{2}\varepsilon\right),
 \qquad k\ge1,
\end{equation}
and that all odd moments of \(\nu_B\) vanish. We use the
deletion and derivative identities
\[
 m_B(w)=w m_{B\setminus i}(w)
       -\sum_{j\ne i}b_{ij}m_{B\setminus\{i,j\}}(w),
 \qquad m'_B(w)=\sum_i m_{B\setminus i}(w),
\]
which follow by sorting the matchings according to the partner of
\(i\) and by counting unmatched vertices. Applied recursively, the first identity shows that
\(m_{B\setminus i}/m_B\) is the root resolvent \(R_i(w)\) of the tree
of simple paths starting at \(i\), in which the edge extending a path
from \(u\) to \(v\) has weight \(\sqrt{b_{uv}}\)
\cite{ep-godsil,ep-ach}. Together with the derivative identity, this
gives
\[
 \int\frac{\mathrm d\nu_B(x)}{w-x}
 =\frac{m'_B(w)}{n m_B(w)}=\frac1n\sum_iR_i(w).
\]
Expanding in powers of \(w^{-1}\) shows that each moment
of \(\nu_B\) is an average over \(i\) of weighted counts of closed
walks from the root of the path tree of \(i\). A closed walk of
length \(2k\) has one of the \(\mathrm{Cat}_k\) Dyck shapes, and each
shape encodes a rooted plane tree \(T\) with \(k\) edges. Allowing
arbitrary labels in \([n]\) on the vertices of \(T\) gives the
weight
\[
 t_n(T;B)=\frac1n\sum_{\phi:V(T)\to[n]}
                  \prod_{\{u,v\}\in E(T)}b_{\phi(u)\phi(v)}.
\]
The path-tree count omits exactly the labelings in which some vertex
of \(T\) repeats the label of one of its ancestors. For such a
labeling, identifying the two vertices creates a cycle; we bound the
weight of one of its edges by \(\varepsilon\) and delete that edge. Summing over the labels of the remaining tree,
leaf by leaf, then shows that the labelings with a given coincidence have
normalized weight at most \(\varepsilon\). At most \(\binom{k+1}{2}\) coincidences between
a vertex and an ancestor are possible, so the path-tree count differs
from \(t_n(T;B)\) by at most \(\binom{k+1}{2}\varepsilon\).

It remains to compare \(t_n(T;B)\) with \(a^k\). Summing over the
label of a leaf of \(T\) produces the row sum of \(B\) at the label of
its parent. Replacing this row sum by \(a\) costs at most \(\Delta\),
because after rerooting \(T\) at the parent, the symmetry of \(B\) and the
row-sum bound show that the sum of all remaining factors is at most
one. Removing the leaves in this way, once for each edge, gives
\(|t_n(T;B)-a^k|\le k\Delta\). Adding the two errors and summing over
the \(\mathrm{Cat}_k\) shapes proves \eqref{ep:matching-moments}. The
odd moments vanish because path trees are bipartite.

We apply the moment bound with \(a=\rho\lambda\) to
\(B=(\rho s_{ij}{\bf1}_{i\ne j})_{i,j\in J}\), for which
\(m_B=m_J^\rho\) and \(\varepsilon\le\rho\epsilon_N\). The row error
\(\bigl|\sum_jb_{ij}-a\bigr|\) is at most \(\rho\) for a bad row and at
most \(\rho(\delta+\epsilon_N)\) for a good row, where the term
\(\epsilon_N\) accounts for the deleted diagonal entry. The average row
error \(\Delta_J\) of \(B\) therefore satisfies
\begin{equation}\label{ep:matching-row-error}
 \Delta_J\le\rho(\delta+\epsilon_N)+\frac{\rho N\beta_N}{|J|}
 \longrightarrow0
 \qquad\text{uniformly for }J\in\mathcal J_N.
\end{equation}
By the path-tree bound preceding \eqref{ep:matching-deletion}, the
roots of \(m_J^\rho\) lie in \(K_\rho\). As \(m_J^\rho\) is monic,
\(|J|^{-1}\log|m_J^\rho(z)|\) is the integral of \(\log|z-x|\)
against the normalized root measure of \(m_J^\rho\), and for \(z\notin K_\rho\) this
function of \(x\) is a uniform limit of polynomials on \(K_\rho\). By
\eqref{ep:matching-moments} and \eqref{ep:matching-row-error}, the
moments of the root measures converge, uniformly for
\(J\in\mathcal J_N\), to those of the semicircle law of variance
\(\rho\lambda\). Therefore
\begin{equation}\label{ep:matching-uniform}
 \sup_{J\in\mathcal J_N}
 \left|\frac1{|J|}\log|m_J^\rho(z)|
                  -U_{\rm sc,\rho\lambda}(z)\right|\longrightarrow0.
\end{equation}
For the full set \(J=[N]\), the same argument applies with \(a=\rho\)
and \(\Delta\le\rho\epsilon_N\).

\par\smallskip\noindent\textbf{3. Diagonal repair and the complementary permanents.}
Let \(\widehat S=S+\epsilon_N I_N\). Along the path \(S+rI_N\),
\(0\le r\le\epsilon_N\), the matching polynomials do not change,
since they involve only off-diagonal entries. The derivative of
\(Z_{S+rI_N}^\rho(z)\) in one diagonal entry is the partition function
with that vertex deleted. By the deletion bound for partition
functions following \eqref{ep:matching-deletion}, its ratio to the
full partition function is at most \(h^{-2}\). Summing these
logarithmic derivatives over the diagonal entries and integrating
over \(r\) gives
\begin{equation}\label{ep:repair-cost-main}
 0\le\log Z_{\widehat S}^\rho(z)-\log Z_S^\rho(z)
 \le N\epsilon_N/h^2=o(N).
\end{equation}

Fix \(J\in\mathcal J_N\), and let \(I=J^c\), \(\widehat T=\widehat S[I,I]\)
and \(c=t+\delta+\epsilon_N\). We split \(I\) into the sets
\(I_{\rm good}\) and \(I_{\rm bad}\) of good and bad rows, and define
\begin{equation}\label{ep:trial-matrix}
 \begin{aligned}
 P_{ij}&=s_{ij}/c &&(i,j\in I_{\rm good},\ i\ne j),\\
 P_{ii}&=1-c^{-1}\sum_{j\in I_{\rm good}\setminus\{i\}}s_{ij}
                          &&(i\in I_{\rm good}),\\
 P_{ii}&=1 &&(i\in I_{\rm bad}),
 \end{aligned}
\end{equation}
and set all other entries equal to zero. Since the rows of
\(S\) sum to one, every good row \(i\) satisfies
\(\sum_{j\in I}s_{ij}\in[t-\delta,t+\delta]\). Hence \(P\) is
nonnegative, symmetric and row-stochastic, so it is doubly stochastic;
each of its positive entries is also a positive entry of \(\widehat T\). For a good row \(i\),
\(cP_{ii}\le2\delta+2\epsilon_N+\sum_{j\in I_{\rm bad}}s_{ij}\).
Because the column sums of \(S\) are one, the last sum, added over all
good rows, is at most \(|I_{\rm bad}|\le N\beta_N\). The diagonal mass
\(\mathfrak d=\sum_{i\in I}P_{ii}\) therefore obeys
\begin{equation}\label{ep:diagonal-mass}
 \mathfrak d\le N\beta_N+\frac{2N(\delta+\epsilon_N)+N\beta_N}{c},
 \qquad \mathfrak d|\log\epsilon_N|=o(N).
\end{equation}
The second bound holds because \(c\to t>0\) and
\((\delta+\epsilon_N+\beta_N)|\log\epsilon_N|\to0\).
Where \(P\) is positive, \(\widehat T_{ij}/P_{ij}=c\) off the diagonal and
\(\widehat T_{ii}/P_{ii}\ge\epsilon_N\) on the diagonal. Gurvits' bound
\cite[Theorem~2.2]{vp-bethe}, combined with the inequality
\((1-p)\log(1-p)\ge-p\), therefore gives
\begin{equation}\label{ep:principal-permanent}
 \begin{aligned}
 \log\per\widehat T
 &\ge(|I|-\mathfrak d)(\log c-1)+\mathfrak d(\log\epsilon_N-1)\\
 &\ge |I|(\log t-1)-\mathfrak d|\log\epsilon_N|
   =|I|(\log t-1)-o(N),
 \end{aligned}
\end{equation}
where the last line uses \(t\le c<1\) for large \(N\).
All these errors are uniform over \(J\in\mathcal J_N\). We bound
\(Z_{\widehat S}^\rho(z)\) from below by the terms of \eqref{ep:Z} for
\(\widehat S\) that are indexed by \(I=J^c\) with \(J\in\mathcal J_N\).
After taking logarithms and dividing by \(N\), the count of
\(\mathcal J_N\) from Step~1 contributes \(h_{\rm bin}(\lambda)\), the
permanent bound \eqref{ep:principal-permanent} the term
\(t(\log t-1)\), and \eqref{ep:matching-uniform} the term
\(2\lambda U_{\rm sc,\rho\lambda}(z)\). Removing the diagonal repair by
\eqref{ep:repair-cost-main} then gives
\[
 \liminf_N N^{-1}\log Z_S^\rho(z)
 \ge h_{\rm bin}(\lambda)+t(\log t-1)+2\lambda U_{\rm sc,\rho\lambda}(z),
\]
which is \eqref{ep:lambda-lower}, since the terms \(\pm t\log t\)
cancel and \(t=1-\lambda\).

\par\smallskip\noindent\textbf{4. The semicircle term.}
For \(a>0\), the logarithmic potential and the Cauchy transform of the
centered semicircle law of variance \(a\) are
\[
 U_{\rm sc,a}(z)=\int_{-2\sqrt a}^{2\sqrt a}
      \log|z-w|\frac{\sqrt{4a-w^2}}{2\pi a}\dd w,
 \qquad
 g_a(z)=\int_{-2\sqrt a}^{2\sqrt a}
      \frac{\sqrt{4a-w^2}}{2\pi a(z-w)}\dd w.
\]

Off the support \([-2\sqrt a,2\sqrt a]\), the Cauchy transform
satisfies \(z=g_a^{-1}+ag_a\), and \(g_a(z)\sim z^{-1}\) as
\(z\to\infty\). For fixed \(a\), we have
\(\mathrm dU_{\rm sc,a}=\Re(g_a\dd z)\), and the relation between \(z\)
and \(g_a\) gives \(g_a\dd z=\mathrm d(-\log g_a+\tfrac a2g_a^2)\).
Taking real parts and normalizing at infinity, we obtain
\begin{equation}\label{ep:semipot}
 2U_{\rm sc,a}(z)=\Re\bigl(zg_a(z)\bigr)-1-2\log|g_a(z)|.
\end{equation}
For \(z\) inside the ellipse, the proof of Lemma~\ref{ep:lower}
chooses \(\lambda\) and \(\vartheta\) with
\(z=\sqrt\lambda(e^{i\vartheta}+\rho e^{-i\vartheta})\). Then
\(g_{\rho\lambda}(z)=e^{-i\vartheta}/\sqrt\lambda\). Indeed, this value solves
the quadratic equation \(z=g^{-1}+\rho\lambda g\), and since
\(\sqrt{\rho\lambda}\,|e^{-i\vartheta}/\sqrt\lambda|=\sqrt\rho<1\), it
is the root selected by the Cauchy transform.
Substituting it into \eqref{ep:semipot} gives
\(2U_{\rm sc,\rho\lambda}(z)=\rho\cos(2\vartheta)+\log\lambda\),
which is the identity used in that proof.
\end{proof}

\subsection{A quadratic polynomial on one star}\label{app:star-proof}
We prove the small-ball estimate that Section~\ref{sec:reciprocal}
uses both for concentration and for star replacement. Its
constants must not depend on the coefficients of the conditional
determinant, which are fixed by the other stars, nor on the number of
pairs in a star.

We fix one star, condition on all the other stars, and index the pairs of this star by \(j\ge1\), using \(0\) for its diagonal variable. Let \(x_j\) be the real orthonormal coordinates of
\((U_j,V_j)\) on \(L\subset\C^2\), and \(x_0\) the coordinates of the
diagonal variable on \(L_0\subset\C\). In these coordinates, the
product of the two coordinate linear forms is
\(U_jV_j=x_j^{\mathsf T}Q_Lx_j\), where \(Q_L\) is a nonzero complex
symmetric matrix that depends only on \(L\). We call the vectors \(x_j\) the pair blocks and \(x_0\) the diagonal
block.
The conditional determinant then has the form
\begin{equation}\label{ep:star}
 P=\mu_P+\ell_0^{\mathsf T}x_0+
     \sum_j\bigl\{\ell_j^{\mathsf T}x_j+\beta_j(U_jV_j-\tau)\bigr\}
             +\sum_{j<k}x_j^{\mathsf T}K_{jk}x_k.
\end{equation}
In \eqref{ep:star}, \(x_0\) appears only linearly, because a
permutation that selects the diagonal entry of the star selects no
other entry of it. The products \(U_jV_j\) come from permutations
that select both entries of a pair. All coefficients, including the
variance-profile factors, are fixed by the other stars and \(z\).
Vector norms are Euclidean, and
\(\|K\|_{\rm F}^2=\sum_{a,b}|K_{ab}|^2\). The scales of \(P\) are
\begin{equation}\label{ep:star-scale}
 \sigma_P^2=\|\ell_0\|^2+
     \sum_j\bigl(\|\ell_j\|^2+\|\beta_jQ_L\|_{\rm F}^2\bigr)
             +\sum_{j<k}\|K_{jk}\|_{\rm F}^2,\qquad
 s_P^2=|\mu_P|^2+\sigma_P^2.
\end{equation}
Since the blocks are centered and the product terms are
centered by \(\tau\), \(\mu_P=\E P\) is the conditional mean of \(P\).
Given the other stars and \(z\), coupled laws on the same supports
with the same \(\tau\) have the same \(\mu_P,\sigma_P,s_P\).

\begin{lemma}[Small balls on a star]\label{ep:stars}
Let the pair blocks and the diagonal block in \eqref{ep:star} be
independent and centered, with unit marginal second moments. Assume that they have uniformly bounded densities on the fixed real supports \(L\) and
\(L_0\), that \(\E(U_jV_j)=\tau\) for every \(j\), and that
\(\E|U_jV_j|^2\) is bounded uniformly in \(j\). The pair laws may
differ from block to block. Then, for \(0<t<1\),
\[
 \begin{aligned}
 \sup_{w\in\C}\Pr(|P-w|\le t\sigma_P)&\le Ct^{1/24}\quad(\sigma_P>0),\\
 \Pr(|P|\le ts_P)&\le Ct^{1/48}\quad(s_P>0).
 \end{aligned}
\]
Moreover, \(c_1s_P^2\le\E|P|^2\le C_1s_P^2\).
The constants depend only on the supports and on the density and
mixed-moment bounds; in particular, they are uniform in the
coefficients and in the number of pairs.
\end{lemma}

\begin{proof}
\noindent\textbf{1. Single blocks and linear forms.}
Consider first a single block of dimension \(d\le4\) and a real
polynomial \(P_{\rm blk}(x)=\ell^{\mathsf T}x+x^{\mathsf T}Q_{\rm blk}x\)
on it, with \(Q_{\rm blk}\) symmetric and
\(\|\ell\|^2+\|Q_{\rm blk}\|_{\rm F}^2=1\). Write \(\E_\square\) for the
expectation under the uniform probability on the cube \([-R,R]^d\),
where \(R\ge1\). Under this law the linear monomials and the centered
quadratic monomials are orthogonal, which gives
\(\E_\square[(P_{\rm blk}-w)^2]\ge c\) for every real \(w\). After
rescaling to the unit cube, the equivalence of the \(L^1\) and \(L^2\)
norms on polynomials of degree at most two gives
\(\E_\square|P_{\rm blk}-w|\ge c\) as well. The Carbery--Wright
inequality \cite[Theorem~8]{ep-cw} therefore bounds the uniform
probability of the event \(|P_{\rm blk}-w|\le s\) on the cube by
\(Cs^{1/2}\). Under the actual block law, the density bound controls
the part of this event inside the cube, and the second moment controls
the probability that the block leaves the cube, which gives
\[
 \Pr(|P_{\rm blk}-w|\le s)\le C(R^4s^{1/2}+R^{-2})\le Cs^{1/6},
 \qquad R=s^{-1/12}.
\]
The diagonal block is the special case \(Q_{\rm blk}=0\).

For a linear form with unit vector \(\theta\), the slab
\(\{|\theta^{\mathsf T}x-w|\le s\}\) meets the cube in a set of volume
at most \(CsR^{d-1}\). The same argument, now with \(R=s^{-1/5}\), shows that
the probability of this slab is at most \(Cs^{2/5}\). By the rank-one
projection bound \cite[Corollary~1.4]{ep-rv}, the same estimate holds
for every unit linear combination of independent block projections.
Integrating it gives a uniform bound on the negative moment of order
\(1/3<2/5\). Consequently, for a concatenated block vector \(Y\) and
every nonzero real matrix \(K\),
\begin{equation}\label{ep:neg}
 \E\|KY\|^{-1/3}\le C\|K\|_{\rm F}^{-1/3}.
\end{equation}
Indeed, by the singular-value decomposition,
\(\|KY\|^2/\|K\|_{\rm F}^2\) is a convex combination of squared unit
projections of \(Y\), and \eqref{ep:neg} follows from the convexity of
\(x\mapsto x^{-1/6}\). Conditioning on an independent copy, we obtain
the same linear estimates for the difference of a block vector and
this copy.

For the first assertion, we reduce to a real polynomial. Choose
whichever of \(\Re P\) and \(\Im P\) has the larger nonconstant
coefficient norm, which is at least \(\sigma_P/\sqrt2\), and divide it
by that norm to obtain \(P_{\rm re}\).
If \(P\) lies in a disk of radius \(t\sigma_P\), then \(P_{\rm re}\)
lies in an interval of radius at most \(\sqrt2\,t\). It therefore
suffices to bound the probability that \(P_{\rm re}\) lies in a given
interval of radius \(s\le\sqrt2\,t\). Let \(K_*\) be the combined Frobenius norm of the
cross-block matrices in \(P_{\rm re}\).

\par\smallskip\noindent\textbf{2. Large cross terms.}
Suppose that \(K_*\ge\sqrt s\). For a uniformly random bipartition of the
pair blocks, the cross terms between the two groups have expected
squared coefficient norm \(K_*^2/2\). Fix a bipartition attaining at
least this value, and put the diagonal block, which has no cross terms, into either group.
Let \(Y\) and \(Z\) be the coordinate vectors of the two groups and
\(K\) the matrix of the cross terms between them, so that
\(\|K\|_{\rm F}\ge K_*/\sqrt2\).
Let \(p=\Pr(|P_{\rm re}-w|\le s)\), and let \(Y',Z'\) be independent
copies of \(Y,Z\). Applying the Cauchy--Schwarz inequality twice, once
in each group of variables, we bound \(p^4\) by the probability that
the four events \(|P_{\rm re}-w|\le s\) at
\((Y,Z),(Y',Z),(Y,Z'),(Y',Z')\) all hold. The alternating sum of the
four values of \(P_{\rm re}\) is
\[
 P_{\rm re}(Y,Z)-P_{\rm re}(Y',Z)-P_{\rm re}(Y,Z')+P_{\rm re}(Y',Z')
   =(Y-Y')^{\mathsf T}K(Z-Z');
\]
all terms internal to either group, and the shift \(w\), cancel. On
the intersection of the four events, this sum has modulus at most
\(4s\). Hence
\[
 p^4\le\Pr(|(Y-Y')^{\mathsf T}K(Z-Z')|\le4s)
 \le Cs^{1/3}\E\|K(Z-Z')\|^{-1/3}
 \le C(s/\|K\|_{\rm F})^{1/3}\le Cs^{1/6}.
\]
For the middle bounds, we condition on \(Z-Z'\), apply the linear
estimate, with exponent weakened to \(1/3\), to the resulting linear
form in \(Y-Y'\), and then use \eqref{ep:neg}.
Taking the fourth root gives \(p\le Cs^{1/24}\).

\par\smallskip\noindent\textbf{3. Small cross terms.}
Suppose now that \(K_*<\sqrt s\). Cross terms belonging to distinct
pairs of blocks are orthogonal, because a block that appears in only
one of them supplies a centered factor. Hence the cross part of
\(P_{\rm re}\) has squared \(L^2\) norm at most \(Cs\), and by
Chebyshev's inequality it exceeds \(s^{1/4}\) in modulus with
probability at most \(Cs^{1/2}\). If \(|P_{\rm re}-w|\le s\) and the
exceptional event does not occur, the remaining part of
\(P_{\rm re}\), a sum of independent one-block polynomials, lies
within \(s+s^{1/4}\le2s^{1/4}\) of \(w\) when \(s\le1\). Since the
nonconstant coefficients of \(P_{\rm re}\) have norm one, the combined
squared coefficient norm of these one-block polynomials is
\(1-K_*^2\), which is bounded below for small \(s\). After normalizing each nonconstant one-block polynomial by its coefficient norm, we apply the one-block bound
\(Cr^{1/6}\) at scale \(r=Cs^{1/4}\); the rank-one projection bound
then gives \(Cs^{1/24}\) for their weighted sum. Adding the probability of the exceptional event gives
\(p\le Cs^{1/24}\) in this case too, and since \(s\le\sqrt2t\), this
proves the first assertion for complex disks.

\par\smallskip\noindent\textbf{4. The constant term.}
For the moment comparison, recall that \(\E P=\mu_P\). Moreover, \(\E|P-\mu_P|^2\le C\sigma_P^2\), since terms from different
blocks, and cross terms from distinct pairs of blocks, are orthogonal,
while within a single block the Cauchy--Schwarz inequality and the bound
\(\E|U_jV_j|^2\le C\) control the linear and quadratic terms. For the
quadratic term we use
\(|\beta_j|=\|\beta_jQ_L\|_{\rm F}/\|Q_L\|_{\rm F}\), where
\(\|Q_L\|_{\rm F}>0\) depends only on \(L\). Conversely, the first
small-ball bound, applied with \(w=\mu_P\) at a small fixed value of
\(t\), gives \(\E|P-\mu_P|^2\ge c_1\sigma_P^2\). Since
\(\E|P|^2=|\mu_P|^2+\E|P-\mu_P|^2\), these two bounds prove
\(c_1s_P^2\le\E|P|^2\le C_1s_P^2\).

For the second small-ball bound, suppose first that
\(\sigma_P<\sqrt t\,s_P\). Then \(|\mu_P|^2=s_P^2-\sigma_P^2>(1-t)s_P^2\),
so for small \(t\) the event \(|P|\le ts_P\) forces
\(|P-\mu_P|\ge s_P/2\). By Chebyshev's inequality and the upper
variance bound, this event has probability \(O(t)\); when
\(P\) is a nonzero constant, it is empty. Now suppose that
\(\sigma_P\ge\sqrt t\,s_P\). Then \(ts_P/\sigma_P\le\sqrt t\), and the
first assertion with
\(w=0\) gives the bound \(Ct^{1/48}\). This proves the bound for small
\(t\), and enlarging the constants covers the remaining
\(t\in(0,1)\).
\end{proof}

\subsection{Concentration and star replacement}\label{app:star-resampling}
We adapt the resampling argument of Section~\ref{sec:planar-dini} to
stars, using the conditional determinant \(P\) of \eqref{ep:star} and
the scales \(\sigma_P\) and \(s_P\) of \eqref{ep:star-scale}.

\begin{proof}[Proof of Proposition~\ref{ep:log-concentration}]
\noindent\textbf{Integrability.}
A nonzero polynomial on the product of the real supports vanishes only on a
Lebesgue null set. Call a monomial allowed if it contains from each pair either no
variable, a single coordinate or the product \(UV\), and from each
diagonal variable at most a single coordinate. Let \(F\) be a nonzero polynomial built from allowed
monomials, such as \(D_N(z)\). Single out one pair and expand \(F\) in
the functions \(1\), the real coordinates of that pair, and
\(UV-\tau\); the coefficients are polynomials in the remaining
blocks. Since \(Q_L\ne0\), these functions are linearly independent,
so some coefficient \(f\) is a nonzero polynomial. Conditionally on
the remaining blocks, \(F\) becomes a one-block polynomial \(P\) with
\(s_P\ge c_L|f|\). Integrating the small-ball bound of
Lemma~\ref{ep:stars} over the chosen pair therefore gives
\(\E_{\rm pair}\log^-|P|\le C+\log^-|f|\).
Since coefficient extraction preserves the allowed monomials in the
remaining blocks, we can repeat the argument with \(f\) in place of
\(F\); for a diagonal block, the
expansion uses \(1\) and the coordinates on \(L_0\). By induction on
the number of blocks, \(\log^-|D_N|\) is integrable, and the same
induction applies to hybrids of the coupled laws. The finite second
moment controls \(\log^+|D_N|\).

\noindent\textbf{Resampling one star.}
Conditionally on all stars but one, the determinant is the polynomial
\(P\) of \eqref{ep:star} in the variables of that star. If \(s_P>0\), Lemma~\ref{ep:stars} bounds both tails
of \(\log(|P|/s_P)\):
\[
 \Pr\!\left(\left|\log\frac{|P|}{s_P}\right|>s\right)
 \le Ce^{-cs}.
\]
The lower tail follows from the small-ball bound, and the upper tail
from \(\E|P|^2\le C_1s_P^2\) and Markov's inequality. Since \(s_P\)
is determined by the other stars, it is unchanged by resampling, so
\(\log s_P\) cancels in the difference of the logarithms before and
after resampling, and the tail bound gives the conditional
exponential-moment hypothesis of Lemma~\ref{pd:resampling}. That
lemma then proves \eqref{ep:fluct} and \eqref{ep:mean}.
\end{proof}

\begin{proof}[Proof of Proposition~\ref{ep:log-transport}]
\noindent\textbf{Comparing a coupled star on the same scale.}
We replace the stars one at a time. Fixing the other stars, let
\(P\) and \(\widetilde P\) be the conditional determinants before and
after one replacement, and let \(d\) and \(d_0\) be the coupling
errors of the proposition. The common value \(\tau\) makes the
conditional mean \(\mu_P\) and the scales \(\sigma_P,s_P\) the
same for \(P\) and \(\widetilde P\). Because the coordinates are
orthonormal, \(\E\|x_j-\widetilde x_j\|^2\le d^2\); moreover,
\(UV-\widetilde U\widetilde V\) is centered. The replacement changes a
cross term by
\[
 (x_j-\widetilde x_j)^{\mathsf T}K_{jk}x_k+
 \widetilde x_j^{\mathsf T}K_{jk}(x_k-\widetilde x_k).
\]
This change has mean zero conditionally on either coupled block. Consequently, the changes in distinct cross terms are
orthogonal to each other and to the changes in the one-block terms,
which are centered and independent across blocks. Together with the identity
\(|\beta_j|=\|\beta_jQ_L\|_{\rm F}/\|Q_L\|_{\rm F}\)
for the product term, these facts give
\(\E|P-\widetilde P|^2\le C(d+d_0)^2\sigma_P^2\).

\noindent\textbf{Clipping.}
For \(0<\delta\le1\), we clip both logarithms below \(\log(\delta s_P)\),
replacing \(\log|P|\) by \(\log\max\{|P|,\delta s_P\}\) and
\(\log|\widetilde P|\) likewise. The clipped logarithm is
\(1/(\delta s_P)\)-Lipschitz, so the bound on
\(\E|P-\widetilde P|^2\) and the inequality \(\sigma_P\le s_P\) bound
the expected absolute difference of the clipped logarithms by
\(C(d+d_0)/\delta\). As in the proof of
Proposition~\ref{pd:transport}, Lemma~\ref{ep:stars} bounds the
expected error caused by clipping each logarithm by
\(C\delta^{1/48}\). Summing over the \(N\) star replacements and
dividing by \(N\), we obtain
\(\E|u_N(z)-\widetilde u_N(z)|\le C(d+d_0)/\delta+C\delta^{1/48}\).
If \(d+d_0>0\), the choice \(\delta=(d+d_0)^{48/49}\) balances the two
errors and gives \eqref{ep:transport}. If \(d+d_0=0\), the two
determinants agree.
\end{proof}

\begin{samepage}
\subsection{A bounded law with the prescribed covariance}\label{app:coupling}
\begin{proof}[Proof of Lemma~\ref{ep:bounded-pairs}]
\noindent\textbf{Truncation and centering.}
Let \(W\) be the vector of real orthonormal coordinates of \((U,V)\)
on \(L\), and let \(\Gamma=\E WW^{\mathsf T}\). Since \((U,V)\) has a
density on \(L\), \(\Gamma\) is positive definite; write \(\lambda_*=\lambda_{\min}(\Gamma)>0\). All the modifications below produce pairs
supported on \(L\), so preserving \(\Gamma\) also preserves the real
covariance of the pair in \(\C^2\).

On the event \(E_m=\{|U|\le m,\ |V|\le m\}\) we keep the pair, and on
its complement we redraw it from its law conditioned on \(E_m\).
For \(Y=U,V,UV\), the mean-square change is at most
\[
 2\E[|Y|^2{\bf1}_{E_m^c}]
 +2\Pr(E_m^c)\E[|Y|^2\mid E_m]\longrightarrow0.
\]
We then center the resulting pair. Its coordinate means tend to zero
and its second moments stay bounded, so by Cauchy--Schwarz the
centering changes the product by \(o(1)\) in \(L^2\).

\noindent\textbf{A positive lower bound for the density.}
Choose a cube in the coordinates \(W\), centered at the origin, with
radius \(R_m\to\infty\) and containing the support of the centered
pair. Mixing in the uniform law on this cube with positive
probability \(o(R_m^{-4})\) changes the product by
\(o(1)\) in \(L^2\), keeps the density bounded above uniformly in
\(m\), and makes the density positive throughout the cube.

\noindent\textbf{Exact covariance.}
Let \(\Gamma_m\) be the covariance of \(W\) after these steps, and
write \(h_m=\|\Gamma_m-\Gamma\|_{\rm op}\to0\).
If \(h_m=0\), the covariance is already \(\Gamma\). Otherwise, for large \(m\),
set \(\omega_m=4h_m/\lambda_*<1\). We mix in, with probability
\(\omega_m\), a replacement law; for the mixture to have covariance
\(\Gamma\), this law must have covariance
\[
 \Gamma_m^\sharp=\Gamma_m+(\Gamma-\Gamma_m)/\omega_m.
\]
For large \(m\) we have \(\lambda_{\min}(\Gamma_m)\ge3\lambda_*/4\),
and \(\|\Gamma_m^\sharp-\Gamma_m\|_{\rm op}=\lambda_*/4\), so the
matrices \(\Gamma_m^\sharp\) are uniformly positive definite and
bounded. As the replacement law we take the uniform law on an ellipsoid centered at the origin with covariance \(\Gamma_m^\sharp\). Its support and density are
bounded uniformly in \(m\), and for large \(m\) the ellipsoid lies in
the cube. The mixture then has covariance exactly
\[
 (1-\omega_m)\Gamma_m+\omega_m\Gamma_m^\sharp=\Gamma.
\]
The expected squared changes in \(U\), \(V\) and \(UV\) are
\(O(\omega_m)\), because their second moments are bounded both before
this step and under the ellipsoid law.

\noindent\textbf{The remaining properties.}
The same construction on \(L_0\) approximates the diagonal variable.
For the resulting laws, the coupling errors \(d\) and \(d_0\) tend to
zero, and the density upper bounds and mixed-moment bounds are uniform
in \(m\). The means and real covariances are exactly those of the
original laws; in particular, \(\tau\) and \(\kappa\) are unchanged.
For fixed \(m\), the density of the approximating law is bounded
above and away from zero on its cube, so the Holley--Stroock
perturbation argument \cite{pd-hs-perturb} transfers to it the
logarithmic Sobolev inequality of the uniform law on the cube, with a
constant depending on \(m\). This suffices, since \(m\) is held fixed
as \(N\to\infty\).
\end{proof}
\end{samepage}

\section{An atomic entry law without the circular law}
\label{app:atomic-profile}

We construct the profiles of Remark~\ref{rem:atomic-law}. Fix
\(0<\alpha<1/3\), and set
\[
 d=\lfloor N^\alpha\rfloor+1,\qquad k=\lceil4\log_2 d\rceil .
\]
Since \(kd^3=O(N^{3\alpha}\log N)=o(N)\), for all sufficiently
large \(N\) we can write
\[
 N=Ld+kd^3+r,\qquad d\le r<2d,\qquad L\ge2 .
\]
For the finitely many remaining dimensions, take \(s_{ij}=1/N\)
for all \(i,j\).

To build the support graph, we start from \(L\) layers
\(V_1,\ldots,V_L\) of \(d\) vertices each, with all directed edges
from \(V_\ell\) to \(V_{\ell+1}\), where \(V_{L+1}=V_1\). We then
replace each of the \(d^2\) edges from \(u\in V_L\) to \(v\in V_1\)
by a return path through \(k\) new blocks \(B_1,\ldots,B_k\) of size
\(d\). With one vertex \(x_j\) chosen in each \(B_j\), this path is
\[
 u\longrightarrow x_1\longrightarrow\cdots\longrightarrow x_k
 \longrightarrow v ,
\]
and we call its \(k+1\) edges its connecting edges. Each \(B_j\) carries
all directed edges, including loops, except the loop at \(x_j\). All
blocks used by the return paths are pairwise disjoint and disjoint
from the layers, so the \(d^2\) paths use \(kd^3\) vertices outside
the layers. Related chains of
bipartite blocks were used by Goel, Kapralov and Khanna to study
perfect matchings under edge sampling
\cite[Section~4]{vp-gkk-sampling}.

Replacing an edge by a path preserves the degrees of its endpoints.
At each \(x_j\), the path supplies the incoming and outgoing edges
lost when its loop is removed. Every vertex in the layers and blocks
therefore has indegree and outdegree \(d\). Set \(s_{ij}=1/d\) on
each directed edge \(i\to j\). On the remaining \(r\) vertices take
a separate \(r\times r\) block with all entries \(1/r\), and set all
unspecified entries to zero. The resulting profile is doubly stochastic.
Since \(r\ge d>N^\alpha\), its largest entry is \(1/d<N^{-\alpha}\).
Thus the largest standard deviation is below \(N^{-\alpha/2}\), and
\(N^\gamma(\max_{i,j}s_{ij})^{1/2}\to0\) for every
\(0<\gamma<\alpha/2\), as in Tikhomirov's question.

For the three-point law of Remark~\ref{rem:atomic-law}, the entries of
\(X_N\) on the \(k+1\) connecting edges of a fixed return path are
nonzero with probability \(1/2\) each, independently of one another.
Let \(E_N\) be the event that every return path has a connecting edge
whose entry is zero.
By the union bound,
\[
 \Pr(E_N^c)\le d^2\,2^{-(k+1)}\le\frac{1}{2d^2}.
\]
Every edge leaving \(V_\ell\) with \(\ell<L\) ends in \(V_{\ell+1}\),
and every edge leaving \(V_L\) is the first connecting edge of a return
path. A directed cycle through a layer vertex therefore enters some
return path, and since the only edge leaving each block is the next
connecting edge, it uses all \(k+1\) connecting edges of that path. On
\(E_N\) one of their entries vanishes, so no directed cycle in the
support of \(X_N\) meets a layer.

Let \(W\) be the set of vertices outside the layers. On \(E_N\),
every permutation contributing a nonzero term to \(\det(zI-X_N)\)
fixes all \(Ld\) layer vertices. Indeed, any layer vertex moved by
such a permutation would belong to a directed cycle in the support
of \(X_N\). Layer vertices carry no loops, so \((X_N)_{ii}=0\) there, and each
fixed layer vertex contributes the factor \(z\). Summing over the permutations of \(W\) gives
\[
 \det(zI-X_N)=z^{Ld}\det\bigl(zI_{|W|}-X_N[W,W]\bigr).
\]
Thus zero has algebraic multiplicity at least \(Ld\) on \(E_N\), and
\[
 \mu_N(\{0\})\ge\frac{Ld}{N}=1-\frac{kd^3+r}{N}.
\]
As \(kd^3+r=o(N)\) and \(\Pr(E_N)\to1\), the mass at zero tends to
one in probability. This proves \(\mu_N\Rightarrow\delta_0\)
in probability.

\section*{AI tool disclosure}
The authors developed the main ideas and led the writing of the
manuscript. They used GPT-6 Astra to assist with proof development,
some derivations and exposition. The authors have verified the
arguments and take full responsibility for the paper.

\printbibliography[title={References}]
\end{document}